\documentclass{amsart}
\usepackage{amsmath}
\usepackage{amssymb}

\def\texorpdfblock#1#2{%
\ifx\pdfoutput\undefined%
#1%
\else%
#2%
\fi%
}

\ifx\pdfoutput\undefined
  
\fi

\texorpdfblock{
\usepackage[dvips]{graphicx}
\usepackage{psfrag}
}{
\usepackage[pdftex]{graphicx}
\usepackage[pdftex,bookmarksnumbered,colorlinks,a4paper,hypertexnames=true,linkcolor=blue,plainpages=false,pdfpagelabels,citecolor=red,urlcolor=blue]{hyperref}
\pdfcompresslevel9 \pdfadjustspacing=1 \hypersetup{
    pdftitle={Menger's theorem for infinite graphs},
    pdfauthor={Ron Aharoni and Eli Berger},
    pdfsubject={},
    pdfkeywords={{Menger's theorem}},
    pdfstartview={FitH}
  }
}

\title{Menger's theorem for infinite graphs}
\author{Ron Aharoni}
\address{Department of Mathematics\\Technion, Haifa\\ Israel 32000}
\email[Ron Aharoni]{ra@tx.technion.ac.il}

\author{Eli Berger}
\address{Department of Mathematics\\Haifa University, Haifa\\ Israel}
\email[Eli Berger]{berger@math.haifa.ac.il}
\begin{document}
 \thanks{\noindent The research of the first author was
supported by grant no. 780-04 of the Israel Science Foundation, by
GIF grant no. 2006311, by the Technion's research promotion fund,
and by the Discont Bank chair.}
\thanks{The research of the second author was supported by the National
Science Foundation, under agreement No. DMS-0111298. Any opinions,
findings and conclusions or recommendations expressed in this
material are those of the authors and do not necessarily reflect
the view of the National Science Foundation.}

\maketitle

\newcommand{\pf}{{\it Proof.}~}
\newcommand{\hetz}{{^\curvearrowright}}
\newcommand{\quo}{/}
\newcommand{\extends}{{\succcurlyeq}}
\newcommand{\fextends}{{\vec{\extends}}}
\newcommand{\extended}{{\preccurlyeq}}
\newcommand{\fextended}{{\vec{\extended}}}

\newtheorem{claim}{Claim}

\newcommand{\cig}{{\mathcal I}(G)}

\newcommand{\ca}{\mathcal A}
\newcommand{\cb}{\mathcal B}
\newcommand{\cc}{\mathcal C}
\newcommand{\cd}{\mathcal D}
\newcommand{\ce}{\mathcal E}
\newcommand{\cf}{\mathcal F}
\newcommand{\cg}{\mathcal G}
\newcommand{\ci}{\mathcal I}
\newcommand{\cj}{\mathcal J}
\newcommand{\ck}{\mathcal K}
\newcommand{\cl}{\mathcal L}
\newcommand{\cm}{\mathcal M}
\newcommand{\cn}{\mathcal N}
\newcommand{\co}{\mathcal O}
\newcommand{\cp}{\mathcal P}
\newcommand{\cq}{\mathcal Q}
\newcommand{\cs}{\mathcal S}
\newcommand{\cgr}{\mathcal R}
\newcommand{\ct}{\mathcal T}
\newcommand{\cu}{\mathcal U}
\newcommand{\cv}{\mathcal V}
\newcommand{\cx}{\mathcal X}
\newcommand{\cy}{\mathcal Y}
\newcommand{\cz}{\mathcal Z}
\newcommand{\cw}{\mathcal W}
\newcommand{\tlp}{T_{\lambda^+}}
\newcommand{\cyta}{\mathcal Y \langle T_\alpha \rangle}

\newcommand{\enp}{\hfill \Box}
\newcommand{\tn}{\tilde{N}}
\newcommand{\ch}{\mathcal H}

\newcommand{\etabar}{\overline{\eta}}
\newcommand{\heta}{\hat{\eta}}
\newcommand{\tg}{\tilde{\gamma}}
\newcommand{\ttau}{\tilde{\tau}}
\newcommand{\rn}{\mathbb{R}^n}

\theoremstyle{plain}
\newtheorem{theorem}{Theorem}[section]
\newtheorem{lemma}[theorem]{Lemma}
\newtheorem{fact}[theorem]{Fact}
\newtheorem{conjecture}[theorem]{Conjecture}
\newtheorem{corollary}[theorem]{Corollary}
\newtheorem{assertion}[theorem]{Assertion}
\newtheorem{proposition}[theorem]{Proposition}
\newtheorem{observation}[theorem]{Observation}
\newtheorem{problem}[theorem]{Problem}

\theoremstyle{remark}
\newtheorem{notation}[theorem]{Notation}
\newtheorem{remark}[theorem]{Remark}
\newtheorem{definition}[theorem]{Definition}
\newtheorem{convention}[theorem]{Convention}
\newtheorem{assumption}[theorem]{Assumption}

\theoremstyle{remark}
   \newtheorem{example}[theorem]{Example}

\begin{abstract}
We prove that Menger's theorem is valid for infinite graphs, in
the following strong version:  let $A$ and $B$ be two sets of
vertices in a possibly infinite digraph. Then there exist a set
$\cp$ of disjoint $A$--$B$ paths, and a set $S$ of vertices
separating $A$ from $B$, such that $S$ consists of a choice of
precisely one vertex from each path in $\cp$. This settles an old
conjecture of Erd\H{o}s.

\end{abstract}

\section{History of the problem}
In 1931 D\'enes K\"{o}nig \cite{konig31} proved a min-max duality
theorem on bipartite graphs:

\begin{theorem}\label{finitekonig}
In any finite bipartite graph, the maximal size of a matching
equals the minimal size of a cover of the edges by vertices.
\end{theorem}

Here a {\em matching} in a graph is a set of disjoint edges, and a
{\em cover} (of the edges by vertices) is a set of vertices
meeting all edges. This theorem was the culmination of a long
development, starting with a paper of Frobenius in 1912. For
details on the intriguing history of this theorem, see
\cite{lovaszplummer}. Four years after the publication of K\"onig's paper
 Phillip Hall
\cite{hall} proved a result which he named ``the marriage
theorem''. To formulate it, we need the following notation: given
a set $A$ of vertices in a graph, we denote by $N(A)$ the set of
its neighbors.

\begin{theorem}
In a finite bipartite graph with sides $M$ and $W$ there exists a
marriage of $M$ (that is, a matching meeting all vertices of $M$)
if and only if $|N(A)| \ge |A|$ for every subset $A$ of $M$.
\end{theorem}

The two theorems are closely related, in the sense that they are
 easily derivable from each other. In fact, K\"{o}nig's
theorem is somewhat stronger, in that the derivation of Hall's
theorem from it is more straightforward than vice versa.

At the time of publication of K\"{o}nig's theorem, a theorem
generalizing it considerably was already known.

\begin{definition}\label{separation} Let $X,Y$ be two sets of vertices in a digraph $D$. A set $S$ of
vertices is called $X$--$Y$-{\em separating} if every
 $X$--$Y$-path meets $S$, namely if the deletion of $S$ severs all
 $X$--$Y$-paths.
\end{definition}

Note that, in particular, $S$ must contain $X \cap Y$.

\begin{notation} The minimal size of an $X$--$Y$-separating set
is denoted by $\sigma(X,Y)$. The  maximal size of a family of
vertex-disjoint paths from $X$ to $Y$ is denoted by $\nu(X,Y)$.
\end{notation}

In 1927 Karl Menger \cite{menger} published the following:

\begin{theorem}
For any two sets $A$ and $B$ in a finite digraph there holds:
$$\sigma(A,B)=\nu(A,B) ~ .$$
\end{theorem}

This was probably the first casting of a combinatorial result in
min-max form. There was a gap in Menger's proof: he assumed,
without proof, the bipartite case of the theorem, which is Theorem
\ref{finitekonig}. This gap was filled by K\"{o}nig. Since then
other ways of deriving Menger's theorem from K\"onig's theorem
have been found, see,  e.g., \cite{aharonifinitepaths}.

Soon thereafter Erd\H{o}s, who was K\"onig's student, proved that,
with the very same formulation, the theorem is also valid for
infinite graphs. This appeared in K\"{o}nig's book
\cite{konigbook}, the first book published on graph theory. The
idea of the proof is this: take a maximal family $\cp$ of
$A$--$B$-disjoint paths. The set $S=\bigcup\{V(P):~P \in \cp\}$ is
then  $A$--$B$-separating, since an
 $A$--$B$-path avoiding it could be added to $\cp$, contradicting the
maximality of $\cp$. Since every path in $\cp$ is finite, if $\cp$
is infinite then $|\cp|=|S|$. Since $\nu(A,B) \ge |\cp|$ and
$\sigma(A,B) \le |S|$, this implies the non-trivial inequality
$\nu(A,B) \ge \sigma(A,B)$ of the theorem. If $\cp$ is finite,
then so is $S$. The size of families of disjoint $A-B$ paths is
thus finitely bounded (in fact, bounded by $|S|$), and hence there
exists a finite  family of maximal cardinality of disjoint
$A$--$B$ paths. In this case one can apply one of many proofs
known for the finite case of the theorem (see, e.g., Theorem
\ref{blockingsimply} below, or \cite{diestel}).

Of course, there is some ``cheating" here. The separating set
produced in the case that $\cp$ is infinite is obviously too
``large". In the finite case the fact that $|S|=|\cp|$ implies
that there is just one  $S$-vertex on each path of $\cp$, while in
the infinite case the equality of cardinalities does not imply
this.  Erd\H{o}s conjectured that, in fact, the same relationship
between $S$ and $\cp$ can be obtained also in the infinite case.
Since it is now proved, we state it as a theorem:

\begin{theorem}\label{erdos}\label{main}
Given two sets of vertices, $A$ and $B$, in a (possibly infinite)
digraph, there exists a family $\cp$ of disjoint  $A$--$B$-paths,
and a separating set consisting of the choice of precisely one
vertex from each path in $\cp$.
\end{theorem}

\begin{figure}[htb]
\begin{center}
\texorpdfblock{
\psfrag{A}{$A$}
\psfrag{B}{$B$}
\psfrag{S}{$S$}
\psfrag{F}{$\mathcal{P}$}
\includegraphics{menger-fig1.eps}
}{
\includegraphics{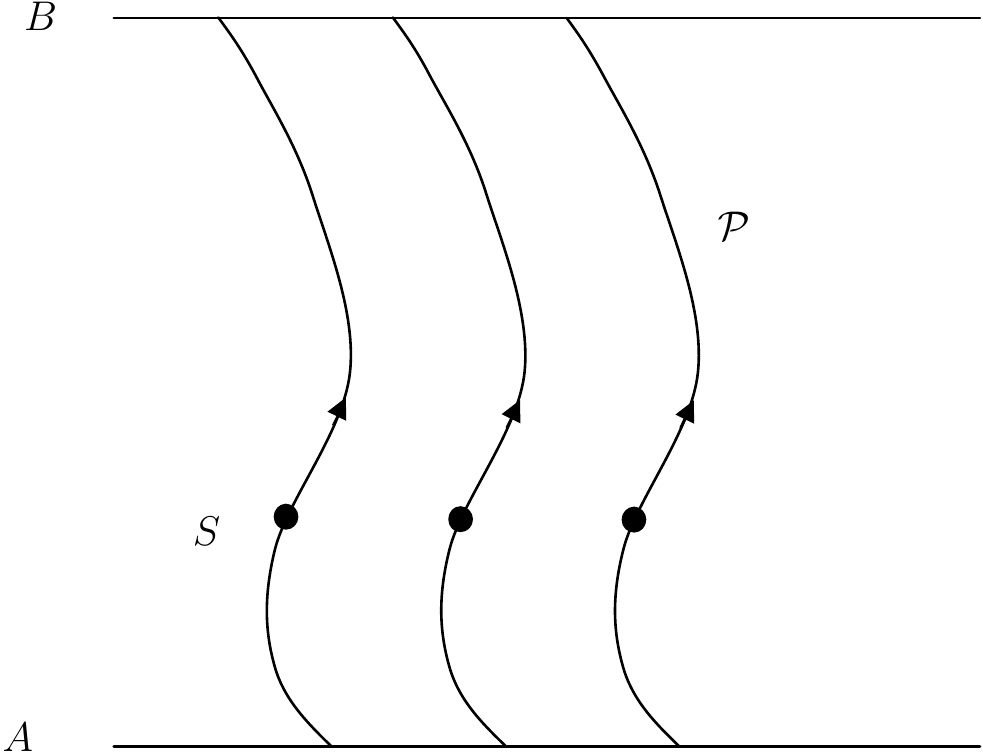}
}
\end{center}
\caption{Illustration of Theorem \ref{main}}
\end{figure}

The earliest reference in writing to this conjecture is
\cite{smolenice} (Problem 8, p. 159. See also \cite{nwproblems}).

The first to be tackled was of course the bipartite case, and the
first breakthrough was made by Podewski and Steffens
\cite{podsteff74}, who proved the countable bipartite case of the
conjecture, namely the countable case of K\"{o}nig's theorem. That
paper established some of the basic concepts that were used in
later work on the conjecture, and also set the basic approach:
introducing an a-symmetry into the problem. In the conjecture (now
theorem) the roles of $A$ and $B$ are symmetrical; the proof in
\cite{podsteff74} starts with asking the question of when can a
given side of a bipartite graph be matched into the other side,
namely the problem of extending Hall's theorem to the infinite
case. Known as the ``marriage problem", this question was open
since the publication of Hall's paper, and Podewski and Steffens
solved its countable case. Around the same time, Nash-Williams
formulated two other necessary criteria for matchability (the
existence of marriage), and he \cite{nw75, nw78} and Damerell and
Milner \cite{dm} proved their sufficiency for countable bipartite
graphs. These criteria are more explicit, but in hindsight the
concepts used in \cite{podsteff74} are more fruitful.

Podewski and Steffens \cite{podsteffmenger} made yet another
important progress: they proved the conjecture for countable
digraphs containing no infinite paths. Later, in
\cite{aharonifinitepaths}, it was realized that this case can be
easily reduced to the bipartite case, by the familiar device of
doubling vertices in the digraph, thus transforming the digraph
into a bipartite graph.

At that point  in time there were two obstacles on the way to the
proof of the conjecture - uncountability and the existence of
infinite paths. The first of the two  to be overcome was that of
uncountability. In 1983 the marriage problem was solved for
general cardinalities, in \cite{aharoninwshelah}. Soon thereafter,
this was used to prove the infinite version of K\"{o}nig's theorem
\cite{aharonikonig}. Namely, the bipartite case of Theorem
\ref{erdos} was proved. Let us state it explicitly:

\begin{theorem}\label{infkonig}
In any bipartite graph there exists a matching $F$ and a cover
$C$, such that $C$ consists of the choice of precisely one vertex
from each edge in $F$.
\end{theorem}

As is well known, Hall's theorem fails in the infinite case. The
standard example is that of the ``playboy'': take a graph with
sides $M=\{m_0, m_1, m_2, \ldots\}$ and $W=\{w_1, w_2, \ldots\}$.
For every $i>0$ connect $m_i$ to $w_i$, and connect $m_0$ (the
playboy) to all $w_i$. Then every subset of $M$ is connected to at
least as many points in $W$ as its size, and yet there is no
marriage of $M$. This is just another indication that in the case
of infinite matchings, cardinality is too crude
 a measure.

But Theorem \ref{infkonig} has an interesting corollary: that if
``cardinality'' is interpreted {\em in terms of the graph}, then
Hall's theorem does apply also in the infinite case. Given two
sets, $I$ and $J$, of vertices in a graph  $G$, we say that $I$ is
{\em matchable into $J$}  if there exists an injection of $I$ into
$J$ using edges of $G$. We write $I <_G J$ if $I$ is matchable
into $J$, but $J$ is not matchable into $I$. (The ordinary notion
of $|I|<|J|$ is obtained when $G$ is the complete graph on a
vertex set containing $I \cup J$.) A {\em marriage} of a side of a
bipartite graph is a matching covering all its vertices. From
Theorem \ref{infkonig} there follows:

\begin{theorem}
\label{infmarriage}
 Given a bipartite graph $\Gamma$ with sides $M$ and $W$,
there does {\em not}  exist a marriage of $M$ if and only if there
exists $A \subseteq M$, such that $N(A) <_\Gamma A$.
\end{theorem}

To see how Theorem \ref{infmarriage} follows from Theorem
\ref{infkonig}, assume that there is no marriage of $M$, and let
$F$ and $C$ be as in Theorem \ref{infkonig}. Let $I=M \setminus
C$. Then the set of points connected to $I$ is obviously $F[I]$
(the set of points connected by $F$ to $I$), which is matchable by
$F$ into $I$. If there existed a matching $K$ of $I$, then $K \cup
(F\upharpoonright (M \cap C))$ would be a marriage of $M$,
contrary to assumption. Thus $I$ is unmatchable. The other
implication in the theorem is obvious.

Proof-wise, the order is in fact reverse: Theorem
\ref{infmarriage} is proved first, and from it Theorem
\ref{infkonig} follows, in a way that will be explained later, in
Section 5.

By the result of \cite{aharonifinitepaths}, there follows  from
Theorem \ref{infkonig} also Theorem \ref{erdos} for all graphs
containing no infinite (unending or non-starting) paths. Thus
there remained the problem of infinite paths. The difficulty they
pose is that when one tries to ``grow'' the disjoint paths desired
in the conjecture, they may end up being infinite, instead of
being  $A$--$B$-paths. In fact, in \cite{aharonifinitepaths} it is
proved that Theorem \ref{erdos} is true, if one allows in $\cp$
not only
 $A$--$B$-paths, but any paths that if they start at all, they do so
at $A$, and if they end they do so at $B$.

The first breakthrough in the struggle against infinite paths was
made in \cite{countablemenger}, where the countable case of the
conjecture was proved.  An equivalent, Hall-type, conjecture, was
formulated, and the latter was proved for countable digraphs. The
core of the proof was in a lemma, stating that if the Hall-like
condition is satisfied, then any point in $A$ can be linked to $B$
by a path, whose removal leaves the Hall-like condition intact.
The lemma is quite easy to prove in the bipartite case and also in
graphs containing no unending paths, but in the general countable
case it requires new tools and methods. Later, the sufficiency of
the Hall-like condition for linkability (linking $A$ into $B$ by
disjoint paths) was proved for graphs in which all but countably
many points of $A$ are linked to $B$ \cite{countablelike}, and
Theorem \ref{erdos} was proved for such graphs in
\cite{aharonidiestel}.

In \cite{aharoni97}  a reduction was shown of the $\aleph_1$ case
of the conjecture to the above mentioned lemma. Namely, a proof of
the conjecture was given for digraphs of size $\aleph_1$, assuming
that the lemma is true for such digraphs.  Combined with a proof
of the lemma for graphs with no unending paths, and for graphs
with countable outdegrees, this settled the conjecture for
digraphs of size at most $\aleph_1$, satisfying one of those
properties. Optimistically, \cite{aharoni97} declares that this
reduction should probably work for general graphs.

The breakthrough leading to the solution of the general case was
indeed the proof of this lemma for general graphs. As claimed in
\cite{aharoni97}, the way from the lemma to the  proof of the
theorem indeed follows the same outline as in the $\aleph_1$ case.
But the general case demands quite a bit more effort.

For the sake of relative self containment of the paper, most
results from previous papers will be re-proved.

\section{Notation}

\subsection{Graph-theoretic notation}
 One non-standard
notation that we shall use is this: for a directed edge $e=(x,y)$ in
a digraph we write $x=tail(e)$ and $y=head(e)$. The rest of the
notation is mostly standard, but here are a few reminders. Given a
digraph $D$ and a subset $X$ of $V(D)$ we write $D[X]$ for the graph
induced by $D$ on $X$. Given a set $U$ of vertices in an undirected
graph, we denote by $N(U)$ the set of neighbors of vertices of $U$.
In a digraph
 we write $N^+(U)$ (respectively $N^-(U)$) for the set of out-neighbors
 (respectively in-neighbors)
 of $U$. Adopting a common abuse of notation, when $U$ consists of a
single vertex $u$, we write $N(u),N^+(u),N^-(u)$ for $N(\{u\}),
N^+(\{u\}),N^-(\{u\})$, respectively. Similar abuse of notation will
apply also to other notions, without explicit mention.

\subsection{Webs}
A {\em web} $\Gamma$ is a triple $(D,A,B)$, where $D=D(\Gamma)$ is a
digraph, and $A=A(\Gamma),B=B(\Gamma)$ are subsets of
$V(D)=V(\Gamma)$. We usually write $V$ for $V(D)$ and $E$ for
$E(D)$. If the identity of a web is not specified, we shall tacitly
assume that the above notation - namely $\Gamma, D, A$ and $B$ -
applies to it.

\begin{assumption} Throughout the paper we shall assume that there are no edges
going out of $B$, or into $A$.
\end{assumption}


Given a digraph $D$, we write $\overleftarrow{D}$ for the graph
having the same vertex set as $D$, with all edges reversed. For a
web $\Gamma =(D,A,B)$ we denote by $\overleftarrow{\Gamma}$ the web
$(\overleftarrow{D},B,A)$.

\subsection{Paths}
Following customary definitions, unless otherwise stated, a ``path"
in this paper is assumed to be simple, i.e. not self intersecting.
{\em All paths $P$ considered in the paper are assumed (unless empty) to have an
initial vertex, denoted by $in(P)$}. If $P$ is finite then its
terminal vertex is denoted by $ter(P)$. The vertex set of a path $P$
is denoted by $V(P)$, and its edge set by $E(P)$. The (possibly
empty) path obtained by removing $in(P)$ and $ter(P)$ from $P$ is
denoted by $P^\circ$.

Given a path $P$, we write $\overleftarrow{P}$ for the path in
$\overleftarrow{D}$ obtained by traversing $P$ in reverse order.

Given two vertices $u,v$ on a path $P$, we write $u \le_P v$ (resp.
$u <_P v$) if $u$ precedes $v$ on $P$ (resp. $u$ precedes $v$ on $P$
and $u \neq v$).


Given a set $\cp$ of paths, we write $\cp^f$ for the set  of finite
paths  in $\cp$, and $\cp^\infty$ for the set of infinite paths in
$\cp$. We also write $V[\cp]=\bigcup\{V(P): ~P \in \cp\}$,~
$E[\cp]=\bigcup\{E(P): ~P \in \cp\}$,~ $in[\cp]=\{in(P):~P \in
\cp\}$,~and $ter[\cp]=\{ter(P):~P \in \cp^f\}$.

For a vertex $x$, we denote by $(x)$ the path whose vertex set is
$\{x\}$, having no edges.

For $X,Y \subseteq V$, a finite path $P$ is said to be an
$X$--$Y$-{\em path} if $in(P) \in X$ and $ter(P) \in Y$.

 Given a path $P$ and a vertex $v \in V(P)$, we write $Pv$ for
the part of $P$ up to and including $v$, and $vP$ for the part of
$P$ from $v$ (including $v$) and on. If $Q=Pv$ for some $v \in V(P)$
we say that $P$ is a {\em forward extension} of $Q$ and write $P
\fextends Q$.

Given two paths, $P$ and $Q$ with $ter(P)=in(Q)$,
we write $P*Q$, or sometimes just $PQ$,
for the concatenation of $P$ and $Q$, namely the (not necessarily
simple) path, whose vertex set is $V(P)\cup V(Q)$ and whose edge set
is $E(P) \cup E(Q)$. If  $V(P) \cap V(Q) = \{ter(P)\}=\{in(Q)\}$
then $P*Q$ is a simple path. In this case clearly
 $P*Q \fextends P$. Given paths
$P,Q$ sharing a common vertex $v$, we write $PvQ$ for the (not
necessarily simple) path  $Pv*vQ$.

\subsection{Warps}
 A set of vertex disjoint paths is called a {\em
warp} (a term taken from weaving). If all paths in a warp are
finite, then we say that the warp is of {\em finite character}
(f.c.). A warp $\cw$ is called $X$-{\em starting} if $in[\cw]
\subseteq X$. Given two sets of vertices, $X$ and $Y$, a warp $\cw$
is called an $X$--$Y$-warp if for every $P \in \cw$ we have $in(P)
\in X, ~ter(P) \in Y$ and $V(P) \cap (X \cup Y) =\{in(P),ter(P)\}$.
We say that a warp $\cw$ {\em links} $X$ to $Y$ if for every $x \in
X$ there exists some $P \in \cw$ such that $V(P) \cap X = \{x\}$ and
$V(xP) \cap Y \neq \emptyset$. Note that a warp linking $X$ to $Y$
needs not be an $X$--$Y$ warp, namely the initial points of its
paths need not lie in $X$, and the terminal points do not
necessarily lie in $Y$. An $X$--$Y$-warp linking $X$ to $Y$ is
called an {\em $X$--$Y$- linkage}. An $A$--$B$-linkage in a web
$\Gamma=(D,A,B)$ is called a {\em linkage of} $\Gamma$. A web having
a linkage is called {\em linkable}. We write $\overleftarrow{\cw}$
for the warp $\{\overleftarrow{P} \mid P \in \cw\}$  in
$\overleftarrow{D}$.

For a set $X \subseteq V$, we denote by $\langle X \rangle$  the
warp consisting of all vertices of $X$ as singleton paths. For every
warp $\cw$ we write $ISO(\cw)$ (standing for ``isolated vertices of
$\cw$") for the set of vertices appearing in $\cw$ as singleton
paths.

\begin{notation}
 Given a warp $\cw$ and a set of vertices $X$,
we write $\cw[X]$ for the unique warp whose vertex set is $X \cap
V[\cw]$ and whose edge set is $\{(u,v) \in E[\cw] \mid ~ u,v \in
X\}$. Paths in $\cw[X]$ are sub-paths of paths in $\cw$. Note that a
path in $\cw$ may break into more than one path in $\cw[X]$. We also
write $\cw -X$ for $\cw[V\setminus X]$.
\end{notation}

\begin{definition}\label{extendingwarps}
A warp $\cu$ is said to be an {\em extension} of a warp $\cw$ if
$V[\cw] \subseteq V[\cu]$ and $E[\cw] \subseteq E[\cu]$. We write
then $\cw \extended \cu$. Note that $\cu$ may amalgamate paths in
$\cw$. If in addition $in[\cw] = in[\cu]$ then we say that $\cu$ is
a {\em forward extension} of $\cw$ and write $\cu \fextends \cw$.
Note that in this case each path in $\cu$ is a forward extension of
some path in $\cw$.
\end{definition}

\begin{notation}\label{warpx}
Given a warp $\cw$ and a set $X \subseteq V$, we write $\cw\langle X
\rangle$ for the set of paths in $\cw$ intersecting $X$, and
$\cw\langle \sim X \rangle$ for $\cw \setminus \cw\langle X \rangle$. Given two sets of vertices, $X$ and $Y$, we write
$\cw\langle  X,Y  \rangle$ for $\cw\langle X \rangle \cap \cw\langle
Y \rangle$ and $\cw\langle  X, \sim Y  \rangle$ for $\cw\langle X
\rangle \cap \cw\langle \sim Y \rangle$.

Given a vertex $x \in V[\cw]$ we write $\cw (x)$ for the path in
$\cw$ containing $x$ (to be distinguished from $\cw \langle
x \rangle$,  the   set
consisting of the single path $\cw (x)$).
\end{notation}

Given a warp $\cw$ in a web $(D,A,B)$, we write $\cw_G$ for $\cw
\langle A \rangle$ and $\cw_H$ for $\cw \setminus \cw_G$ (the
subscript ``G" stands for ``ground" - these are the paths in $\cw$
that start ``from the ground", namely at $A$. The subscript ``H"
stands  for ``hanging in air". These terms  originate in the way the
authors are accustomed to draw webs - with the ``$A$" side at the
bottom, and the ``$B$" side on top).
\\

A set $\cf$ of paths is called a {\em fractured warp} if its edge
set is the edge set of a warp and every two paths $P,Q \in \cf$ may
intersect only if none of them is a trivial path and $in(P) =
ter(Q)$ or $in(Q) = ter(P)$. If $\cw$ is a warp and $X$ is a set of
vertices, we write $\cw \downharpoonright X$ for the fractured warp
consisting of all paths of the form $xPy$ where $P \in \cw$, $x \in
X \cup \{in(P)\}$, $y \in X \cup \{ter(P)\}$, $V(xPy) \not \subseteq
X$ and $V(xPy) \cap X \subseteq \{x,y\}$. A somewhat more comprehensible definition
is given by the following properties: $E[\cw
\downharpoonright X] = E[\cw] \setminus E[\cw[X]]$,
no path in $\cw
\downharpoonright X$ amalgamates two paths in $\cw$, and all singleton paths $(y) \in \cw$,
where $y \not \in X$, belong to $\cw
\downharpoonright X$.

A set of pairwise disjoint paths and directed cycles is called a
{\em cyclowarp}. If $\cc$ is a cyclowarp, we denote by
$\cc^{path}$ the warp obtained from $\cc$ by removing all
its cycles.

\subsection{Operations between warps} \nopagebreak [4] \begin{notation}\label{star}
\nopagebreak [4] Let $\cu$ and $\cw$ be warps such that $V[\cu] \cap
V[\cw] \subseteq ter[\cu] \cap in[\cw]$. Denote then by $\cu*\cw$
the warp $\{P*Q \mid ~P \in \cu^f,~Q \in \cw,~in(Q)=ter(P)\} \cup
\{P \in \cu \mid ter(P) \not \in in[\cw]\}$. In particular,
$\cu*\cw$ contains $\cu^\infty$. Denote by $\cu \diamond \cw$  the
warp whose vertex set is $V[\cu]\cup V[\cw]$ and whose edge set is
$E[\cu]\cup E[\cw]$.
\end{notation}

Thus $\cu \diamond \cw \supseteq \cu*\cw$. The difference is that
$\cu \diamond \cw$ may contain also paths in $\cw$ not meeting any
path from $\cu$.

\begin{figure}[htb]
\begin{center}
\texorpdfblock{
\psfrag{A}{$A$}
\includegraphics[scale=0.8]{menger-fig2.eps}\hspace*{3em}
\includegraphics[scale=0.8]{menger-fig3.eps}
}{
\includegraphics[scale=0.8]{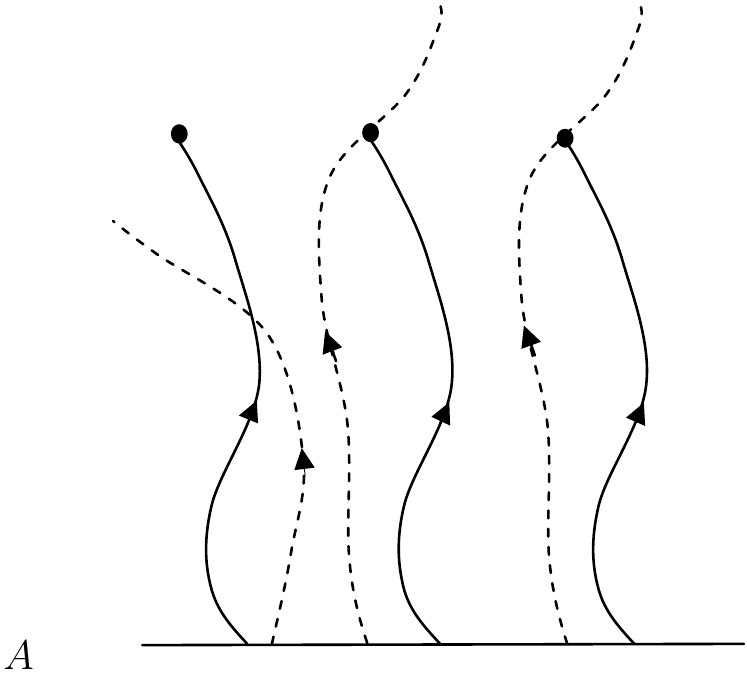}\hspace*{3em}
\includegraphics[scale=0.8]{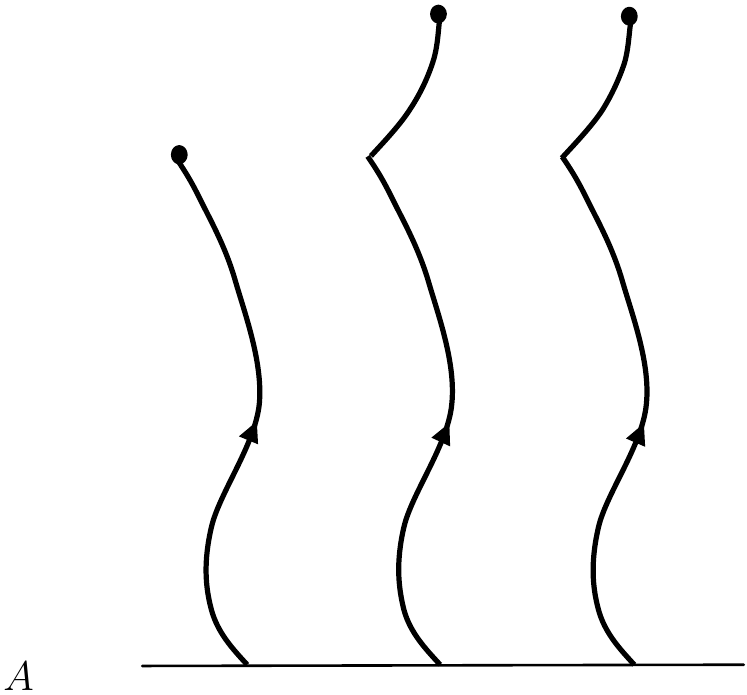}
}
\end{center}
\caption{On the left there are drawn a warp $\cu$ (solid line) and a
warp $\cw$ (dashed line). On the right is drawn their ``arrow", $\cu
\hetz \cw$.} \label{figurehetz}
\end{figure}

There is also a binary operation defined on {\em all} pairs of
warps. Given warps $\cu$ and $\cw$, their ``arrow" $\cu \hetz \cw$
is obtained by taking each path in $\cu$ and ``carrying it along
$\cw$", if possible, and if not keeping it as it is. Formally, this
is defined as follows:

\begin{notation}\label{hetz}
Let $\cu$ and $\cw$ be two warps and let $P$ be a path in $\cu$. We
define the $\cu$-$\cw$-{\em extension} $Ext_{\cu-\cw}(P)$ of $P$ as
follows. Consider first the case that $P$ is finite. Let $u =
ter(P)$. If there exists a path $Q \in \cw$ satisfying $u \in V(Q)$
and $V(uQ) \cap V[\cu] = \{u\}$ let $Ext_{\cu-\cw}(P) = PuQ$. In any
other case (i.e. if either $P$ is infinite or $u \not \in V[\cw]$ or
$V(u\cw(u))$ meets $\cu$ at a vertex other than $u$) we take
$Ext_{\cu-\cw}(P) = P$. Let
$$\cu \hetz \cu = \{Ext_{\cu-\cw}(P) : ~ P \in \cu \}.$$

(See Figure \ref{figurehetz}.)
\end{notation}

Note that $\cu \hetz \cw$ is a warp and $\cu \hetz \cw \fextends
\cu$.
\begin{observation}\label{fextends_iff_hetsisitself}
$\cw \fextends \cu$ if and only if $\cu \hetz \cw=\cw$.
\end{observation}

Next we wish to define the ``arrow" of a sequence of warps. As a
first step, we define the limit of an ordinal-indexed sequence of
warps.

\begin{definition}\label{liminfofwarps}
Let $(S_\alpha:~ \alpha < \theta)$ be a sequence of sets.
The {\em limit} (actually, $\lim\inf$)
of the sequence, denoted by  $\lim_{\alpha<\theta}S_\alpha$, is defined as
$\bigcup_{\beta < \theta} \bigcap_{\beta \le \alpha < \theta}
S_\alpha$. Let $(\cw_\alpha:~ \alpha < \theta)$ be a sequence of
warps. The {\em limit} $\lim_{\alpha<\theta}\cw_\alpha$ of the
sequence is the warp whose edge set is $\lim_{\alpha < \theta}
E[\cw_\alpha]$ and whose vertex set is $\lim_{\alpha < \theta}
V[\cw_\alpha]$.
\end{definition}

As noted, $\lim_{\alpha<\theta}\cw_\alpha$ is in fact the ``lim inf"
of the warps. The fact that it is indeed a warp is straightforward.
Note that by this definition if $\theta$ is not a limit ordinal,
namely $\theta=\psi+1$, then $\lim_{\alpha<\theta}\cw_\alpha$ is
just $\cw_\psi$.

\begin{observation}\label{limter}
Let $(\cw_\alpha:~ \alpha < \theta)$ be a sequence of warps. Then
$ter[\lim_{\alpha<\theta}\cw_\alpha] \supseteq
\lim_{\alpha<\theta}ter[\cw_\alpha]$.
\end{observation}




\begin{definition}\label{uparrow}
Let $(\cw_\alpha:~ \alpha < \theta)$ be an ordinal-indexed sequence
of warps. Define a sequence $\cw'_\alpha, ~\alpha < \theta$, by:
$\cw'_0=\cw_0$, $\cw'_{\psi+1}=\cw'_{\psi}\hetz \cw_{\psi+1}$ (where
$\psi+1 <\theta$), and for limit ordinals $\alpha \le \theta$ define
$\cw'_\alpha= \lim_{\psi < \alpha}\cw'_\psi$. Let $\uparrow_{\alpha
< \theta} \cw_{\alpha}$ be defined as $\cw'_\theta$ if $\theta$ is a
limit ordinal, and as $\cw'_\beta$ if $\theta=\beta+1$.
\end{definition}

Note that if $(\cw_\alpha:~ \alpha < \theta)$ is
$\fextended$-ascending, then this definition coincides with the
``limit" definition. If $\{\cw_i, ~i \in I\}$ is an unordered  set
of warps, then  $\uparrow_{i \in I} \cw_i$ can be defined by first imposing
an arbitrary well-order on $I$. Of course, the resulting warp
depends on the order chosen, but when applied we shall use a fixed well order.




\subsection{Almost disjoint families of paths}

Given a set $X$ of vertices, a set $\cp$ of paths is called $X$-{\em
joined} if the intersection of the vertex sets of any two paths from
$\cp$ is contained in $X$ (so, a warp is just a $\emptyset$-joined
family of paths). For a single vertex $x$, we write simply
``$x$-joined" instead of ``$\{x\}$-joined". A family of $x$-joined
paths starting at $x$ is called a {\em fan}. A family of $x$-joined
paths terminating at $x$ is called an {\em in-fan}.

Given a set $X \subseteq V$ and a vertex $u$, a $u$-fan  $\cf$ is said to be a
$u$--$X$-fan  if  $ter[\cf] \subseteq X$.
An $X$--$u$-in-fan is defined similarly.
  A $u$-fan consisting of
infinite paths is called a $(u,\infty)$-fan.

\subsection{Separation}
\begin{definition}\label{separating}
An $A$--$B$-{\em separating} set of vertices in a web
$\Gamma=(D,A,B)$ is  plainly said to be {\em separating}.
\end{definition}

\begin{definition}\label{essentialweb}
 Given a (not necessarily separating)
subset $S$ of $V(D)$, a vertex $s \in S$ is said to be {\em
essential} (for separation) in $S$ if it is not separated from $B$
by $S \setminus \{s\}$. The set of essential elements of $S$ is
denoted by $\ce(S)$, and the set $S \setminus \ce(S)$ of inessential
vertices by $\ci\ce(S)$. If $S = \ce(S)$ then we say that $S$ is
{\em trimmed}.
\end{definition}

\begin{convention} By removing those vertices of $A$ from which $B$ is unreachable, we
may assume that $A$ is trimmed. We shall tacitly make this
assumption.
\end{convention}

\begin{lemma}
\label{essentialpoints} If $S$ is an $A$--$B$ separating set of
vertices, then so is $\ce(S)$.
\end{lemma}

\begin{proof}
Let $Q$ be an  $A$--$B$-path. Since by assumption $S$ is $A$--$B$
separating, $V(Q) \cap S \neq \emptyset$. The last vertex $s$ on $Q$
belonging to $S$ is essential in $S$, since the path $sQ$ shows that
$s$ is not separated from $B$ by $S \setminus \{s\}$.
\end{proof}

A path $P$ in a warp $\cw$ is said to be {\em essential} (in $\cw$)
if $P$ is finite and $ter(P) \in \ce(ter[\cw])$. The set of
essential paths in $\cw$ is denoted by $\ce(\cw)$, and the set of
inessential paths by $\ci\ce(\cw)$. If $\cw = \ce(\cw)$ we say that
$\cw$ is {\em trimmed}.

To Definition \ref{separation} we add the following. Given a set $X$
of vertices, a vertex set $S$ is called $X$-$\infty$-{\em
separating} if it contains a vertex on every infinite path starting
in $X$. The minimal size of an $X$-$\infty$-separating set is
denoted by $\sigma(X,\infty)$.




\begin{notation}\label{roofnotation}
For a set $S$ of vertices in a web $\Gamma=(D,A,B)$ we denote by
$RF(S)=RF_\Gamma(S)$ the set of all vertices separated by $S$ from
$B$. We also write $RF^\circ(S)=RF(S)\setminus \ce(S)$.
\end{notation}

The letters ``$RF$'' stand for ``roofed'', a term originating again
in the way the authors draw their webs, with the ``$A$" side at the
bottom, and the ``$B$" above. Note that in particular, $S \subseteq
RF(S)$ and $\ci\ce(S) \subseteq RF^\circ(S)$. Given a warp $\cw$, we
write $RF(\cw)=RF(ter[\cw])$, $RF^\circ(\cw)=RF^\circ(ter[\cw])$. A
warp $\cw$ is said to be {\em roofed} by a set of vertices $S$ if
$V[\cw] \subseteq RF(S)$.

\begin{lemma}\label{lastonpathbelongstoroof}
Let $S$ be a set of vertices and $P$ any path. If $V(P) \cap
RF(S)\neq \emptyset$ then the last vertex on $P$ belonging to
$RF(S)$ belongs to $\ce(S) \cup \{ter(P)\}$.
\end{lemma}
\begin{proof}
Let $v$ be the last vertex on $P$ belonging to $RF(S)$. Suppose that
$v \neq ter(P)$. We have to show that  $v \in \ce(S)$. Let $u$ be
the vertex following $v$ on $P$. Then $u \not \in RF(S)$, meaning
that  there exists an $S$-avoiding path $Q$ from $u$ to $B$. Since
$v \in RF(S)$ the path $vuQ$ meets $\ce(S)$. Since this meeting can
occur only at $v$, it follows that $v \in \ce(S)$.
\end{proof}

\begin{lemma}\label{essentialsandwiched}
If $C,D$ are sets of vertices such that $\ce(D) \subseteq C
\subseteq D$ then $\ce(C)=\ce(D)$.
\end{lemma}

\begin{proof} Let $x \in \ce(D)$. Then there exists an $x$--$B$ path
avoiding $D \setminus \{x\}$, and thus avoiding $C \setminus \{x\}$,
showing that $x \in \ce(C)$. On the other hand, if $x \in \ce(C)$
then there exists an $x$--$B$ path $P$ avoiding $C \setminus \{x\}$.
If $P$ does not avoid $D \setminus \{x\}$ then its last vertex
belongs to $\ce(D)$, and thus to $C$, a contradiction. Thus $x \in
\ce(D)$.
\end{proof}

\begin{observation}\label{STXY}
Let $S,T,X,Y$ be four sets of vertices, with $X \cap Y = \emptyset$.
If $X \subseteq RF(T \cup Y)$ and $Y \subseteq RF(S \cup X)$ then $X
\cup Y \subseteq RF(S \cup T)$ (otherwise stated as: $\ce(S \cup T
\cup X \cup Y) = \ce(S \cup T)$).
\end{observation}

\begin{proof} For an $(X \cup Y)$--$B$ path $P$ consider the last
vertex $z$ on $P$ belonging to $X \cup Y$. By the conditions of the
observation, $zP$ must meet $S \cup T$. \end{proof}

\begin{lemma}\label{sandwich}
If $R,S,T$ are three sets of vertices satisfying $T = \ce(T)$ and
$RF(R) \subseteq RF(S) \subseteq RF(T)$ then $S$ is
$R$--$T$-separating.
\end{lemma}

\begin{proof}
Consider an $R$--$T$ path $P$ and let $x = ter(P)$. Since $T =
\ce(T)$ there exists an $x$-$B$ path $Q$ satisfying $in(Q) = x$ and
$V(Q) \cap T = \{x\}$. Then $PxQ$ is an $R$--$B$ path. (A-priory, $PxQ$
might not be a simple path. However, it obviously {\em contains}
a simple $R$--$B$ path.)
Since $S$ is $R$-$B$ separating, we have $V(PxQ) \cap S \neq
\emptyset$. But since $S \subseteq RF(T)$ and $V(Q) \cap T = \{x\}$,
we have $V(Q) \cap RF(S) \subseteq \{x\}$, and hence $V(PxQ) \cap S
= V(P) \cap S \neq \emptyset$, proving the lemma.
\end{proof}

\begin{notation}\label{restriction_of_web_to_a_roofed_set}
Let $S$ be a  set of vertices in a web $\Gamma=(D,A,B)$, such that
$RF(S)=S$ (which is equivalent to $S$ being equal to $RF(T)$ for
some set $T$). We denote then by $\Gamma[S]$ the web $(D[S], S \cap
A, \ce(S))$. Given a warp $\cw$ we write $\Gamma[\cw]$ for
$\Gamma[RF(\cw)]$.
\end{notation}


\begin{lemma}\label{limroof}
Let $(S_\alpha:~ \alpha < \theta)$ be a sequence of sets, satisfying
$S_\alpha \subseteq RF(S_\beta)$ for $\alpha < \beta < \theta$. Then
$RF(\lim_{\alpha < \theta} S_\alpha) \supseteq \bigcup_{\alpha <
\theta} RF(S_\alpha)$.
\end{lemma}

\begin{proof}
Let $x \in \bigcup_{\alpha < \theta} RF(S_\alpha)$. We may assume that $x
\in RF(S_0)$ and thus $x \in \bigcap_{\alpha < \theta}
RF(S_\alpha)$. Let $P$ be an $x$--$B$ path and let $t$ be the last
vertex on $P$ belonging  to $\bigcup_{\alpha < \theta} S_\alpha$.
Say, $t \in S_\beta$. Since $t \in RF(S_\alpha)$ for all $\beta <
\alpha < \theta$ the vertex $t$ must be in $S_\alpha$ for all such
$\alpha$, and hence $t \in \lim_{\alpha < \theta} S_\alpha$.
\end{proof}

\subsection{Deletion and quotient}

A basic operation on webs is that of removing vertices. In fact,
there are two ways of doing this. One is plain deletion: for a
subset $X$ of $V$ we denote by   $\Gamma-X$ the web $(D-X, A
\setminus X,~B \setminus X)$. For a path $P$ we abbreviate and write
$\Gamma-P$ instead of $\Gamma-V(P)$.

\begin{lemma} \label{roofofunnion}

 $RF(X \cup Y) = X \cup RF_{\Gamma-X}(Y)$.

\end{lemma}

\begin{proof}
Note that $X$ is contained in the sets appearing on both sides of
the equality. Hence it suffices to show that for a vertex $v$ not
belonging to $X$, namely a vertex of $\Gamma-X$, a $v$--$B$ path
avoids $X \cup Y$ in $\Gamma$ if and only if it avoids $Y$ in
$\Gamma - X$. But this is almost a tautology.
\end{proof}

The other type of removal is
 taking a quotient. The difference from deletion is that
 taking a quotient with respect to a set $X$ of
 vertices means deleting the vertices of $X$ as vertices through which paths
 can go from $A$ to $B$, but also adding $X$ to $A$, indicating a commitment to link $X$ to $B$. If indeed
 such linking is possible, then the possibility arises of linking vertices of $A$ to $B$ by first linking them to $X$,
 and then linking vertices of $X$
 to $B$.

 \begin{definition}\label{gammaoverx}
 Given a
subset $X$ of $V \setminus A$, write $D\quo X$ for the digraph
obtained from $D$ by deleting all edges going into vertices of $X$,
and all vertices in $RF^\circ(X)$, including those of $\ci\ce(X)$.
Define $\Gamma \quo X$ as the web $(D \quo X, \ce(A \cup X), B)$.
\end{definition}

\begin{observation} Since we are assuming that $A$ is trimmed,
$A(\Gamma \quo X)=(A \cup X) \setminus RF^\circ(X)$.
\end{observation}

\begin{remark}\label{quotientinbipartite}
In bipartite webs deleting a vertex $b \in B$ and taking a quotient
with respect to it are the same, as far as linkability is concerned,
since taking a quotient with respect to $b$ means that $b$ is added
to $A$, and is linked automatically to itself. This is the reason
why the quotient operation is not needed in the proof of the
bipartite case of the theorem.
\end{remark}

\begin{lemma}\label{roofofunnionb}
For any two sets $X$ and $Y$ of vertices, $RF^\circ_\Gamma(X \cup
Y)=RF^\circ(X) \cup RF_{\Gamma \quo X}(Y \setminus RF^\circ(X))$.
\end{lemma}

\begin{proof}
Let $v$ be a vertex in $RF^\circ_\Gamma(X \cup Y)$. Suppose that $v
\not \in
 RF^\circ(X) \cup RF_{\Gamma \quo X}(Y
\setminus RF^\circ(X))$. Then there exists a path $P$ from $v$ to
$B$ in $\Gamma \quo X$,
 avoiding $(Y
\setminus RF^\circ(X))\setminus \{v\}$. Since   $P$ is contained in
$V(\Gamma \quo X)$, it is disjoint from $RF^\circ(X)$. Hence the
fact that it avoids $(Y \setminus RF^\circ(X))\setminus \{v\}$ means
that in fact it avoids  $Y \setminus \{v\}$, and since there are no
edges in $\Gamma \quo X$ going into $X$, it  also avoids $(X \cup Y)
\setminus \{v\}$, contradicting the assumption on $v$. This proves
that the set on the left hand side is contained in the set on the
right hand side. The other containment relation is proved similarly.
\end{proof}

\begin{lemma}\label{quotientbyunion}
For any two sets $X$ and $Y$ of vertices, if $Y \cap RF^\circ(X) = \emptyset$ then $\Gamma \quo (X \cup Y)=
(\Gamma \quo X) \quo Y$.
\end{lemma}

\begin{proof}
Let us first show that the two webs share the same vertex set. By
the definition of the quotient, we need to show: $$V
\setminus RF^\circ_\Gamma(X \cup Y)= ((V \setminus
RF^\circ_\Gamma(X)) \setminus
RF^\circ_{\Gamma \quo X}(Y).$$ This follows from
Lemma \ref{roofofunnionb}.

We also need to show that the two webs share the same source set,
namely that $\ce_\Gamma(A \cup X \cup Y)=\ce_{\Gamma \quo
X}(\ce_\Gamma(A \cup X) \cup Y)$. But this
also follows from Lemma \ref{roofofunnionb}.

It remains to show equality of  the edge sets of the two webs.
Write $V' = V(\Gamma \quo (X \cup Y)) = V((\Gamma \quo X) \quo Y)$
and $A' = A(\Gamma \quo (X \cup Y)) = A((\Gamma \quo X) \quo Y)$.
Then $E(\Gamma \quo (X \cup Y)) = E((\Gamma \quo X) \quo Y)
= \{(u,v) \in E(\Gamma) \mid ~u \in V', ~ v \in V' \setminus A'\}$.


\end{proof}

\begin{corollary}\label{quotientbyunioncor}
For any two sets $X_1$ and $X_2$ of vertices, if $Y = \ce(X_1 \cup X_2)$
then $(\Gamma \quo X_1) \quo Y = (\Gamma \quo X_2) \quo Y = \Gamma \quo Y$.
\end{corollary}

Given a warp $\cw$, we write $\Gamma \quo \cw$ for $\Gamma \quo
ter[\cw]$. If  $\cu$ and $\cw$ are two warps, we write $\cu \quo \cw$ for $\cu
\quo ter[\cw]$.


\begin{definition}
\label{quotientofwarps} Given a warp $\cw$ and a set $X$ of
vertices, we define the quotient $\cw \quo X$
by $V[\cw \quo X] = (V[\cw] \cup X) \setminus RF^\circ(X)$ and
$E[\cw \quo X] = \{(u,v) \in E[\cw] \mid~ u \not \in RF^\circ(X) ,v \not \in
RF(X)\}$.
\end{definition}

We end this section with a few lemmas, some of which are obvious and
some have similar proofs to those above, and hence we list them
without proofs:

\begin{lemma}
\label{quotientisawarpinquotent} $\cw \quo X$ is a warp in $\Gamma
\quo X$.
\end{lemma}

\begin{lemma}
$\langle \ce(X) \setminus V[\cw] \rangle \subseteq \cw \quo X$.
\end{lemma}

\begin{lemma}\label{quotientisastartinginquotient}
If $in[\cw] \subseteq A(\Gamma)$ then $in[\cw \quo X] \subseteq
A(\Gamma \quo X)$
\end{lemma}

\begin{lemma}
If $\cw \extended \cw'$ then $\cw \quo X \extended \cw' \quo X$. If
$\cw \fextended \cw'$ then $\cw \quo X \fextended \cw' \quo X$.
\end{lemma}

\begin{lemma}\label{terminalofquotientofwarp}
$in[\cw \quo X] = (in[\cw] \cup X) \setminus RF^\circ(X)$ and
$ter[\cw \quo X] \supseteq (ter[\cw] \setminus RF^\circ(X)) \cup (\ce(X)
\setminus V[\cw])$.
\end{lemma}


\begin{lemma}\label{missinglemmaonroofing}
For a subset $Z$  of $V(\Gamma)$ and a warp $\cv$ in $\Gamma$ we
have $RF^\circ_{\Gamma}(\cv) \cap V(\Gamma \quo Z)\subseteq
RF_{\Gamma \quo Z}^\circ(\cv \quo Z)$.
\end{lemma}

\begin{lemma}\label{quotientstrongerinseparation}
If $S,T$ are disjoint sets of vertices, then $RF_{\Gamma-T}(S)
\setminus RF^\circ(T) \subseteq RF_{\Gamma \quo T}(S \setminus
RF^\circ(T))$.
\end{lemma}

\section{Waves and hindrances}

\begin{definition}
 An  $A$-starting warp $\cw$ is called a {\em wave} if  $ter[\cw]$ is  $A$--$B$-separating.
  \end{definition}

 Clearly,
$\langle A \rangle$ (namely, the set of singleton paths, $\{(a) \mid
~a \in A\}$), is a wave. It is called {\em the trivial wave}.


\begin{observation}
\label{waveinduced} If $S = RF(S) \supseteq A$ and $\cw$ is a wave
in $\Gamma[S]$ then $\cw$ is also a wave in $\Gamma$.
\end{observation}

Lemma \ref{essentialpoints} implies:

\begin{lemma}\label{above}
  If
$\cw$ is a wave then so is $\ce(\cw)$.
\end{lemma}

This gives

\begin{lemma}
\label{essential} A path $W$ belonging to a wave $\cw$ is essential
in $\cw$ if and only if $\cw \setminus \{W\}$ is not a wave.
\end{lemma}

\begin{proof}
If $W$ is inessential, then by Lemma \ref{above}, $A \subseteq
RF(\ce(\cw)) \subseteq RF(\cw \setminus \{W\})$.

If, on the other hand, $W$ is essential, then $W$ is finite. Let
$t=ter(W)$. Since $t \in \ce(ter[\cw])$, there exists a path $P$
from $t$ to $B$ avoiding $ter[\cw]\setminus \{t\}$, and then $WtP$
is an $A$--$B$ path avoiding $ter[\cw \setminus \{W\}]$, showing
that $\cw \setminus \{W\}$ is not a wave.
\end{proof}

One nice property of waves is that they stay waves upon taking quotients.

\begin{lemma}\label{waveinquotient} If $\cu$ is a wave and $X
\subseteq V$ then $\cu\quo X$ is a wave in $\Gamma\quo X$.
\end{lemma}

\begin{proof}
By Lemmas \ref{quotientisawarpinquotent} and \ref{quotientisastartinginquotient},
The warp $\cu \quo X$ is indeed an $A(\Gamma \quo X)$-starting warp
in $\Gamma \quo X$.

Let $Q$ be a path in $\Gamma \quo X$ from $A(\Gamma \quo X)$, namely
$(A \cup X) \setminus RF^\circ(X)$, to $B$. We have to show that $Q$
meets $ter[\cu \quo X]$.

If $in(Q) \in A$ then, since $\cu$ is a wave, $in(Q) \in RF_\Gamma[\cu]$.
Otherwise $in(Q) \in \ce(X)$. Thus in both cases $in(Q) \in RF_\Gamma[\cu] \cup \ce(X)$.
Let $t$ be the last vertex on $Q$ belonging to $RF_\Gamma[\cu] \cup
\ce(X)$. From the choice of $t$ it follows that $t \not \in
RF_\Gamma^\circ(X) \cup RF_\Gamma^\circ(\cu)$, and hence $t \in (ter[\cu]
\setminus RF_\Gamma^\circ(X)) \cup (\ce(X) \setminus RF_\Gamma^\circ(\cu))$. By
Lemma \ref{terminalofquotientofwarp}  $t \in ter[\cu \quo X]$.
\end{proof}

A wave $\cw$ is called a {\em hindrance} if $in[\cw] \neq A$. The
origin of the name is that in finite webs a hindrance is an
obstruction for linkability. In the infinite case this is not
necessarily so. A web containing a hindrance is said to be {\em
hindered}.

As a corollary of Lemma \ref{waveinquotient} we have

\begin{corollary}\label{hindranceinquotient}
If $A \cap RF(S) = \emptyset$ and $\ch$ is a hindrance in $\Gamma$
then $\ch \quo S$ is a hindrance in $\Gamma \quo S$.
\end{corollary}

For, if $a \in A \setminus in[\ch]$ then $a \in A(\Gamma \quo S)
\setminus in[\ch \quo S]$.

Clearly, a hindrance is a non-trivial wave. A web not containing any
non-trivial wave is called {\em loose}.


\begin{lemma}[the self roofing lemma]\label{separatingverticesofwave}
If $\cw$ is a wave then $V[\cw]\subseteq RF(\cw)$.
\end{lemma}

\begin{proof}
Suppose, for contradiction, that there exists a  path $Q$ avoiding
$ter[\cw]$, from some vertex $x$ on a path $P \in \cw$ to $B$.
Taking a sub-path of $Q$, if necessary, we can assume that $PxQ$ is
a path. Then $PxQ$  avoids $ter[\cw]$, contradicting the fact that
$\cw$ is a wave.
\end{proof}

\begin{corollary}\label{separatingverticesofwavedeletion}
Let $X \subseteq V$ and let $\cw$ be a wave in $\Gamma - X$. Then
$V[\cw] \setminus ter[\cw] \subseteq RF^\circ(ter[\cw] \cup X)$
\end{corollary}

\begin{proof}
Let $u \in V[\cw] \setminus ter[\cw]$. By Lemma
\ref{separatingverticesofwave} we have $V[\cw] \subseteq RF_{\Gamma
- X}(\cw) \subseteq RF_\Gamma(ter[\cw] \cup X)$. Since $u \not \in
ter[\cw] \cup X$, we get $u \in RF^\circ(ter[\cw] \cup X)$.
\end{proof}

\begin{definition}
A warp $\cw$ is called {\em self roofing} if $V[\cw] \subseteq
RF(\cw)$.
\end{definition}

Lemma \ref{separatingverticesofwave} implies that every wave is self
roofing. In fact, an easy corollary of this lemma together with Lemma
\ref{waveinquotient} extends it to
waves in quotient webs.

\begin{corollary}\label{quo_self_roofing}
If $\cw$ is a wave in $\Gamma \quo X$ for some set $X$ then $\cw$ is
a self roofing warp in $\Gamma$.
\end{corollary}



For two waves $\cw$ and $\cw'$ we write $\cw \equiv \cw'$ if
$ter[\ce(\cw)] = ter[\ce(\cw')]$. Also write $\cw \le \cu$ if
$RF(\cw) \subseteq RF(\cu)$. Clearly, this is equivalent to the
statement that $ter[\cw] \subseteq RF(\cu)$. The relation $\le$ is a
partial order on the equivalence classes of the $\equiv$ relation.
Namely, if $\cw \le \cu$ and $\cw \equiv \cw'$, $\cu \equiv \cu'$
then $\cw' \le \cu'$, while if $\cw \le \cu$ and $\cu \le \cw$ then
$\cu \equiv \cw$. We write $\cu > \cw$ if $\cw \le \cu$ and $\cw
\not \equiv \cu$, i.e., $RF(\cw) \subsetneqq RF(\cu)$. We say that a
wave $\cw$ is $\le$-{\em maximal} if there is no wave $\cu$
satisfying $\cu > \cw$.

By the self roofing lemma (Lemma \ref{separatingverticesofwave}) we have:

\begin{corollary} \label{largerroofsmore} For two waves $\cu$ and
$\cw$, if $\cw \extended \cu$ then $\cw \le \cu$.
\end{corollary}




The next lemma  is formulated in great generality (hence its
complicated statement), so as to avoid repeating the same type of
arguments again and again:

\begin{lemma}
\label{hetzgeneral} Let $X$ and $Y$ be two  sets of vertices in
$\Gamma$, and let $\cu, \cw$ be warps, satisfying the following
conditions:

(1)~$\cu$ is a wave in $\Gamma - X$.

(2)~$Y \subseteq RF_{\Gamma - X}(\cu)$.

(3)~ $\cw$ is a self roofing warp in $\Gamma - Y$.

(4)~$X \subseteq RF_{\Gamma - Y}(\cw)$ and $X \cap V[\cw] \subseteq in[\cw]$.

(5)~Every path in $\cw$ meets $RF_{\Gamma - X}(\cu)$.

Then $\ce_\Gamma(ter[\cu \hetz \cw]) = \ce_\Gamma(ter[\cu] \cup
ter[\cw]) = \ce_\Gamma(ter[\cu] \cup ter[\cw] \cup X \cup Y)$.

\end{lemma}

(The last equality means of course that $X \cup Y \subseteq
RF(ter[\cu] \cup ter[\cw])$.)

\begin{proof}
By (1) and (2) we have $Y \subseteq RF(X \cup ter[\cu])$ and
by (3) and (4) we have $X \subseteq RF(Y \cup ter[\cw])$.
This together with Observation \ref{STXY} yields $\ce(ter[\cu] \cup ter[\cw]) =
\ce(ter[\cu] \cup ter[\cw] \cup X \cup Y)$, so  we only need to show
the first equality. Since $ter[\cu \hetz \cw] \subseteq ter[\cu]
\cup ter[\cw]$, by Lemma \ref{essentialsandwiched} it suffices to
show that $ter[\cu \hetz \cw] \supseteq \ce(ter[\cu] \cup
ter[\cw])$.

Let $z \in \ce(ter[\cu] \cup ter[\cw])$. We need to show that $z \in
ter[\cu \hetz \cw]$.

Consider first the case that $z \in ter[\cu]$. If $z \not \in
V[\cw]$ then $\cu(z) \in \cu \hetz \cw$ and we are done. Thus we may
assume that $z \in V[\cw]$, which by (3) entails that $z \in
RF(ter[\cw] \cup Y)$ and $z \not \in Y$. The fact that $z \in
\ce(ter[\cu] \cup ter[\cw] \cup X \cup Y)$ implies therefore that $z
\in ter[\cw]$, again implying $\cu(z) \in \cu \hetz \cw$.

We are left with the case that $z \in ter[\cw] \setminus ter[\cu]$.
Let $W = \cw(z)$ and let $u$ be the last vertex in $W$ which is in
$RF_{\Gamma - X}(\cu)$. The existence of such $u$ is
certified by (5).
Note that by (4) we know that the path $uW$ does not meet $X$.
Therefore we may apply Lemma \ref{lastonpathbelongstoroof}
in the web $\Gamma - X$ and get $u \in ter[\cu] \cup \{z\}$. Suppose that
$u=z$. Then by by the choice of $u$ and by (2), we have $z \not \in
X \cup Y$. Since $z \in \ce(ter[\cu] \cup ter[\cw] \cup X \cup Y)$,
there exists a $z$--$B$ path avoiding $ter[\cu] \cup ter[\cw]$. By
the choice of $u$ and by (1), this path must meet $ter[\cu]$, and
the only vertex at which this can happen is $u$ itself,
contradicting the assumption of the present case. We have thus
proved that $u \neq z$, and thus $u \in ter[\cu]$. This implies that
$\cu(u)uW \in \cu \hetz \cw$, proving $z \in ter[\cu \hetz \cw]$.

\end{proof}

The most frequently used  case of this lemma will be that of $Y = X
= \emptyset$ :

\begin{lemma}\label{hetzofwaves}
If $\cu$ and $\cw$ are waves then so is $\cu \hetz \cw$.
\end{lemma}

\begin{proof}
Combine the lemma with the fact that $ter[\cu]$, and hence {\em a
fortiori} $ter[\cu] \cup ter[\cw]$,
 is $A$--$B$-separating.
 \end{proof}

Another case we shall use is in which $X = \emptyset$ but $Y$ is not
necessarily empty.

\begin{lemma}\label{roofminus}
If $\cu$ is a wave in $\Gamma$, $Y \subseteq RF(\cu)$ and $\cw$
is a wave in $\Gamma - Y$, then $\cu \hetz \cw$ is a wave in
$\Gamma$.
\end{lemma}

Taking $Y = \emptyset$ but $X$ not necessarily empty, and using
Corollary \ref{quo_self_roofing} we get:

\begin{lemma}
\label{grounding} Let  $X,Z$  be subsets of $V(\Gamma)$ such that $X
\subseteq Z$. Let $\cu$ be a wave in $\Gamma - X$ and let $\cw$ be a
wave in $\Gamma \quo  Z$. If every path in $\cw$ meets $RF_{\Gamma
-X}(\cu)$ then $\cu \hetz \cw$ is a wave in $\Gamma$.
\end{lemma}

\begin{remark}\label{remarkgrounding}
Since $\cu$ is a wave in $\Gamma - X$,
every path in $\cw$ starting at $A$ meets
$RF_{\Gamma - X}(\cu)$.
Therefore the only paths in $\cw$ for which
the assertion of meeting $RF_{\Gamma -X}(\cu)$
really needs to be checked are those starting at $Z$.
\end{remark}




By Corollary \ref{largerroofsmore} if $\cu$ and $\cw$ are waves,
then $\cu \le \cu \hetz  \cw$. Lemma \ref{hetzgeneral} implies more:

\begin{lemma} \label{hetzbeatsboth} For any two waves $\cu$ and
$\cw$ we have:
 $\cu, \cw \le \cu \hetz  \cw$.
\end{lemma}

\begin {lemma}\label{essentialofterhetz}
$\ce(ter[\cu \hetz \cw]) \cap RF^\circ(\cu) =\emptyset$.
\end{lemma}

\begin{proof}
$\ce(ter[\cu \hetz \cw]) \cap RF^\circ(\cu) \subseteq \ce(ter[\cu]
\cup ter[\cw]) \cap RF^\circ(ter[\cu] \cup ter[\cw]) = \emptyset$
\end{proof}

\begin{lemma}\label{zornconditionforwaves}
If $(\cw_\alpha:~ \alpha< \theta)$ is a $\fextended$-ascending
sequence of waves, then $\uparrow_{\alpha< \theta}\cw_\alpha$ is a
wave and $\uparrow_{\alpha< \theta}\cw_\alpha \geq \cw_\alpha$ for every
$\alpha < \theta$.
\end{lemma}

\begin{proof}
This is a direct corollary of Observation \ref{limter} and Lemma
\ref{limroof}.
\end{proof}

Since clearly $\uparrow_{\alpha< \theta}\cw_\alpha \fextends
\cw_\alpha$ for all $\alpha< \theta$, by Zorn's lemma  this implies:

\begin{lemma}\label{maximalwave}
In every web there exists a $\fextended$-maximal wave. Furthermore,
every wave can be forward extended to a $\fextended$-maximal wave.
\end{lemma}

One corollary of this lemma is that a hindered web contains a
maximal hindrance.

\begin{corollary}\label{maximalhindrance} If there exists in
$\Gamma$ a hindrance  then there exists in $\Gamma$ a
$\fextended$-maximal wave that is a hindrance.
\end{corollary}

Next we show that there is no real distinction between
$\fextended$-maximality and $\le$-maximality.

\begin{lemma}\label{maximalismaximal}
Any  $\fextended$-maximal wave (and hence also any
$\extended$-maximal wave) is  $\le$-maximal. If $\cv$ is a
$\le$-maximal wave then there does not exist a trimmed wave $\cw$
such that $\ce(\cv) \precneqq \cw$.
\end{lemma}

\begin{proof}
Assume first that  $\cv$ is a $\le$-non-maximal wave, i.e., there
exists a wave $\cw > \cv$, meaning that $RF(\cw) \supsetneqq
RF(\cv)$. By Lemma \ref{hetzbeatsboth} it follows that $\cv \hetz
\cw \neq \cv$, and since $\cv \hetz \cw \fextends \cv$ it follows
that $\cv$ is not $\fextended$-maximal and hence also not
$\extended$-maximal. This proves the first part of the lemma.

Assume next that $\cv$ is a $\le$-maximal wave. Let $\cu=\ce(\cv)$.
Suppose, for contradiction, that $\cu \precneqq \cw$ for some
trimmed wave $\cw$. This means that there exists some path
$W \in \cw \setminus \cu$.
Since $\cw$ is trimmed, $W$ is
finite. Write $t=ter(W)$.
Since $in[W] \subseteq A$ and we assume no edges enter $A$,
the only two possibilities are that  either $W$ is a proper forward extension
of some path in $\cu$ or $W$ does not meet $V[\cu]$ at all.
In both case we have $t \not \in ter[\cu]$.
Since $\cw$ is trimmed we have $t \not \in RF^\circ(\cw)$ and hence
$t \not \in RF^\circ(\cu)$. Thus $t \not \in RF(\cu)$,
which
implies that $\cw > \cv$, a contradiction.
\end{proof}

\begin{corollary}\label{allmaximalroofthesame}
If $\cu,\cv$ are each either $\extended$-maximal, or
$\fextended$-maximal, or $\le$-maximal waves, then $\cu \equiv \cv$.
\end{corollary}

\begin{proof}
By the lemma, in all cases $\cu$ and $\cv$ are $\le$-maximal. By
Lemma \ref{hetzbeatsboth}~~ $\cu \hetz \cv \ge \cu, \cv$, which, by
the $\le$-maximality of $\cu$ and $\cv$, implies that $RF(\cu \hetz
\cv) = RF(\cu)=RF(\cv)$. The last equality means that $\cu \equiv
\cv$.
\end{proof}

Thanks to Corollary \ref{allmaximalroofthesame}, we may speak about ``maximal waves",
without specifying whether we mean $\le$ or $\extended$ or
$\fextended$-maximality, as long as do this only in contexts involving
vertices roofed by the waves, or quotient over the waves, or other
properties that do not distinguish between equivalent waves.

\begin{observation}\label{aofgammaquowave}
If $\cw$ is a wave, then $A(\Gamma \quo \cw) = \ce(ter[\cw])$.
\end{observation}
\begin{proof}
Recall that $\Gamma \quo \cw$ is defined as $\Gamma \quo ter[\cw]$,
which in turn means that $A(\Gamma \quo \cw)=(A \cup ter[\cw])
\setminus RF^\circ(ter[\cw])$.  Since $\ce(ter[\cw])=ter[\cw]
\setminus RF^\circ(ter[\cw])$ we have $\ce(ter[\cw]) \subseteq (A
\cup ter[\cw]) \setminus RF^\circ(ter[\cw])$. Since $\cw$ is a wave,
$A \subseteq RF(\cw)$, implying that $A \setminus RF^\circ(\cw)
\subseteq ter[\cw]$, and hence $(A \cup ter[\cw]) \setminus
RF^\circ(ter[\cw]) \subseteq ter[\cw] \setminus RF^\circ(ter[\cw]) =
\ce(ter[\cw])$.
\end{proof}

\begin{lemma} \label{starisawave} If $\cw$ is a wave in $\Gamma$ and $\cv$ is a wave
in $\Gamma \quo \cw$ then $\cw*\cv$ is a wave in $\Gamma$.
\end{lemma}

\begin{proof} Let $P$ be a path from $A$ to $B$. We have to show that $P$ meets
$ter[\cw*\cv]$. Since $\cw$ is a wave, $P$ meets $ter[\cw]$. Let $t$
be the last vertex on $P$ belonging to $ter[\cw]$. Then clearly $t
\in \ce(ter[\cw])$ and $V(tP) \cap RF^\circ(\cw) = \emptyset$, and hence by Observation \ref{aofgammaquowave}
$tP$ is an $A(\Gamma \quo \cw)$--$B(\Gamma \quo \cw)$ path in $\Gamma \quo \cw$. Thus $tP$ meets $ter[\cv]$, and
since clearly $ter[\cv] \subseteq ter[\cw*\cv]$ it follows that $tP$ meets
$ter[\cw*\cv]$ and so does $P$, as required.
\end{proof}

\begin{lemma}\label{quotientovermaximal}
If $\cw$ is a $\extended$-maximal wave then $\Gamma\quo \cw$ is
loose.
\end{lemma}

\begin{proof}
Assume, for contradiction, that there exists a non-trivial wave
$\cv$ in $\Gamma\quo \cw = \Gamma \quo \ce(\cw)$. If all paths in
$\cv$ are singletons then, since $\cv$ is non-trivial, $\cv
\subsetneqq \langle ter[\ce(\cw)] \rangle$, contradicting the
definition of $\ce(\cw)$. Thus not all paths in $\cv$ are
singletons, and hence $\cw*\cv \succneqq \cw$, and since by Lemma
\ref{starisawave} $\cw*\cv$ is a wave this contradicts the
maximality of $\cw$.
\end{proof}

By Lemma \ref{allmaximalroofthesame}, the $\extended$-maximality in
the above lemma can be replaced by $\fextended$- or
$\le$-maximality.

\begin{lemma} \label{waveingammaminusx}
Let $X$ be a subset of $V \setminus A$, and let $\cu$ be a warp in
$\Gamma$ avoiding $X$, such that $\cu$ is a wave in $\Gamma -X$.
Then $\cu \quo X$
is a wave in $\Gamma \quo X$. Furthermore,

\begin{equation}\label{equationwaveingammaminusx}
RF_{\Gamma-X}(\cu) \setminus RF^\circ(X) \subseteq RF_{\Gamma \quo X
}(\cu \quo X).
\end{equation}
\end{lemma}

\begin{proof}


Note that $\langle \ce(X) \rangle \subseteq \cu \quo X$. Since
$A(\Gamma \quo X) \subseteq (RF_{\Gamma-X}(\cu) \setminus
RF^\circ(X)) \cup \ce(X)$, in order to prove that $\cu \quo X$ is a
wave in $\Gamma \quo X$  it suffices to prove
(\ref{equationwaveingammaminusx}). Let $Q$ be a path in $\Gamma \quo
X$ starting at a vertex $z \in RF_{\Gamma-X}(\cu) \setminus
RF^\circ(X)$ and ending in $B$. We have to show that $Q$ meets
$ter[\cu \quo X]$.
If $Q$ meets $X$ then it meets $\ce(X)$ and we are done. If not,
then the desired conclusion follows from the fact that $z \in
RF_{\Gamma-X}(\cu)$.
\end{proof}

A corollary of this lemma is that $\Gamma \quo X$ contains more
``advanced" waves than $\Gamma -X$:

\begin{corollary}\label{quotientisstronger}
If  $X$ and $\cu$ are as above, and if $\cv$ is a maximal wave in
$\Gamma \quo X$, then $RF_{\Gamma}(\cv) \supseteq RF_{\Gamma-X}(\cu)$
and $RF^\circ_{\Gamma}(\cv) \supseteq RF^\circ_{\Gamma-X}(\cu)$.
\end{corollary}

One advantage that the quotient operation has over deletion is the
following. Given two sets of vertices, $X_1$ and $X_2$, there is no
natural way of combining a wave in  $\Gamma -X_1$ with a wave in
$\Gamma -X_2$, so as to yield a third wave in some web. By contrast,
there does exist a natural definition of a combination of a wave
$\cw_1$ in $\Gamma\quo X_1$ with a wave $\cw_2$  in $\Gamma\quo
X_2$. Writing $X = \ce(X_1 \cup X_2)$, we can combine $\cw_1$ and $\cw_2$
by taking the warp $(\cw_1\quo X) \hetz (\cw_2\quo X)$.

\begin{lemma}
\label{hetzinquotient} Let $X_1,X_2 \subseteq V$, and write $X=\ce(X_1
\cup X_2)$. If $\cw_1$ is a wave in $\Gamma\quo X_1$ and $\cw_2$ is a
wave in $\Gamma\quo X_2$, then $(\cw_1 \quo X) \hetz (\cw_2 \quo X)$
is a wave in $\Gamma\quo X$. Moreover,

$$ RF_{\Gamma/X}((\cw_1 \quo X) \hetz (\cw_2 \quo X)) \supseteq
RF_{\Gamma/X}(\cw_1 \quo X) \cup RF_{\Gamma/X}(\cw_2 \quo X).
$$

\end{lemma}

\begin{proof}
Corollary \ref{quotientbyunioncor} and Lemma \ref{waveinquotient} imply that
$\cw_1 \quo X$ and $\cw_2 \quo X$ are both waves in $\Gamma \quo X$,
and hence by Lemma \ref{hetzofwaves} so is $(\cw_1 \quo X) \hetz
(\cw_2 \quo X)$. The second part of the lemma follows from Lemma \ref{hetzbeatsboth}.
\end{proof}

The next lemma is a special case of Lemma
\ref{zornconditionforwaves} that we will need.

\begin{lemma}\label{limitofwavesinquotient}
Let $(X_i ~:~ 0 \le i < \omega)$ be a $\subseteq$-ascending sequence
of subsets of $V \setminus A$. For each $i < \omega$, let $\cw_i$ be
a wave in $\Gamma \quo X_i$.
Write $X=\ce(\bigcup_{i < \omega}X_i)$. Then $\uparrow_{i < \omega}(\cw_i
\quo X)$ (taken as an up-arrow  of waves in $\Gamma \quo X$) is a
wave in $\Gamma \quo X$.
\end{lemma}

We conclude this section with two lemmas taken from
\cite{countablemenger}, whose proofs are rather technical and hence
will not be presented here:

\begin{lemma}\label{hinderedminusfinitelymanypointsishindered}
If $\Gamma$ is hindered and $X$ is a finite subset of $V\setminus A$
then $\Gamma-X$ is hindered. \end{lemma}

This is not necessarily true   if  $X$ is infinite.

\begin{lemma}\label{hinderedimplieswave} If
$\Gamma$ is unhindered, and $\Gamma-v$ is hindered for a vertex $v
\in V \setminus A$, then there exists a wave $\cw$ in $\Gamma$ such
that $v \in ter[\cw]$. \end{lemma}

\section{Bipartite conversion of webs and
warp-alternating paths}

\subsection{Aims of this section}
As already mentioned, Menger's theorem is better understood, in
both its finite and infinite cases, if its relationship to
K\"{o}nig's theorem is apparent. As mentioned in the introductin, a simple transformation,
observed in \cite{aharonifinitepaths} (but probably known
earlier), reduces the finite case of Menger's theorem to
K\"onig's theorem. This ``bipartite conversion" is effective also
for webs containing no infinite paths, but not for general webs.
We chose to describe it here since it inspired many of the ideas
of the present proof, and some points in the proof are illuminated
by it. The bipartite conversion is also the most natural source
for definitions involving alternating paths. As is common in
matching theory, the latter will constitute one of
our main tools.

\subsection{The bipartite conversion of a web}\label{conversion} The ``bipartite
conversion" turns a digraph into a bipartite graph. Every vertex
of the digraph is replaced  by two copies,
one sending arrows and the other receiving them. The graph becomes
then bipartite, with one side consisting of the ``sending" copies,
and the other consisting of the ``receiving" copies.

 For webs the construction is a little
 different: $A$-vertices are given only ``sending" copies, and
 $B$-vertices are given only ``receiving" copies. Thus the web $\Gamma=(G,A,B)$
 turns into a bipartite web
$\Delta=\Delta(\Gamma)=(G_\Delta,A_\Delta,B_\Delta)$, in the
 following way. Every vertex $v \in V\setminus A$ is assigned a
 vertex $w(v) \in B_\Delta$, and every vertex $v \in V \setminus B$ is
 assigned a vertex $m(v) \in A_\Delta$. Thus, vertices in $V
 \setminus (A \cup B)$ are assigned two copies each. The edge set
 $E_\Delta=E(G_\Delta)$ is defined as $\{(m(x),w(y))\mid ~(x,y) \in
 E(G)\} \cup \{(m(x),w(x)) \mid ~x \in V\setminus (A \cup B)\}$.

The above transformation converts a web into a bipartite web,
together with a matching, namely the set of edges $\{(m(x),w(x))
\mid  ~x \in V\setminus (A \cup B)\}$. This transformation can be
reversed: given a bipartite graph $\Delta$ whose two sides are
 $A$ and $B$, together  with a matching $J$ in it, one can construct from it a
  web $\Lambda=\Lambda(J)$ (the reference to
$\Delta$ is suppressed), as follows. To every edge $(x,y) \in J$
we assign a vertex $v(x,y)$. The vertex set $V(\Lambda)$ is
$\{v(x,y) \mid ~(x,y) \in J\} \cup V(\Delta) \setminus \bigcup J$.
(Here $\bigcup J$ is the set of vertices participating in edges
from $J$.)  The ``source" side $A_{\Lambda}$ of $\Lambda$ is
defined as $A_\Delta \setminus \bigcup J$, and the ``destination"
set $B_\Lambda$ is $B_\Delta \setminus \bigcup J$.

For $u \in V(\Lambda)$ define $m(u) = u$ if $u \in A_\Lambda
\setminus J$, and $m(v(x,y)) =  x$ (namely, the $A$-vertex of
$(x,y)$) for every edge $(x,y) \in J$. Let $w(u) = u$ if $u \in
B_\Lambda \setminus J$, and  $w(v(x,y)) =  y$ (namely, the
$B$-vertex of $(x,y)$) for every edge $(x,y) \in J$. The edge set
of $\Lambda$ is defined as $\{(u,v) \mid ~(m(u), w(v)) \in
E[\Delta]\}$.

Let us now return to our web $\Gamma$, and consider a warp $\cw$ in
it. Let $J=J(\cw)$ be the matching in $\Delta(\Gamma)$, defined by
$J =\{(m(u),w(v)) \mid~(u,v) \in E[\cw]\} \cup \{(m(u),w(u)) \mid u
\not \in \bigcup E[\cw]\}$.  We abbreviate and write $\Lambda(\cw)$
for $\Lambda(J(\cw))$. From the definitions there easily follows:

\begin{lemma}\label{linkabilityandmarriage}
If $\cw$ is a linkage in $\Gamma$, then $J(\cw)$ is a marriage of
$A_\Delta$ in $\Delta=\Delta(\Gamma)$. If $\Gamma$ does not contain
unending paths, then the converse is also true.
\end{lemma}

\subsection{Alternating paths}

The definition of one of our main tools, that of $\cy$-alternating paths, where $\cy$ is a warp, is quite involved. To be able to follow its fine points, it may be helpful to keep in mind the main property required of a $\cy$-alternating path: that the symmetric difference of its edge set and the edge set of $\cy$ is the edge set of a warp. For a precise definition, see Definition \ref{triangle} below.

\begin{definition}\label{alternatingpaths}
Let $\cy$ be a warp in $\Gamma$. A $\cy${\em -alternating path} is
a  sequence $Q$ having one of the following forms:

(i)~   an infinite sequence $(u_0,F_0, w_1,
R_1,u_1,F_1,w_2,R_2,u_2,\ldots)$,

(ii)~  an infinite sequence $(w_1,
R_1,u_1,F_1,w_2,R_2,u_2,\ldots)$,

(iii)~  $(u_0,F_0, w_1,
R_1,u_1,F_1,w_2,R_2,u_2,\ldots, R_k,u_k)$,

(iv)~ $(u_0,F_0, w_1,
R_1,u_1,F_1,w_2,R_2,u_2,\ldots, R_k,w_k,F_{k+1},w_{k+1})$,

(v)~   $(w_1,
R_1,u_1,F_1,w_2,R_2,u_2,\ldots, R_k,u_k)$,

(vi)~  $(w_1,
R_1,u_1,F_1,w_2,R_2,u_2,\ldots, R_k,w_k,F_{k+1},w_{k+1})$,\\

and satisfying the following
conditions:
\begin{enumerate}
\item $u_i,w_i$ are vertices, and $F_i,R_i$ are paths having at least one edge each. Furthermore, $in(F_i)=u_i,~ter(F_i)=w_{i+1}$ for all relevant values
of $i$. The paths $R_i$ are subpaths of paths from $\cy$, and
$in(R_i)=u_i,~ter(R_i)=w_i$ for all relevant values of $i$.

\item For paths of types (i),(iii) or (iv) $u_0 \not \in V[\cy]$.

\item For paths of types (iv) or (vi) $w_{k+1} \not \in V[\cy]$.

\item
If $v \in V(R_i) \cap V(R_j)$ for  $i \neq j$, then either
$v=u_i=w_j$ or $v=w_i=u_j$.
\item
If $v \in V(F_i) \cap V(F_j)$ for  $i \neq j$, then either
$v=u_i=w_{j+1}$ or $v=w_{i+1}=u_j$.
\item
If $V(F_i) \cap V(R_j) \neq \emptyset$ then either:\\
\indent (i) $j=i+1$, and $V(F_i) \cap V(R_j) =\{w_j\}$, or:\\
\indent (ii) $i>j$ and  $V(F_i) \cap V(R_j)\cap \{u_i,w_{i+1},u_j,w_j\} =\emptyset$, namely the paths $F_i$ and $R_j$ meet only at their interiors, or:\\
\indent  (iii) $j=i$, and  $w_i, w_{i+1} \not \in V(F_i) \cap V(R_j)$,
namely the paths $F_i$ and $R_j$ meet only at $u_i$ and possibly also at their interiors.
\end{enumerate}
The notation ``$F_i$" and ``$R_i$" stands for ``forward" and
``reverse", respectively - we think of $Q$ as going forward on
$F_i$, and reversely on $R_i$. The links $F_i$ and $R_i$ are
called ``forward links" and ``backward links" of $Q$, respectively.
The last three requirements in the definition mean that links can
only meet at their endpoints, with one exception: a forward link can go through an internal vertex of a backward link, if the latter precedes it in the path.
Allowing this may seem redundant, since if this happens then the alternating path can be replaced
by a shorter one having the same initial and terminal vertices. But there is one place, namely Lemma \ref{alpath_from_u_andback} below, in which
this type of paths must be permitted.

The first vertex ($u_0$ or $w_1$) on $Q$ is denoted by $in(Q)$. If $in(Q)=u_0 \in A$ then
$Q$ is said to be $A$-{\em starting}. Note that by condition (2)
this implies that $u_0 \not \in V[\cy]$.
If  $Q$ is infinite, then
$Q$ is said to be an $(in(Q), \infty)$-$\cy$- alternating path. If $Q$
is finite, we write  $ter(Q)$ for its last vertex, and say that $Q$ is an $(in(Q),ter(Q))$-$\cy$-alternating path. If
 $in(Q)=u_0 \in A \setminus
V[\cy]$ and $ter(Q) \in B\setminus V[\cy]$,  we say that $Q$ is
{\em augmenting}. The source of this name is that in this case
the application of
$Q$ to $\cy$ adds one more path to $\cy$ than it removes from it (see Definition \ref{triangle} below for the meaning of ``application" of an alternating path to a warp). This meaning of ``augmentation" does not depend on the condition $in(Q) \in A, ~ter(Q)\in B$, but since this is the only case we shall use the notion of augmentation, we add this condition.

If  $in(Q)=w_1 \in ter[\cy]$ and $ter(Q)=u_k \in in[\cy]$ then $Q$ is said to
be {\em reducing}. In this case the application of
$Q$ to $\cy$ removes one more path from $\cy$ than it adds to it.

If $Q$ is infinite, or it is finite and $ter(Q) \not \in V[\cy]$ then $Q$
is said to be $\cy$-{\em leaving}.

\end{definition}

\begin{definition}\label{triangle}
For a $\cy$-alternating path $Q$ as above, $\cy \triangle Q$ is
the cyclowarp whose edge set is $E[\cy] \triangle E(Q)$, namely $E[\cy]
\setminus \bigcup E(R_i) \cup \bigcup E(F_i)$, with $ISO(\cy
\triangle Q) = ISO(\cy)$.
\end{definition}

(Recall that $ISO(\cy)$ denotes the
set of singleton paths in $\cy$.)
The cyclowarp $\cy \triangle Q$ is also said to be the result of {\em
applying $Q$ to $\cy$}.

\begin{definition} Let $\cu, \cy$ be warps. A $\cy$-alternating
path is said to be $[\cu,\cy]$-alternating if all paths $F_i$ in
Definition \ref{alternatingpaths} are subpaths of paths in $\cu$.
A $[\cu,\cy]$-alternating path is said to be $\cu$-{\em comitted}
if no $R_i$ contains a point from $V[\cu]\setminus ter[\cu]$ as an
internal point. Namely, if the alternating path switches to $\cu$
whenever possible.
\end{definition}

Every  $\cy$-alternating path in $\Gamma$ corresponds in a natural
way to a $J(\cy)$-alternating path in $\Delta(\Gamma)$, which, in
turn, corresponds to a path  in $\Lambda(\cy)$. Moreover, an
augmenting $\cy$-alternating path corresponds to an
$A_\Lambda$--$B_\Lambda$ path in $\Lambda$. We summarize this in:

\begin{lemma}\label{alternatingpathsturneddirect}
Let $\cy$ be a warp in $\Gamma$, and let $\Lambda=\Lambda(\cy)$.
Then there exists an augmenting $\cy$-alternating path if and only
if there exists an $A_\Lambda$--$B_\Lambda$ path in $\Lambda$.
\end{lemma}


An  $A$--$B$-warp $\cy$ is called {\em strongly maximal} if $|\cy
\setminus \cu| \ge |\cu \setminus \cy|$ for every
 $A$--$B$-warp $\cu$.  The following is well known (see, e.g., \cite{mcdiarmid}):

\begin{lemma}\label{strongmax}
An  $A$--$B$-warp $\cy$ is   strongly maximal if and only if there
does not exist an augmenting $\cy$-alternating path.
\end{lemma}

 Note that in the
finite case ``strong maximality" means just ``having maximal
size", and hence obviously there
 exist strongly maximal warps. Hence the following
result implies Menger's theorem:

\begin{theorem} \label{blockingsimply} Let $\cy$ be a strongly maximal  $A$--$B$-warp. For
every $P \in \cy$ let $bl(P)$ be the last vertex on $P$
participating in an $A$-starting $\cy$-alternating path if such a vertex exists,
and $bl(P)=in(P)$ if there is no $A$-starting $\cy$-alternating path meeting
$P$. Then the set $BL=\{bl(P):~P \in \cy\}$ is
$A$--$B$-separating.
\end{theorem}

(The letters ``bl" stand for ``blocking".) This result also yields
an equivalent formulation of Theorem \ref{main}, noted in
\cite{mcdiarmid}: in every web there exists a strongly maximal
$A$--$B$-warp.

Theorem  \ref{blockingsimply} was proved by Gallai \cite{gallai}.
A detailed proof is given in Chapter 3 of \cite{diestel}. We
give here an outline of the proof, since it yields one of
the simplest proofs of the finite case of Menger's theorem, and
since the idea will recur in Section 8.

{\em Proof of Theorem \ref{blockingsimply}.}~Let $T$ be an
 $A$--$B$-path. Let $P$ be the first path from $\cy$ it meets, say at
a vertex $z$. Assuming that $z\neq bl(P)$, it must precede $bl(P)$
on $P$, since it lies on the  alternating path $Tz$. Assuming that
$T$ avoids $BL$, it follows that either:

(i)~ $T$ meets a path $R \in \cy$ at a vertex $u \in V(R)$
preceding $bl(R)$ on $R$, and  $uT-u$ is disjoint from
$V[\cy]$,~~or:

(ii)~$T$ meets a path $R \in \cy$ at a vertex $u \in V(R)$
preceding $bl(R)$ on $R$, and the next vertex $w$ on $T$ belonging
to $V(W)$ for some $W \in \cy$ comes after $bl(W)$ on$W$.

Assume that (i) is true. Let $Z$ be a $\cg$-alternating path from
$bl(R)$ to $Y \setminus S$. If $Z$ does not meet $T$, then
$Tu\overleftarrow{R}bl(R)Z$  is an augmenting $\cg$-alternating
path, contradicting Lemma \ref{strongmax}. If $Z$ meets $T$, let
$z$ be the last vertex on $Z$ belonging to $V(T)$. Then the path
$TzZ$ is again an augmenting $\cg$-alternating path, again
yielding a contradiction.

 On the other hand, (ii) is impossible
since
 the alternating path reaching $bl(R)$
can be extended by adding to it $\overleftarrow{R}uTw$, so as to
form an  alternating path meeting $W$ beyond $bl(W)$. $\enp$

\subsection{Safe alternating paths}

\begin{definition}\label{saps}
A $\cy$-alternating path $Q$ is called {\em safe} if:
\begin{enumerate}
\item
For every $P \in \cy$ the intersection $E[Q] \cap E(P)$ (which, in the notation of Definition \ref{alternatingpaths}, is
$\bigcup E(R_i) \cap E(P)$)
 is the edge set of a
subpath (that is, a single interval) of $P$, and:
\item
$E(Q) \setminus E[\cy]$ does not contain an infinite path or a cycle.
\end{enumerate}
\end{definition}

 We use the abbreviation ``$\cy$-s.a.p" for ``safe
$\cy$-alternating path". A $\cy$-s.a.p whose forward links $F_i$ are
fragments of a warp $\cw$ is called a $[\cw,\cy]$-s.a.p.

\begin{figure}[htb]
\begin{center}
\texorpdfblock{
\psfrag{A}{$A$}
\includegraphics[scale=0.6]{painting_top.eps}\hspace*{3em}
\includegraphics[scale=0.6]{painting_bot.eps}
}{
\includegraphics[scale=0.6]{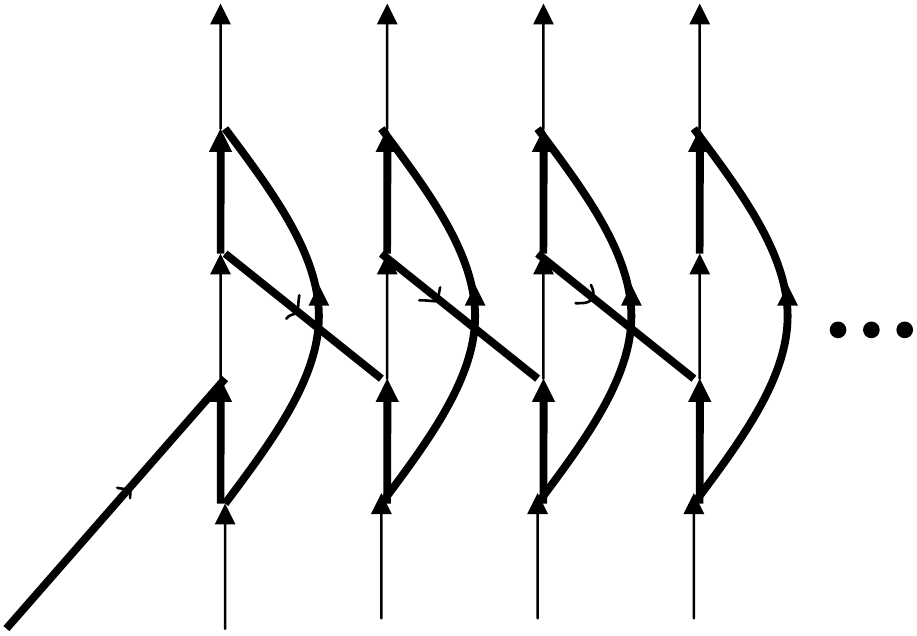}\hspace*{3em}
\includegraphics[scale=0.6]{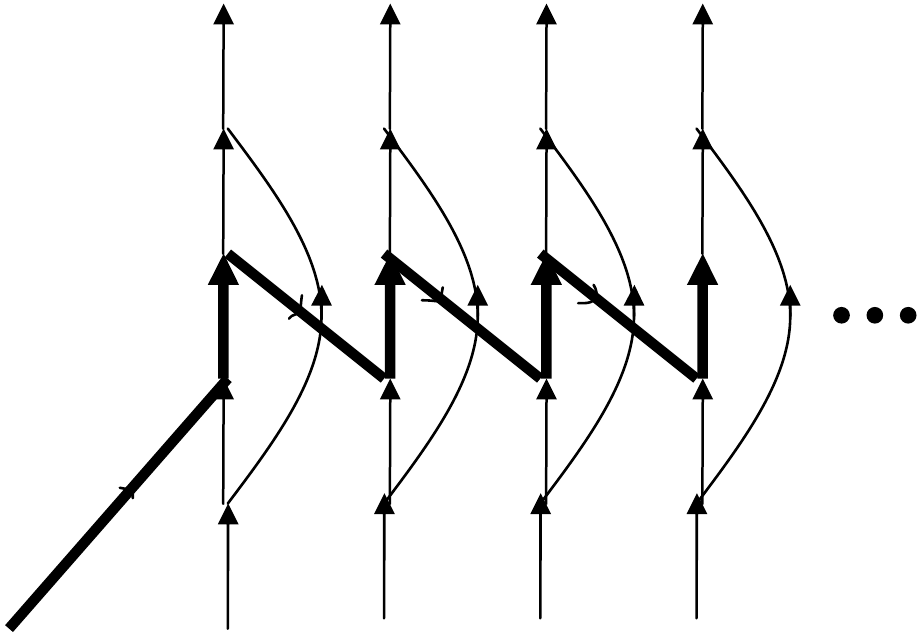}
}
\end{center}
\caption{An example of an alternating path (bolded on the left)
whose application results in a warp including an infinite path
(bolded on the right).} \label{badaltpath}
\end{figure}

If $Q$ is an infinite $\cy$-alternating path then $\cy \triangle Q$
may contain infinite paths, even if $\cy$ itself is f.c (reminder -
``f.c." means ``of finite character", namely having no infinite
paths). See Figure \ref{badaltpath}.

The name ``safe" originates in the fact that this cannot occur if
$Q$ is safe. For, each path in $\cy \triangle Q$ consists then of
only three parts (one or two of which may be empty) - a subpath of a
path of $\cy$, followed by a path lying outside $\cy$, followed then
by another subpath of a path of $\cy$. For the same reason, $\cy
\triangle Q$ does not contain cycles. We summarize this in:

\begin{lemma}\label{safeisindeedsafe}
If  $\cy$ is warp of f.c. and $Q$ is a $\cy$-s.a.p, then also $\cy
\triangle Q$ is a warp of f.c.
\end{lemma}

\begin{definition}\label{degenerate}
A $(u,v)$-$\cy$-alternating path $Q$ (where possibly $v=\infty$)
is called {\em degenerate} if $\cy \triangle Q$ contains a path
from $u$ to $v$.
\end{definition}

The definition of ``safeness" implies:

\begin{lemma}\label{safeisnondegenerate}
If a $(u,v)$-$[\cw,\cy]$-s.a.p $Q$ is degenerate, then the path
connecting $u$ to $v$ in $\cy \triangle Q$ is contained in a path
from $\cw$.
\end{lemma}

A fact that we shall use about s.a.p's is:

\begin{theorem}\label{czcysaps} Let $\cz$ and $\cy$ be f.c. warps,
such that $in[\cz] \supseteq in[\cy]$. Then there exists a choice
of a $z$-starting $\cy$-leaving maximal s.a.p $Q(z)$ for each $z
\in in[\cz] \setminus in[\cy]$, such that those s.a.p's $Q(z)$ that are
finite end at distinct vertices of $ter[\cz]$ (namely, $ter(Q(z))\neq ter(Q(z'))$ whenever $z\neq
z'$ and $Q(z),Q(z')$ are finite. Note that the paths $Q(z)$ themselves are not required to be disjoint).
\end{theorem}

The maximality of the paths $Q(z)$ means that each  $Q(z)$ is continued either indefinitely or
until a vertex of $ter[\cz]\setminus V[\cy]$ is reached.   For the proof of the theorem we
shall need the following lemma:

\begin{lemma}\label{alpath_from_u_andback}
Let $\cz$ and $\cy$ be f.c warps
such that $in[\cz] \supseteq in[\cy]$,
and let $u \in in[\cz] \setminus V[\cy]$. Then at
least one of the following possibilities occurs:
\begin{enumerate}
\item
There exists a $(u,\infty)$-$[\cz,\cy]$-s.a.p, or:

\item
There exists a vertex $v \in ter[\cz]\setminus V[\cy]$ for which there exist both
a $(u,v)$-$[\cz,\cy]$-s.a.p
and a $(v,u)$-$[\cy,\cz]$ alternating path.
\end{enumerate}
\end{lemma}

Note that in case (2), the $(u,v)$-$[\cz,\cy]$-s.a.p must be of type (iv) and
the $(v,u)$-$[\cy,\cz]$ alternating path must be reducing of type (v).
To follow the logic of the proof, keep in mind that $[\cy,\cz]$-alternating paths
are  $\cz$-alternating, but not necessarily $\cy$-alternating. Namely, they are ``$\cz$-committed", meaning that  whenever they meet a $\cz$-path they must switch to it, but they are not ``$\cy$-committed". In contrast, $[\cz,\cy]$-alternating paths are $\cy$-committed, while not necessarily $\cz$-committed. The following two examples illustrate this point:

\begin{example}\label{safeexample}
Suppose that $\cy$ consists of one path, $Y=(a,b,c,d)$,
while $\cz$ consists of the paths $(a,d), (s,b,t)$ and $(x,c,y)$.
Consider first a case in which $u=x$.
Since the graph is finite, (2) is impossible, and hence (1) should hold.
Indeed,  the easiest way to show that this is true is to take $v=t$.
The safe $(u,v)-$$[\cz,\cy]$-alternating path  (written by the order of its vertices)
whould then be $(x,c,b,t)$. Since this alternating path is $\cz$-committed,
we can choose its reverse $(t,b,c,x)$ to be the $(v,u)$-$[\cy,\cz]$-alternating path
required in (1). Note, however, that the choice of $v=t$ is not unique.
We could also take $v=y$, with the safe $[\cz,\cy]$-$(u,v)$-alternating path
$(x,c,b,a,d,c,y)$ and the $(v,u)$-$[\cy,\cz]$-alternating path $(y,c,x)$.

Consider next another case, in which $u=s$.
In this case, if we try to construct a $\cz$-committed
alternating path, we end up with the alternating path $(s,b,a,d,c,y)$, which is not safe.
The only way to obtain (1) is then taking $v=t$. The safe $(u,v)$-$[\cz,\cy]$-alternating path is $(s,b,a,d,c,b,t)$,
and the $[\cy,\cz]$-alternating path  is $(t,b,s)$.
\end{example}


{\em Proof of the lemma:}
By duplicating edges when necessary we may assume that $E[\cz] \cap
E[\cy] = \emptyset$. It is clear that vertices on paths from $\cy$ not belonging to $V[\cz]$ do not play any role in the proof. They can be ignored, meaning that subpaths having them
as internal points can be made to be single edges, and then terminal points of the resulting warp not belonging to $V[\cz]$ can be removed. Hence we shall assume that $V[\cy] \subseteq V[\cz]$.

Let $SR$ (standing for ``safely reachable")  be the set of vertices $v \in ter[\cz] \setminus V[\cy]$
for which there exists a $(u,v)$-$[\cz,\cy]$-s.a.p, and let $C$ be the set of vertices $x$ for which
there exists a $(v,x)$-$[\cy,\cz]$- alternating
path $T(x)$ for some $v \in SR$. Note that $T(x)$ is not necessarily unique, but to avoid cumbersome phrasing we shall sometimes pretend that it is. Thus we shall refer by $T(x)$ to {\em some}
alternating path satisfying the above conditions.

Assuming negation of possibility (2) of the lemma,
we have $u \not \in C$. Our aim is to show that this implies possibility (1) of the lemma.
To that end, we construct a $u$-starting
$[\cz,\cy]$-s.a.p $S$. This is done in stages, where at the $i$-th stage
we shall have at hand a $u$-starting
$[\cz,\cy]$-s.a.p $S_i$ extending $S_{i-1}$, whose last link is
a backward link on some path  $Y_i \in \cy$,
ending at a vertex $u_i \not \in C$ (the paths $Y_i$ need not be distinct).
The construction will be shown to go on indefinitely, meaning that possibility (1) of the lemma is true.

Let $u_0=u$. The path  $Z_1=\cz(u)$ must meet some  path in $\cy$, or else $ter(Z_1) \in SR$, meaning that $Z_1$
 can serve as  $T(u)$ to show that $u \in C$.
Let $w_1$ be the first vertex on $Z_1$ lying on some
path $Y_1\in \cy$. Since $u \not \in V[\cy]$, we know that $w_1 \neq u$.
Since $w_1 \not \in in[\cz]$ and $in[\cy] \subseteq in[\cz]$ we also know that
$w_1 \neq in(Y_1)$.
Let $y$ be the vertex preceding $w_1$ on $Y_1$.

\begin{assertion} \label{w1notinc}
Neither $w_1$ nor $y$ are in $C$. \end{assertion}

If $w_1 \in C$ then extending $T(w_1)$ by the path $u\cz(u)$ (either
as an extension of a link, in case the last link on $T(w_1)$ is a $\cz$-link, or as a separate link, if
the last link on $T(w_1)$ is a $\cy$-link,)
would show that $u \in C$. Similarly, if $y \in C$ then $T(y)$ can be extended by
the forward link $yY_1w_1$ followed by the backward link
$u\cz(u)$ to show that $u \in C$.

Returning to the construction of $S$, we go back on $Y_1$ to the first vertex $u_1$ on $Y_1$  not belonging to $C$ (possibly $u_1=in(Y_1)$). By the assertion,
this means going at least one edge back on $Y_1$, meaning that the path $S_1$ obtained is $\cy$-alternating.

Let $B_1=V((u_1Y_1w_1)^\circ)$ - this is the set of vertices on which $S_1$ goes backwards.
(Recall that $P^\circ$ is obtained from a path $P$ by removing from it $in(P)$ and $ter(P)$.)

Assume now that $i \ge 1$, and that $S_i=(u_0,F_1,w_1,R_1,\ldots R_i,u_i)$ has been already defined, where each
forward link $F_i$ is
part of a path $Z_i=\cz(u_{i-1})$, the
backward
link $R_i$ is part of a path $Y_i=\cy(w_i)$ and $u_i=ter(S_i) \not \in C$. We also assume that $S_i$ is safe.
Denote by $B_i$ the set of inner points of the backward links of $S_i$, namely $B_i = \bigcup_{j \le i}V(R_i^\circ)$.

Let $Z_{i+1}=\cz(u_i)$.

\begin{assertion}
$u_iZ_{i+1}$ meets $V[\cy] \setminus B_i$.
\end{assertion}

Assuming negation of the assertion, the alternating path obtained from $S_i$ by adding to it the link $u_iZ_{i+1}$
shows that $ter(Z_{i+1}) \in SR$ (it is this argument for which we need to allow
 alternating paths to go through previous backward links). Then the path $u_iZ_{i+1}$ can serve as $T(u_i)$, to show that $u_i \in C$, a contradiction.

Let $w_{i+1}$ be the  first vertex on $u_iZ_{i+1}$ belonging to $V[\cy] \setminus B_1$. Let $Y_{i+1}=\cy(w_{i+1})$.
Since $w_{i+1}\neq in(Z_{i+1})$ and $in[\cy]\subseteq in[\cz]$, we have $w_{i+1}\neq in(Y_{i+1})$. Let $y_{i+1}$ be the vertex preceding $w_{i+1}$ on $Y_{i+1}$.

\begin{assertion} \label{wnotinc} $w_{i+1} \not \in C$.\end{assertion}
Assuming for contradiction that $w_{i+1} \in C$, concatenating $T(y_{i+1})$ with $u_iZ_{i+1}w_{i+1}$ would yield a path $T(u_i)$, showing that $u_i \in C$. Note that $T(u_i)$ ``ignores" the meeting with
$Y_i$ at $w_i$, but this is fine, since it needs not be $\cy$-committed.

\begin{assertion} \label{ynotinc} $y_{i+1} \not \in C$.\end{assertion}

Assuming for contradiction that $y_{i+1} \in C$, concatenating $T(y_{i+1})$ with the single edge link $(y_{i+1},w_{i+1})$ and then with $u_iZ_{i+1}w_{i+1}$ would yield a path $T(u_i)$, showing that $u_i \in C$ (again, remember that $T(u_i)$ needs not be $\cy$-committed).

The last assertion, the fact that $w_{i+1} \not \in B_i$, and the choice of $u_i$ as the first vertex on $Y_i$ not belonging to $C$, imply:
\begin{assertion}\label{w2after}
If $Y_{i+1}=Y_j=Y$ for some $j\le i$ then $w_{i+1}>_Yw_j$.
\end{assertion}

We continue the construction of the alternating paths $S_i$, adhering to the following two rules:

{\em Rule 1:} If $Y =Y_i \in \cy$ is met for the first time, we go
on it backwards until we reach the first vertex $u_i$ on $Y$ not belonging to $C$.

{\em Rule 2:} If $Y =Y_i \in \cy$ has already been met, we go
backwards on $Y$ until we reach a vertex $w=w_j$ for some $j<i$,
and let $u_i=w_j$.

Since by the induction hypothesis $w_j \not \in C$ for $j <i$, by Assertion \ref{wnotinc} when Rule 2 is applied we still have the condition $u_i \not \in C$.
Rule 2 guarantees that the alternating paths $S_i$ constructed are safe. As noted above,
the condition $ter(S_i) \not \in C$ implies that the construction continues indefinitely, and
generates an infinite  $[\cz,\cy]$-alternating path $S$. In fact, $S$
is safe, since
Condition (1) of Definition \ref{saps} follows from the construction, while Condition (2) is true
since the non-$\cy$ links in $S$ come from $\cz$, which is f.c.
This proves that possibility (1) of the lemma holds.
$\enp$

{\em Proof of Theorem \ref{czcysaps}} The connected components of
the graph whose edge set is $E[\cz]\cup E[\cy]$ are countable.
Hence we may assume that $\cz$ and $\cy$ are countable. Let $z_1,
z_2, \ldots$ be an enumeration of  $in[\cz] \setminus in[\cy]$.
Applying Lemma \ref{alpath_from_u_andback} with $u=z_1$ we obtain
a $z_1$-starting $[\cz,\cy]$-s.a.p $Q_1$, satisfying condition (1)
or (2) of the lemma. If (1) is true, continue by applying the
lemma to $z_2$. If (2) is true, denote the vertex $v$ appearing in
the lemma by $v_1$, and the
$(v,z_1)$-$[\cz,\cy]$- alternating
path by $T_1$. Then $\cz_1=(\cz \triangle
T_1)^{path}$ is a f.c. warp, with $in[\cz_1]=in[\cz] \setminus \{z_1\}$,
$ter[\cz_1]=ter[\cz] \setminus \{v_1\}$.
(Recall that $(\cz \triangle T_1)^{path}$ is the warp obtained from $\cz \triangle
T_1$ by removing its cycles. Such cycles might appear since $T_1$ is not required
to be safe.)
Apply now the lemma to
the pair $(\cz_1,\cy)$, with $u=z_2$.

Continuing this way, we obtain a sequence $Q_i$ of $z_i$-starting
$\cy$-s.a.p's, which are either infinite or end at distinct
vertices of $ter[\cz]$, as promised in the theorem. $\enp$

\begin{remark}\label{fracturedczcy} The theorem applies also when $\cz$ is a fractured warp.
Reducing the fractured case to the non-fractured case is done by
duplicating those vertices which serve as both an initial point and a terminal point of paths from $\cz$,
thus turning $\cz$ into a proper warp.
\end{remark}

\section{A Hall-type equivalent conjecture}
In \cite{countablemenger} Theorem \ref{erdos} was shown to be
equivalent to the following Hall-type conjecture:

\begin{conjecture}\label{halltype}
An unhindered web is linkable.
\end{conjecture}

Both implications in this equivalence are quite easy. To show how
Theorem \ref{erdos} implies Conjecture \ref{halltype}, suppose
that Theorem \ref{erdos} is true, and let $\cp$ and $S$ be as in
the theorem. Then $\{Ps:~ P\in \cp,s \in V(P) \cap S\}$ is a
wave, and unless $\cp$ is a linkage, it is also a hindrance. To
prove the converse implication, take a
$\extended$-maximal wave $\cw$ in $\Gamma$ (see Lemma
\ref{maximalwave}), and let $S=ter[\ce(\cw)]$. By Lemma
\ref{quotientovermaximal}, $\Gamma\quo S$  is loose, and in
particular unhindered. Assuming that Conjecture \ref{halltype} is
true, the web $\Gamma \quo  S$ has therefore a linkage $\cl$.
Taking $\cp = \cw * \cl$ then fulfils, together with $S$, the
requirements of Theorem \ref{erdos}.

In fact, the above argument shows that the following is also
equivalent to Theorem \ref{erdos}:

\begin{conjecture}\label{halltypeloose}
A loose web is linkable.
\end{conjecture}

Here is a third equivalent formulation, generalizing Theorem
\ref{infmarriage}:

\begin{conjecture}\label{halltypeunlikable}
If $\Gamma$ is unlinkable then there exists an $A$--$B$-separating
set $S$ which is linkable into $A$ in $\overleftarrow{\Gamma}$,
but $A$ is not linkable into $S$ in $\Gamma$.
\end{conjecture}

The main result of this paper is that Conjecture \ref{halltype},
and hence also Theorem \ref{erdos}, are true for general graphs.
Let us thus re-state the conjecture, this time as a theorem:

\begin{theorem}\label{halltypetheorem}
An unhindered web is linkable.
\end{theorem}

The rest of the paper is devoted to the proof of Theorem
\ref{halltypetheorem}. The proof is divided into two stages. We
first define a notion of a $\kappa$-{\em hindrance} for every
regular uncountable cardinal $\kappa$, and show that
the existence of a $\kappa$-hindrance
implies the existence of a hindrance.
Then we shall show that if a web
is unlinkable then it contains either a hindrance or a
$\kappa$-hindrance for some uncountable regular $\kappa$.

\section{Safely linking one point}

In this section we prove a result, whose key role was already
mentioned in the introduction:

\begin{theorem}\label{safelinking}
If $\Gamma$ is unhindered then for every $a \in A$ there exists an
 $a$-$B$-path $P$ such that  $\Gamma -P$ is unhindered.
\end{theorem}

Let us first outline the proof of the theorem in the case of
countable graphs. This will serve two purposes: first, the main
idea of the proof appears also in the general case; second, it
will help to clarify the obstacle which arises in the uncountable
case. A main ingredient in the proof is the following:

\begin{lemma}\label{separatingoutneighbors}
Let $Q \subseteq V\setminus (A \cup B)$, and let $\cu$ be a wave
in $\Gamma -Q$, such that
\begin{equation}\label{sepY}
N^+(Q)\setminus Q \subseteq RF_{\Gamma - Q}(\cu) ~ .
\end{equation}
Then $\cu$ is a wave in $\Gamma$.
\end{lemma}

\begin{proof}
Let $P$ be an  $A$--$B$-path. We have to show that $P$ contains a
vertex from $ter[\cu]$. If $P$ is disjoint from $Q$ then, since
$\cu$ is a wave in $\Gamma-Q$, $P$ contains a vertex from
$ter[\cu]$. If $P$ meets $Q$ then, since $Q \cap B =\emptyset$,
there exists a vertex $y \in V(P) \cap N^+(Q)\setminus Q$. Choose
$y$ to be the last such vertex on $P$. By (\ref{sepY}), the path
$yP$ then contains a vertex belonging to $ter[\cu]$, as desired.
\end{proof}

{\em Proof of Theorem \ref{safelinking} for countable webs.~}
Enumerate all  $a$-$B$-paths as $P_1, P_2, \ldots$. Assuming that the
theorem fails, there exists a first vertex $y_1$ on $P_1$, such
that  $\Gamma - P_1 y_1$ is hindered.  Let $T_1=P_1 y_1-y_1$. Then
 $\Gamma -T_1$ is unhindered.  By Lemma \ref{hinderedimplieswave},
there exists a wave $\cw_1$ in  $\Gamma-T_1$ such that  $y_1 \in
ter[\cw_1]$. Let $i_2$ be the first index (if such exists) such
that $P_{i_2}$ does not meet $V[\cw_1]$.
Let $z$ be the last
vertex on $P_{i_2}$ lying on $T_1$, and let $P'_2= T_1zP_{i_2}$.
We may assume $\Gamma - P'_2$ is unhindered and hence
by Lemma \ref{hinderedminusfinitelymanypointsishindered}, the web
$\Gamma - T_1-zP_{i_2}$ is also hindered, since it is obtained from
$\Gamma - P'_2$ by removing finitely many vertices. Let $y_2$ be the first vertex on
$zP_{i_2}$ such that $\Gamma - T_1-zP_{i_2}y_2$ is hindered, and let $T_2=
T_1 \cup (zP_{i_2}y_2-y_2)$. By Lemma \ref{hinderedimplieswave}, there
exists a wave $\cw_2$ in $\Gamma -T_2$, such that $y_2\in
ter[\cw_2]$.

Continuing this way, we obtain an ascending sequence of trees
$(T_i:~i < \mu)$ (where $\mu$ is either finite or $\omega$), all
rooted at $a$ and directed away from $a$, and a sequence of waves
$\cw_i$ in $\Gamma - T_i$ disjoint from all trees $T_j$, such that
every $a$-$B$-path contains a vertex separated by some $\cw_i$
from $B$. Let $T= \bigcup_{i<\mu} T_i$ and $\cw=\uparrow \cw_i$.
Each $\cw_i$ is a wave in $\Gamma - T_i$ and hence also in $\Gamma - T$.
Therefore Lemma \ref{zornconditionforwaves} implies that
$\cw$ is a wave in $\Gamma -T$, separating from $B$ at least one
vertex from each $a$-$B$ path. By Lemma
\ref{separatingoutneighbors}, $\cw$ is a wave in $\Gamma$, and
since $a \not \in in[\cw]$, it is a hindrance, contradicting the
assumption of the theorem, that $\Gamma $ is unhindered. $\enp$

The difficulty in extending this proof beyond the countable case
is that after $\omega$ steps the web  $\Gamma-T_\omega$ may be
hindered, and then we can not proceed with the same construction,
since, for example, Lemma \ref{hinderedimplieswave} is not
applicable. Here is a brief outline of how this difficulty is
overcome.

Why was the construction of the trees $T_i$ necessary, and why
wasn't it possible just to delete the initial parts of the paths
$P_i$, and consider the waves (say) $\cu_i$ resulting from those
deletions? Because then each $\cu_i$ lives in a different web, and
it is impossible to combine the waves $\cu_i$ to form one big
wave. This we shall solve by taking quotient, instead of deleting
vertices - as we saw in Lemma \ref{hetzinquotient} it is then
possible to combine the resulting waves. But then we obtain a wave
which is not a wave in $\Gamma$, but in some quotient of it,
namely it does not necessarily start in $A$, while for the final
contradiction we need a wave (in fact, hindrance) in $\Gamma$
itself. This we overcome by performing the proof in two stages. In
the first we take quotients, and obtain a wave $\cw$ ``hanging in
air" in $\Gamma \quo X$ for some {\em countable} set $X$ (keeping
$X$ countable is a key point in the proof). In the second stage we
use the countability of $X$ to delete its elements one by one, in
a way similar to that used in the countable case, described above.
This process will generate a wave $\cv$, and the ``arrow"
concatenation of $\cv$ and $\cw$ will result in the desired wave
in $\Gamma$.

{\em Proof of Theorem \ref{safelinking}}.~ Construct inductively
trees $T_\alpha$ rooted at $a$ and directed away from $a$, as
follows. The tree $T_0$ consists of the single vertex $a$. For
limit ordinals $\beta$ define $T_\beta=\bigcup_{\alpha
<\beta}T_\alpha$. Assume that $T_\alpha$ is defined. Suppose first
that there exists a vertex $x \in V \setminus (A \cup
V(T_\alpha))$ such that $(u,x) \in E$ for some $u \in
V(T_\alpha)$, and $\Gamma -a - F -x$ is unhindered for every
finite subset $F$ of $V(T_\alpha)$ not including $a$. In this case
we choose such a vertex $x$, and construct $T_{\alpha+1}$ by
adding $x$ to $V(T_\alpha)$ and $(u,x)$ to $E(T_\alpha)$. If no
 vertex $x$ satisfying the above conditions exists, the process of definition is terminated at
$\alpha$, and we write $T=T_\alpha$.

The tree $T$ thus constructed has the property that for every
finite subset $F$ of $V(T)$ not including $a$ the web  $\Gamma -
a-F$ is unhindered, and $T$ is maximal with respect to this
property. Write $Y = N^+(V(T)) \setminus V(T)$. Then for every $y
\in Y$ there exists a finite set $F_y\subseteq V(T) \setminus
\{a\}$ such that $\Gamma-a - F_y -y$ is hindered.
Thus, by Lemmas
\ref{hinderedimplieswave}, there
exists a wave $\cu$ in $\Gamma-a - F_y$ with $y \in ter[\cu]$.
Let us now fix some maximal wave in $(\Gamma-a) \quo F_y$ and call it $\ca_y$.
Corollary \ref{quotientisstronger} yields $y \in RF_\Gamma(\ca_y)$.

Assuming that Theorem \ref{safelinking} fails, we have:

\begin{equation}
V(T) \cap B=\emptyset.
\end{equation}

Call a vertex $t \in V(T)$ {\em bounded} if there exists a
countable subset  $G_t$ of $V(T)$ and a wave $\cb=\cb_t$ in
$(\Gamma-a) \quo  G_t$ such that $t \in RF_\Gamma^\circ(\cb)$. Let $Q$ be
the set of non-bounded elements of $V(T)$. For every bounded
vertex $t \in V(T)$ choose a fixed set $G_t$ and a fixed wave
$\cb_t$ as above.

Let $\Gamma'=\Gamma-Q-a$. The core of the proof of Theorem
\ref{safelinking} is in the following:

\begin{proposition}\label{mainassertion}
 For every $y \in Y$ there exists a wave $\cu_y$ in $\Gamma'$
satisfying $y \in RF_{\Gamma'}(\cu_y)$.
\end{proposition}

{\em Proof of  the proposition:} Let $y$ be a fixed element of
$Y$. We shall construct a countable subset $X$ of $V(T) \setminus
A$, and a wave $\cw$ in $(\Gamma-a)\quo X$, having the following
properties:

(a)~ $y \in RF(\cw) $.

(b)~$F_z \subseteq X$
for every $z \in Y \cap
V[\cw \langle X \rangle]$.

(c)~$G_t \subseteq X$ and $t \in RF^\circ_\Gamma(\cw)$ for every $t \in X
\setminus Q$.

(d)~$V[\cw \langle X \rangle] \cap V(T) \subseteq X$.

The construction is by a ``closing up" procedure. We construct an
increasing sequence of sets $X_i$ whose union is to be taken as
$X$, and waves $\cw_i$ in $(\Gamma-a)\quo X_i$ whose ``$\uparrow$"
limit will eventually be taken as $\cw$, and at each step  we take
care of conditions (b) and (c), alternately, for all vertices $z
\in Y \cap V[\cw_i \langle X_i \rangle]$ and $t \in X_i \setminus
Q$. We shall do this in steps, as follows.

We take $X_0 = F_y$ and let $\cw_0 = \ca_y$.

For every $i < \omega$
we then take
$X_{i+1} = X_i \cup  \bigcup_{z \in Y \cap
V[\cw_i \langle X_i \rangle]} F_z \cup \bigcup_{t \in X_i
\setminus Q} G_t \cup (V[\cw_i \langle X_i \rangle] \cap V(T))$

and let $\cw_{i+1}$ be a maximal wave in
$(\Gamma-a) \quo X_{i+1}$.

Let $X = \bigcup_{i < \omega} X_i$ and $\cw= \uparrow_{i < \omega} (\cw_i \quo X)$.
Note that
for every $z \in Y \cap
V[\cw \langle X \rangle]$ we have
$z \in Y \cap
V[\cw_i \langle X_i \rangle]$ for some $i<\omega$.
This implies $F_z \subseteq X_{i+1} \subseteq X$
proving condition (b).
For every $t \in X \setminus Q$, we have $t \in X_i$ for some
$i < \omega$ and hence $G_t \subseteq X_{i+1}$. Since
$\cw_{i+1}$ is a maximal wave in
$(\Gamma-a) \quo X_{i+1}$, we have $t \in RF^{\circ}_{(\Gamma-a) \quo X_{i+1}}(\cw_{i+1})$.
If $t \in \ce(X)$, this implies
$t \in RF^{\circ}_{(\Gamma-a) \quo X}(\cw_{i+1} \quo X)
\subseteq RF^{\circ}_{(\Gamma-a) \quo X}(\cw)
\subseteq RF^{\circ}_\Gamma (\cw)$,
yeilding condition (c).
Of course, if $t \in \ci\ce(X)$, we still have $t \in RF^{\circ}_\Gamma(\cw)$.
Conditions (d) is obviously taken care of  by
the construction. In view of Corollary \ref{largerroofsmore}, condition
(a) has been taken care of by the fact that $\cw \fextends
\cw_1 \quo X$.


By conditions (c) and (d), we have:

\begin{assertion}\label{terw}
(i)~$ter[\ce(\cw)\langle X \rangle] \cap V(T) \subseteq Q$.\\
(ii)~$V[\ce(\cw)\langle X \rangle] \cap Q  \subseteq
ter[\ce(\cw)\langle X \rangle]$.
\end{assertion}

\begin{proof}~  Let $t$ be a vertex in $ter[\ce(\cw)\langle X \rangle] \cap V(T)$. By
condition (d) above, $t \in X$. Since by assumption $t \not \in
RF^\circ(\cw)$, by condition (c) it follows that $t \in Q$. This
proves (i).

To prove part (ii), assume that $q \in (Q \cap V[\cw\langle X
\rangle]) \setminus  ter[\ce(\cw)\langle X \rangle]$.
By the self roofing lemma (Lemma
\ref{separatingverticesofwave}), it follows that $q \in
RF^\circ(\cw)$. But, since $\cw$ is a wave in $\Gamma \quo X$, and
$X$ is countable, this contradicts the fact that $q \in Q$.
\end{proof}

Let $\cw'$ be obtained from $\ce(\cw)$ by the removal of all paths
ending at $Q$. By Assertion \ref{terw} (ii), $\cw'$ is a wave in
$(\Gamma\quo X) -Q-a$, and by condition (a) we have $y \in RF_{\Gamma'}(\cw')$.
Thus $\cw'$ has almost all properties required from the
wave $\cu$ in the proposition, the only problem being that we are
looking for a wave $\cu$ in $\Gamma-Q-a$, not in $\Gamma\quo X-Q-a$.
We now wish to ``bring $\cw'$ to the ground'', namely make it start
at $A$, not at $A \cup X$.

To achieve this goal, we enumerate the vertices of $X$ as
$x_1,x_2,\ldots$, and start deleting them one by one - this time,
real deletion, not the quotient operation. Let $k_1=1$, delete
$x_{k_1}=x_1$, and choose a maximal wave $\cv_1$ in $\Gamma-a -x_1$.
Next choose the first vertex $x_{k_2}$ not belonging to $V[\cv_1]$
(if such exists), take a maximal wave $\cv'_2$ in $\Gamma -a
-\{x_{k_1}, x_{k_2}\}$, and define $\cv_2 = \cv_1 \hetz \cv'_2$.
Then choose the first $k_3$ such that $x_{k_3} \not \in V[\cv_2]$
(if such exists), take a maximal wave $\cv'_3$ in $\Gamma -a
-\{x_{k_1}, x_{k_2}, x_{k_3} \}$, and define $\cv_3 = \cv_2 \hetz
\cv'_3$. If the process terminates after $m$ steps for some finite
$m$, let $\cv=\cv_m$. Otherwise, let  $\cv = \uparrow_{k <
\omega}\cv_k$. Let $\theta =\omega$ if this process lasts $\omega$
steps, and $\theta=m+1$ if it terminates after $m$ steps for some
finite number $m$. For $i < \theta$ denote the set $\{x_{k_1},
x_{k_2}, \ldots, x_{k_i} \}$ by $R_i$, and write $R = \{x_{k_1}, x_{k_2}, x_{k_3} \ldots\}$.

Our goal now is to show that $\cv \hetz \cw'$ is a
wave in $\Gamma'$. This will be done by applying
Lemma \ref{grounding} with
 $\Gamma$ replaced by $\Gamma'$, the wave $\cu$ replaced by $\cv$,
the wave $\cw$ replaced by $\cw'$, the set $X$ in the lemma replaced by $R$ and the set $Z$
replaced by $X$.
We already know that $\cw'$ is a
wave in $\Gamma' \quo  X$.
We need to show that $\cv$ is a wave in $\Gamma' - R$ and every path in $\cw'$ meets $RF_{\Gamma'
-R}(\cv)$. Note that we already know that $\cv$ is a wave in $\Gamma - R$.
Therefore, in order to show that it is a wave in $\Gamma' - R$, we only need to prove it does not meet $Q$.
Also note that following Remark \ref{remarkgrounding}, in order to show that every path in $\cw'$ meets $RF_{\Gamma'
-R}(\cv)$, it is enough to consider only paths starting at $X$.

Recall that $\cv'_i$ is a $\le$-maximal wave in $\Gamma -a-R_i$
and by Lemma \ref{hetzbeatsboth}, so is $\cv_i$.

We also have

\begin{assertion}\label{ass1}
$V(T) \cap ter[\cv] = \emptyset$.
\end{assertion}%

\begin{proof} If $t \in V(T) \cap ter[\cv]$ then $t=ter(P)$ for some $P \in
\cv_i$ for some $i$. But then, the wave $\cv_i \setminus \{P\}$ is
a hindrance in $\Gamma - \{a, x_{k_1}, x_{k_2},
\ldots,x_{k_i},t\}$, contradicting the fact that the deletion
of any finite subset of $V(T)$ does not generate a hindrance.
\end{proof}

\begin{assertion}\label{vdisjointfromq}
$V[\cv] \cap Q =\emptyset$.
\end{assertion}

\begin{proof}
Suppose, for contradiction, that $V[\cv] \cap Q \neq \emptyset$.
Then there exists $i <\theta$ and $q \in Q$ such that $q \in
V[\cv_i]$. By Assertion \ref{ass1}, $q \not \in ter[\cv_i]$, and
since $\cv_i$ is a wave in $\Gamma- a - R_i$, by the self roofing lemma (Lemma
\ref{separatingverticesofwave}) $q \in RF^\circ_{\Gamma- a -
R_i}(\cv_i)$. By Lemma \ref{quotientisstronger} it follows that
 $q \in RF^\circ_\Gamma(\cu)$, where $\cu$
is a maximal wave in $(\Gamma-a) \quo R_i$. But this contradicts
the definition of $Q$.
\end{proof}

\begin{remark} As pointed out by R. Diestel,
Assertion \ref{vdisjointfromq} is not essential for the argument
that follows, since by the definition of $Q$ we have: $V[\cv] \cap
Q \subseteq ter[\cv]$. Thus we could replace $\cv$ by
$\cv'=\cv\setminus \cv \langle Q \rangle$, and the argument below
would remain valid. But since in fact $\cv'=\cv$, we chose the
longer, but more informative, route.
\end{remark}

By Assertion
\ref{vdisjointfromq} $\cv$ is a wave in $\Gamma' - R$.

\begin{assertion}
\label{cvroofy} If $z \in Y \cap V[\cw\langle X \rangle]$ then $z
\in RF_{\Gamma' - R}(\cv)$.
\end{assertion}

\begin{proof}
By (b) we have $F_z \subseteq X$. Let $n < \omega$ be chosen so that
$X' = \{x_1, \ldots ,x_n\} \supseteq F_z$. Since $\Gamma -a- X'$ is
unhindered and $\Gamma -a- X' - z$ is hindered, by Lemma
\ref{hinderedimplieswave} there exits a wave $\cz$ in $\Gamma -a-
X'$ with $z \in ter[\cz]$. Let $i$ be maximal such that $R_i
\subseteq X'$.  By the maximality property of $i$ we have $X'
\setminus R_i \subseteq V[\cv_i] \subseteq RF_{\Gamma-a-R_i}(\cv_i)$.

We now note that $\cv_i$ is a wave in $\Gamma-a-R_i$ and $\cz$
is a wave in $\Gamma - a -X'$. Hence we can conclude that $\cv_i \hetz \cz$ is a wave in
$\Gamma - a - R_i$, by applying
 Lemma \ref{roofminus} (with $\Gamma$ replaced by
$\Gamma-a-R_i$, the wave $\cu$ replaced by $\cv_i$, the wave $\cw$
replaced by $\cz$ and $Y$ replaced by $X' \setminus R_i$).
By the maximality of $\cv_i$ we have $\cv_i=\cv_i \hetz \cz$. This
implies that $z \in RF_{\Gamma -a- R_i}(\cv_i)$. Since $\cv \extends
\cv_i$ and $R \cup Q \supseteq R_i$ we have $z \in RF_{\Gamma - a -
Q - R}(\cv)$.
\end{proof}

\begin{corollary}\label{wmeetsrfv}
Every path in $\cw'$ meets $RF_{\Gamma'-R}(\cv)$.
\end{corollary}

\begin{proof}
Let $W$ be a path in $\cw'$ and let $w = in(W)$. Then either
$w \in A$ or $w \in X$. If $w \in A$ then since $\cv$ is a wave in
$\Gamma' - R$, we have $w \in RF_{\Gamma'-R}(\cv)$. If $w \in X$
then let $t = ter(W)$.
Since $\cw'$ was obtained from $\cw$ by removing paths ending at $Q$, we have $t \not \in Q$.
By Assertion \ref{terw}(i), we now have $t \not \in V(T)$.
Let $z$ be the first vertex in $W$ outside $V(T)$.
Then $z \in Y$ and by assertion \ref{cvroofy} we have $z \in F_{\Gamma'-R}(\cv)$.
\end{proof}

Define: $\cu_y= \cv \hetz \cw'$. Apply Lemma \ref{grounding} with
 $\Gamma$ replaced by $\Gamma'$, the wave $\cu$ replaced by $\cv$,
the wave $\cw$ replaced by $\cw'$, the set $X$ in the lemma replaced by $R$ and the set $Z$
replaced by $X$.
Corollary \ref{wmeetsrfv} asserts that indeed
every path in $\cw'$ meets $RF_{\Gamma'-R}(\cv)$
as needed to apply the lemma.
The lemma yields that the warp
$\cu_y$ is a wave in $\Gamma'$. This completes the proof of
Proposition \ref{mainassertion}.

To end the proof of Theorem \ref{safelinking}, let $\cu =
\uparrow_{y\in Y}\cu_y$. Then $\cu$ separates  $Y$ from $B$. By
Lemma \ref{separatingoutneighbors} it follows that $\cu$ is a wave
in $\Gamma$, and since it does not contain $a$ as an initial
vertex of a path, it is a hindrance in $\Gamma$. This contradicts
the assumption that $\Gamma$ is unhindered. $\enp$

\section{$\kappa$-ladders and $\kappa$-hindrances}
\subsection{Stationary sets}

As is customary in set theory, an ordinal is taken as the set of
ordinals smaller than itself, and a cardinal $\kappa$ is
identified with the smallest ordinal of cardinality $\kappa$. An
uncountable cardinal $\lambda$ is called {\em singular} if there
exists a sequence $(\nu_\alpha:~\alpha <\mu)$ of ordinals, whose
limit is $\lambda$, where all $\nu_\alpha$, as well as $\mu$, are
smaller than $\lambda$. The smallest singular cardinal is
$\aleph_\omega$, which is the limit of $(\aleph_i:~i <\omega)$. A
singular cardinal is necessarily a limit cardinal, namely it must
be of the form $\aleph_\theta$ for some limit ordinal $\theta$. On
the other hand,  ZFC (assuming its consistency) has models in
which there exist non-singular limit cardinals.

A non-singular cardinal is called {\em regular}.

The main set-theoretic notion we shall use is that of {\em
stationary sets}. A subset of an uncountable regular cardinal
$\kappa$ is called {\em unbounded} if its supremum is $\kappa$,
and {\em closed} if it contains the supremum of
each of its bounded subsets. A subset of $\kappa$ is called {\em stationary} (or
$\kappa$-stationary) if it intersects every closed unbounded
subset of $\kappa$. A function $f$ from a set of ordinals to the
ordinals is called {\em regressive} if $f(\alpha)<\alpha$ for all
$\alpha$ in the domain of $f$. A basic fact about stationary sets
is Fodor's lemma:

\begin{theorem}\label{fodor}
If $\kappa$ is regular and uncountable, $\Phi$ is a
$\kappa$-stationary set, and $f:~\Phi \to \kappa$ is regressive,
then there exist a stationary subset $\Phi'$ of $\Phi$ and an
ordinal $\beta$ such that $f(\phi)=\beta$ for all $\phi \in
\Phi'$.
\end{theorem}

 Fodor's lemma implies that stationary sets are in some
sense ``big". This is expressed also in the following:

\begin{lemma}\label{stationaryisbig} If $\Xi_\alpha,~\alpha <\lambda$ are non-stationary, and
$\lambda < \kappa$, then $\bigcup_{\alpha<\lambda} \Xi_\alpha$ is
non-stationary.
\end{lemma}

This is another way of saying that the intersection of fewer than
$\kappa$ closed unbounded sets is closed and unbounded.

\subsection{$\kappa$-ladders}

The tool used in the proof of Theorem \ref{halltypetheorem} in the
uncountable case is $\kappa$-{\em ladders}, for uncountable
regular cardinals $\kappa$. A $\kappa$-ladder $\cl$ is a sequence
of ``rungs" $(R_\alpha:~\alpha< \kappa)$. At each step $\alpha$ we
are assuming that a warp $\cy_\alpha=\cy_\alpha(\cl)$ in $\Gamma$
is defined, by the previous rungs of $\cl$.  For each $\alpha \ge
0$, assuming $\cy_\alpha$ is defined, we let
 $\Gamma_\alpha=\ce(\Gamma\quo \cy_\alpha)$.

The warp $\cy_0$ is defined as $\langle A \rangle$, and for limit
ordinals $\alpha$, we let $\cy_\alpha=\uparrow_{\theta <
\alpha}\cy_\theta$.

For successor ordinals $\alpha+1$, the warp $\cy_{\alpha+1}$ is
defined by $\cy_\alpha$ and by the rung $R_\alpha$, the latter
being chosen as follows. A first constituent of $R_\alpha$ is a
(possibly trivial) wave $\cw_\alpha$ in $\Gamma_\alpha$. If the
set $V(\Gamma_\alpha)\setminus (A(\Gamma_\alpha)\cup
V[\cw_\alpha])$ is non-empty, then $R_\alpha$ consists also of a
vertex $y_\alpha$ in this set. The warp $\cy_{\alpha+1}$ is
defined in this case as $\cy_\alpha \hetz \cw_\alpha \cup \langle
y_\alpha \rangle$. If $V(\Gamma_\alpha)\setminus
(A(\Gamma_\alpha)\cup V[\cw_\alpha])=\emptyset$, then
$\cy_{\alpha+1}$ is defined as $\cy_\alpha \hetz \cw_\alpha$. In
this case all consecutive rungs will consist just of the trivial
wave, meaning that the ladder will ``mark time", without changing.

We also wish to keep track of the steps in which a new hindrance
emerges in the ladder. This is done by keeping record of subsets
$\ch_\alpha$ of $\cy_\alpha$. These sets are not uniquely defined
by $\cl$, but to simplify notation we assume that the ladder comes
with a fixed choice of such sets, which is subject to the
following conditions.

We define $\ch_0=\emptyset$. If $\ci\ce(\cy_{\alpha+1}) \setminus
\ch_\alpha \neq \emptyset$ we pick  a (possibly unending) path $H$
in this set, write $H_\alpha=H$, and $\ch_{\alpha+1} =
\bigcup_{\theta < \alpha}\ch_\theta \cup \{H\}$.

If $\ci\ce(\cy_{\alpha+1}) \setminus \ch_\alpha \neq \emptyset$
we let $\ch_{\alpha+1}=\ch_\alpha$. For limit $\alpha$ we define
$\ch_\alpha = \bigcup_{\theta < \alpha}\ch_\theta$.

\begin{remark} Note that it is possible that $\bigcup_{\alpha < \kappa}
\ch_\alpha \neq \ci\ce(\cy)$, namely that we never exhaust all of
$\ci\ce(\cy)$. \end{remark}

Since a path in $\ch_\alpha$ is inessential in $\cy_\alpha$, it
will never ``grow" in any later stage $\beta$, and hence we have:

\begin{lemma}\label{halphainessential}
$\ch_\alpha \subseteq \ci\ce(\cy_\beta)$ for all $\beta \ge
\alpha$.
\end{lemma}

The set of ordinals $\alpha$ for which $\ci\ce(\cy_{\alpha+1})
\setminus \ch_\alpha \neq \emptyset$ is denoted by $\Phi(\cl)$. As
noted, $\Phi(\cl)$ is not uniquely defined by $\cl$ itself, and is
dependent on the choice of the sets $\ch_\alpha$.

\begin{example}\label{aisall}
Let $|A|=\aleph_0,~B=\emptyset, ~V(\Gamma)=A$, and choose $\kappa
= \aleph_1$. Since  $\Gamma_1$ is defined as $\ce(\Gamma \quo
\langle A \rangle)$, it is empty (i.e., $\Gamma_1$  has no
vertices), and $\cy_\alpha = \ci\ce(\cy_\alpha)=\langle A \rangle$
for all $1 \le \alpha < \aleph_1$. The  paths $(a),~ a \in  A$ can
be chosen as $H_\alpha$ in any order, and thus $\Phi(\cl)$ can be
any countable ordinal.
\end{example}

 We write $\Phi^\infty(\cl)$ for the set of
those $\alpha \in \Phi(\cl)$ for which $\ci\ce(\cy_{\alpha+1})
\setminus \ch_\alpha$ contains an unending path, and $\Phi^{fin}$
for $\Phi(\cl) \setminus \Phi^\infty(\cl)$.

Let $\Phi_h(\cl)=\{\alpha \mid \cw_\alpha ~\mbox{is a
hindrance}\}$, and $\Phi_h^\infty(\cl)= \{\alpha \mid
\cy_\alpha^\infty \setminus \bigcup_{\theta<\alpha}
\cy_\theta^\infty \neq \emptyset\}$. Unlike $\Phi(\cl)$, the set
$\Phi_h(\cl)$ is determined by $\cl$. The difference between the
two sets is that the ordinals in $\Phi_h(\cl)$ are ``newly
hindered", namely there is a hindered vertex generated at that
stage, whereas the fact that $\alpha \in \Phi(\cl)$ means that not
all hindered vertices generated so far have been ``taken into
account", in the sense of being included in $\ch_\alpha$. In
Example \ref{aisall} $\Phi_h(\cl)=\{0\}$.

\begin{lemma}\label{hinderedimplieshindered}
$\Phi_h(\cl) \subseteq \Phi(\cl)$.
\end{lemma}

\begin{proof}
Suppose that $\alpha \in \Phi_h(\cl)$. We shall show that
$\ci\ce(\cy_{\alpha+1}) \setminus \ch_\alpha \neq \emptyset$,
which will imply the desired inclusion result. Let $x$ be a vertex
in $A(\Gamma_\alpha) \setminus in[\cw_\alpha]$. Then $x = ter(P)$
for some $P \in \ce(\cy_\alpha)$. By the definition of
$\ch_\alpha$, we have $P \not \in \ch_\alpha$. By the definition
of a wave, $ter[\cw_\alpha]$ is separating in $\Gamma_\alpha$ and
thus also in $\Gamma$. The set $ter[\cy_\alpha \hetz \cw_\alpha]
\setminus \{x\}$ contains $ter[\cw_\alpha]$ and is hence
separating as well. Therefore $P \in \ci\ce(\cy_\alpha \hetz
\cw_\alpha)$. Thus $\ci\ce(\cy_\alpha \hetz \cw_\alpha) \setminus
\ch_\alpha \neq \emptyset$, meaning that $R_\alpha$ is hindered.
\end{proof}

\begin{lemma}\label{phiinftycontainedinphi}
$\Phi_h^\infty(\cl) \subseteq \Phi^\infty(\cl)$.
\end{lemma}

\begin{proof} Let $\alpha$ be an ordinal in $\Phi_h^\infty(\cl)$,
and let $P$ be a path witnessing this, namely $P \in
\cy_\alpha^\infty \setminus \bigcup_{\theta<\alpha}
\cy_\theta^\infty$. Then $P \not \in \bigcup_{\theta<\alpha}
\cy_\theta$, and since $\ch_\alpha \subseteq
\bigcup_{\theta<\alpha} \ci\ce(\cy_\theta)$, this implies that $P
\in \ci\ce(\cy_\alpha) \setminus \ch_\alpha$.
\end{proof}

The following is obvious from the way the sets $\ch_\alpha$ are
chosen:

\begin{lemma}\label{fewhindered}
If~\quad $|\ci\ce(\cy_\alpha)|\ge \kappa$ for some $\alpha
<\kappa$, then $\Phi(\cl) \supseteq [\alpha, \kappa)$.
\end{lemma}


\begin{notation}\label{yalphaandxalpha}
Write $T_\alpha=T_\alpha(\cl)$ for $A(\Gamma_\alpha)$. The warp
$\cy_\kappa$ is denoted by $\cy=\cy(\cl)$.
 For
$\alpha\in \Phi^{fin}(\cl)$ denote $ter(H_\alpha)$ by $x_\alpha$.
The set $\{y_\alpha:~\alpha < \kappa\}$ is denoted by $Y(\cl)$,
and for every $\beta \le \kappa$ write $Y_\beta(\cl)$ for
$\{y_\alpha:~\alpha < \beta\}$. The set $\{x_\alpha \mid \alpha
\in \Phi^{fin}(\cl)\}$ is denoted by $X^{fin}(\cl)$.
\end{notation}

The definitions clearly imply:

\begin{lemma}\label{t_alpha_progressive}
$T_\alpha$ is $A$--$B$-separating for all $\alpha <\kappa$. If
$\alpha < \beta$ then $T_\alpha \subseteq RF(T_\beta)$.
\end{lemma}

By the definition of $\Gamma_\alpha$ as $\ce(\Gamma \quo
\cy_\alpha)$ we have:

\begin{lemma}\label{talpha_minimal} $T_\alpha$ is a minimal  $A$--$B$-separating for all $\alpha
<\kappa$.
\end{lemma}

Define $RF(\cl) = \bigcup_{\theta < \kappa} RF(T_\theta)$ and
$RF^\circ(\cl) = \bigcup_{\theta < \kappa} RF^\circ(T_\theta)$.

 Also write
$\Gamma^\alpha=\Gamma[RF(T_\alpha)]$, which means that
$D(\Gamma^\alpha)$ (the digraph of $\Gamma^\alpha$) is
$\Gamma[RF(T_\alpha)]$, $A(\Gamma^\alpha)=A$ and
$B(\Gamma^\alpha)=T_\alpha$.

For $\alpha < \beta$ let $\Gamma_\alpha^\beta$ be the part of
$\Gamma$ between $T_\alpha$ and $T_\beta$, namely
$V(\Gamma_\alpha^\beta)=V(\Gamma_\alpha[RF_{\Gamma_\alpha}(T_\beta)])$,
$D(\Gamma_\alpha^\beta)=D(\Gamma_\alpha[RF_{\Gamma_\alpha}(T_\beta)])$,
$A(\Gamma_\alpha^\beta)=T_\alpha,
~B(\Gamma_\alpha^\beta)=T_\beta$.

\begin{notation}
We shall write $V^\alpha=V^\alpha(\cl)$ for $V(\Gamma^\alpha)$,
and $V_\alpha$ for $V(\Gamma_\alpha)$, namely
$V^\alpha=RF(T_\alpha)$ and $V_\alpha=V(\Gamma) \setminus
RF^\circ(T_\alpha)$.
\end{notation}

\begin{notation}\label{xandphi}
 Let
$\Phi_G(\cl)=\{\alpha \in \Phi(\cl) \mid in(H_\alpha) \in A\}$ and
$\Phi_H(\cl) = \Phi(\cl) \setminus \Phi_G(\cl)$ (The ``G" stands
for ``grounded" and the ``H" stands for ``hanging in air").
\end{notation}

Throughout the proof we shall construct again and again ladders,
which will all be denoted by $\cl$. In all these cases we shall
use the following:

\begin{convention}\label{notationonladders}
We shall denote $\cy(\cl)$, for the ladder $\cl$ considered at
that point, by $\cy$.  We shall also write $T_\alpha$ for
$T_\alpha(\cl)$,  $Y_\alpha$ for  $Y_\alpha(\cl)$, $\Phi$ for
$\Phi(\cl)$, and so on.
\end{convention}

\begin{lemma}\label{onlyaxcounts}
$\Phi_H(\cl)$  is non-stationary.
\end{lemma}

\begin{proof}
For  $\alpha \in \Phi_H(\cl)$ we have $in(H_\alpha) =y_\beta$ for
some $\beta <\alpha$. The function $f(\alpha)=\beta$ defined in
this way is a regressive injection from $\Phi_H(\cl)$ to $\kappa$.
Thus, by Fodor's lemma, $\Phi_H(\cl)$ is not stationary.
\end{proof}

The following  is obvious:

\begin{lemma}
\label{slambdaplus} A vertex $v \in V$ belongs to $RF(\cl)
\setminus RF^\circ(\cl)$ if and only if there exists $\beta<
\kappa$ such that $v \in T_\alpha$ for all $\alpha \ge \beta$.
\end{lemma}

\begin{lemma}\label{stayingwithinroofedtalpha}
Let $Q$ be a $\cy$-alternating path, and assume that  $in(Q) \in
RF^\circ(T_\alpha)$. Then: \begin{enumerate}
\item \label{staying}
 $V(Q) \subseteq RF(T_\alpha)$,~ and:
\item \label{betasmaller}
If $in(Q)=x_\alpha$ and $ter(Q)=y_\beta$, then $\beta< \alpha$.
\end{enumerate}
\end{lemma}

\begin{proof} Write
 $z=in(Q)$. Using the same notation as in Definition
\ref{alternatingpaths}, write $Q$ as
$(z=z_0,F_1,u_1,R_1,z_1,F_2,u_2,R_2,z_2...)$, where $F_i$ are
forward paths, namely using edges not belonging to $E[\cy]$, $R_i$
are backward paths, namely using edges of $E[\cy]$, $u_i$ are
vertices on paths from $\cy$ at which $Q$ switches from forward to
backward direction, and $z_i$ vertices at which $Q$ switches from
backward to forward direction. Since $z \in RF^\circ(T_\alpha)$,
and $T_\alpha$ separates $V[\cl]$ from $B$, $F_1$ is contained in
$RF(T_\alpha)$. Possibly $u_1 \in T_\alpha$, but since $R_1$ goes
backwards, $z_1 \in RF^\circ(T_\alpha)$. Thus $F_2$ is contained
in $RF(T_\alpha)$. By an inductive argument following these steps
we obtain part \ref{staying} of the lemma.

If $ter(Q)=y_\beta$, then by part (\ref{staying}), $y_\beta \in
RF(T_\alpha)$. But $y_\beta \in V(\Gamma_\beta) \setminus
A(\Gamma_\beta) = V(\Gamma) \setminus RF(T_\beta)$. Therefore
$RF(T_\alpha) \setminus RF(T_\beta) \neq \emptyset$, and hence
$\beta <\alpha$.
\end{proof}

Write $\zeta(\alpha)$ for the minimal ordinal at which $H_\alpha$
emerges as an inessential path, namely the minimal ordinal $\beta$
such that $H_\alpha \in \ci\ce(\cy_\beta)$. The choice of
$H_\alpha$ implies:

\begin{lemma}\label{appear_in_time}
$\zeta(\alpha)\le \alpha$ for all $\alpha \in \Phi(\cl)$.
\end{lemma}

Since $H_\alpha \in \ci\ce(\cy_{\zeta(\alpha)})$, we have:

\begin{lemma}
$x_\alpha \in RF^\circ(T_{\zeta(\alpha)})$ for every $\alpha \in
\Phi^{fin}(\cl)$.
\end{lemma}

Combined with Lemma \ref{appear_in_time}, this yields:

\begin{lemma}\label{x_alpha_in_rfcircle_alpha}
$x_\alpha \in RF^\circ(T_\alpha)$ for every $\alpha \in
\Phi^{fin}(\cl)$.
\end{lemma}

\subsection{$\kappa$-hindrances}

Ordinals in $\Phi(\cl)$ are ``troublesome", witnessing as they do
the existence of hindrances. Thus, if $\Phi(\cl)$ is ``large" then
the ladder may pose a problem for linkability of $\Gamma$. And now
we know what ``large" should be: stationary. This is the origin of
the following definition:

\begin{definition}
If $\Phi(\cl)$ is $\kappa$-stationary, then  $\cl$ is called a
$\kappa$-{\em hindrance}.
\end{definition}

Lemmas \ref{onlyaxcounts} and  \ref{stationaryisbig} yield
together:

\begin{lemma}\label{phiastationary}
If $\cl$ is a $\kappa$-hindrance then $\Phi_G(\cl)$ is stationary.
\end{lemma}

\begin{example} Let $A$ be a set of size $\aleph_1$, $B$ a set of
size $\aleph_0$, let $D$ be the complete directed graph on
$(A,B)$, namely  $E(D)=A \times B$, and let $\Gamma=(D,A,B)$. We
define an $\aleph_1$ ladder in $\Gamma$, as follows. Order $B$ as
$(b_\alpha \mid \alpha < \omega)$ and $A$ as
$(a_\alpha \mid \alpha < \omega_1)$.

For $\alpha < \omega$ let $\cw_\alpha$ be the trivial wave, and
$y_\alpha =b_\alpha$. Then for all such $\alpha$ we have
$\Gamma_\alpha = \Gamma \quo \{b_i \mid i < \alpha\}$ and
$\ch_\alpha = \emptyset$. At the $\omega$ step we have $\cy_\omega
= \langle A \cup B \rangle$, $\Gamma_\omega = \Gamma \quo  B =
((B,\emptyset),B,B)$ and $\ch_\omega = \emptyset$. Note that all
the singleton paths in $\langle A \rangle$ are inessential in
$\cy_\omega$.

For $0\le \alpha \le \aleph_1$ let $R_{\omega +\alpha}$ consist of
the inessential singleton path $H_{\omega +\alpha}=(a_\alpha)$. We
then  have $ \cy_{\omega+\alpha} = \langle A \cup B \rangle$,
 $\Gamma_{\omega+\alpha} =
((B,\emptyset),B,B)$ and $\ch_{\omega+\alpha} = \langle \{a_\theta
\mid \theta < \alpha\} \rangle$.

Thus $\Phi(\cl)= [\omega, \aleph_1)$, which is stationary, and
hence $\cl$ is an $\aleph_1$-hindrance.
\end{example}

\begin{example} [accommodated from \cite{aharoninwshelah}] Let $\kappa$ be an
uncountable regular cardinal, and $\Psi$  a $\kappa$-stationary
set. Let $A = \{a_\alpha \mid \alpha \in \Psi\}$, $B = \{b_\alpha
\mid \alpha < \kappa\}$, and let $D$ be the directed graph whose
vertex set is $A \cup B$ and whose edge set is $E = \{(a_\alpha,
b_\beta) \mid \beta < \alpha \}$. Let $\Gamma=(D,A,B)$.

By Fodor's lemma, $\Gamma$ is unlinkable.

Define a $\kappa$-ladder in $\Gamma$ as follows. For all $\alpha<
\kappa$ let $y_\alpha = b_\alpha$ and let $\cw_\alpha$ be the
trivial wave. Define sets $\ch_\alpha$ by adding to $\ch_\alpha$,
for each $\alpha \in \Psi$, the singleton inessential path
$H_\alpha = (a_\alpha)$. Here we have $\cy_\alpha = \langle A \cup
\{b_\theta \mid \theta < \alpha \} \rangle$ and the path
$(a_\beta)$ is inessential in it for every $\beta \leq \alpha$.
Since $\Psi$ is stationary, this is a $\kappa$-hindrance.
\end{example}

\begin{example}
The following example shows the role of infinite paths in
$\kappa$-hindrances. Let $\Psi$ be an $\aleph_1$-stationary set
all of whose element are limit ordinals (e.g., $\Psi$ can be the
set of {\em all} countable limit ordinals). For every $\alpha \in \Psi$,
let $(\eta_i^\alpha \mid i < \omega)$
be an ascending
sequence converging to $\alpha$, where $\eta_0^\alpha = 0$.

Let
$C = \{c_i^\alpha \mid \alpha \in \Psi, ~ i < \omega\}$,
$B = \{b_\alpha : \alpha < \omega_1\}$,
let $A$ be the subset of $C$
$A = \{c_0^\alpha \mid \alpha \in \Psi \}$,
let $D$ be the directed graph whose vertices
are $C \cup B$ and whose edges are
$E = \{(c_i^\alpha, c_{i+1}^\alpha)\mid \alpha \in \Psi, ~ i < \omega\}
\cup \{(c_i^\alpha, c_j^\beta) \mid \alpha, \beta \in \Psi, ~ i,j < \omega, ~
\beta < \alpha, ~ \eta_i^\alpha \geq \eta_j^\beta \}
\cup \{(c_i^\alpha, b_\beta) \mid \alpha \in \Psi, ~ i < \omega, ~
\beta \leq \eta_i^\alpha \}$
and and let $\Gamma=(D,A,B)$.

Again, by Fodor's lemma, $\Gamma$ is unlinkable.

We can construct an $\aleph_1$-ladder $\cl$ on $\Gamma$ by taking
$y_\alpha = b_\alpha$ and $\cw_\alpha =
\{(b_\beta) \mid \beta < \alpha \} \cup
\{(c_i^\beta, c_{i+1}^\beta) \mid \eta_{i+1}^\beta = \alpha \}
\cup \{(c_i^\beta) \mid \eta_i^\beta < \alpha < \eta_{i+1}^\beta \}$.
For $\alpha \in \Psi$, the concatenation of these waves forms an infinite
path $(c_0^\alpha, c_1^\alpha, c_2^\alpha, c_3^\alpha, \ldots)$
in $\cy_\alpha$. We can take this path as $H_\alpha$.

This yields $\Phi(\cl) = \Psi$ and therefore $\cl$ is an
$\aleph_1$-hindrance.

\end{example}

\begin{lemma}\label{most_y_reach_salpha}
If $\Gamma$ does not contain a $\kappa$-hindrance then for every
$\kappa$-ladder $\cl$ and every $\alpha < \kappa$ there holds
$|\cy_\alpha \langle \sim T_\alpha \rangle| < \kappa$.
\end{lemma}

\begin{proof}
A path $P \in \cy_\alpha$ not meeting $T_\alpha$ belongs to
$\ci\ce(\cy_\alpha)$. Hence, if $|\cy_\alpha  \langle  \sim T_\alpha
\rangle| \ge \kappa$ then $|\ci\ce(\cy_\alpha)| \ge \kappa$, and
hence by Lemma \ref{fewhindered} $\cl$ is a $\kappa$-hindrance.
\end{proof}

The following lemma is not essential for the discussion to follow,
but its understanding may clarify the nature of
$\kappa$-hindrances. It says that Lemmas
\ref{hinderedimplieshindered}, \ref{phiinftycontainedinphi} and
\ref{fewhindered} summarize all reasons for $\cl$ to be a
$\kappa$-hindrance:

\begin{lemma}\label{dependingonlyonl}
A $\kappa$-ladder $\cl$ is a $\kappa$-hindrance if and only if
either:

(i)~$\Phi_h(\cl) \cup \Phi_h^\infty(\cl)$ is stationary, ~or:

(ii)~$|\ci\ce(\cy_\alpha)|\ge \kappa$ for some $\alpha <\kappa$.
\end{lemma}

This means, among other things, that although $\Phi(\cl)$ is not
uniquely determined by $\cl$, whether it is stationary or not is
determined by $\cl$ alone. Namely, $\cl$ being a
$\kappa$-hindrance is independent of the order by which the paths
$H_\alpha$ are chosen. The lemma also clarifies why we need to
work with $\Phi(\cl)$ rather than $\Phi_h(\cl)$: because of the
possible occurrence of case (ii).

{\em Proof of Lemma \ref{dependingonlyonl}}: In view of Lemmas
\ref{hinderedimplieshindered}, \ref{phiinftycontainedinphi} and
\ref{fewhindered}, it remains to be shown that if $\Phi(\cl)$ is
stationary, then one of conditions (i) and (ii) is true. By Lemma
\ref{appear_in_time} $\zeta(\alpha) \le \alpha$ for all $\alpha$.
If the set $\{\alpha \mid \zeta(\alpha)=\alpha\}$ is stationary,
then (i) holds. Otherwise, assuming  $\Phi(\cl)$ is stationary, by
Fodor's lemma there exist a stationary subset $\Phi' \subseteq
\Phi(\cl)$ and an ordinal $\beta<\kappa$, such that
$\zeta(\alpha)=\beta$ for every $\alpha \in \Phi'$. By the
definition of $\zeta$ this implies that $|\ci\ce(\cy_\beta)|\ge
\kappa$, proving (ii).
 $\enp$




\begin{lemma}\label{noinfinitepathsiny}
Let $\cl$ be a $\kappa$-ladder that is not a $\kappa$-hindrance,
and let $\Sigma$ be a closed unbounded set avoiding $\Phi(\cl)$.
Then for every $P \in \cy(\cl)$ the set $\Sigma(P)=\{\alpha \in
\Sigma \mid T_\alpha \cap V(P) \neq \emptyset\}$ is closed in
$\kappa$.
\end{lemma}

\begin{proof}

Let $\Psi$ be an infinite subset of $\Sigma(P)$, and assume, for
contradiction, that $\beta =\sup \Psi$ does not belong to
$\Sigma(P)$, namely $V(P) \cap T_\beta =\emptyset$. By assumption,
$T_\alpha \cap V(P) \neq \emptyset$ for some $\alpha < \beta$.
Choose a vertex $x \in T_\alpha \cap V(P)$. Since $\beta \not \in
\Sigma(P)$, we have $x \not \in T_\beta$, and thus $x \in
RF^\circ(T_\beta)$, which together with the assumption that
$\Sigma(P) \cap T_\beta =\emptyset$ implies that $V(P) \subseteq
RF^\circ(T_\beta)$, meaning that $P \in \ci\ce(\cy_\beta)$. Since
$V(P) \cap T_\psi \neq \emptyset$ for every $\psi \in \Psi$, for
each such $\psi$ there exists an initial segments of $P$ belonging
to $\ce(\cy_\psi)$. But this clearly implies that $P \not\in
\bigcup_{\psi \in \Psi}\ci\ce(\cy_\psi)$, and thus $\beta \in
\Phi_h(\cl)$, contradicting the fact that $\Phi(\cl) \cap
\Sigma=\emptyset$.
\end{proof}

Theorem \ref{halltype} will follow from the combination of two theorems:

 \begin{theorem}\label{C}
If $\Gamma$ does not possess a hindrance or a $\kappa$-hindrance
for any uncountable regular cardinal $\kappa$, then it is
linkable.
\end{theorem}

\begin{theorem}\label{kappahindranceimplieshindrance}
If $\Gamma$ contains a $\kappa$-hindrance  for some uncountable
regular cardinal $\kappa$, then it contains a hindrance.
\end{theorem}

Theorem \ref{C} is akin to a version of the infinite ``marriage
theorem", proved in \cite{aharoninwshelah}, hence an appropriate
name for it is ``the linkability theorem". We shall prove Theorem
\ref{kappahindranceimplieshindrance} in the next section, and
Theorem \ref{C} in the last section of the paper.

\section{From $\kappa$-hindrances to hindrances}

In this section we prove Theorem
\ref{kappahindranceimplieshindrance}. Namely, that if $\Gamma$
contains a $\kappa$-hindrance for some uncountable regular
cardinal $\kappa$, then it is hindered. This was, in fact, proved
in \cite{aharoni97}. The proof there is only for
$\kappa=\aleph_1$, but it goes verbatim to all uncountable regular
cardinals $\kappa$. That proof is shorter than the one given
below, since  it relies on previous results. It uses the bipartite
conversion, applies the bipartite version of Theorem
\ref{kappahindranceimplieshindrance} proved in
\cite{aharonikonig}, and shows how to take care of the one problem
that may arise along this route, namely that the paths in the
resulting hindrance are non-starting.

Our proof here does not use the main result of
\cite{aharonikonig}, but rather re-proves it, borrowing as ``black
boxes" only two lemmas. We use this as an opportunity to give the
main theorem of \cite{aharonikonig}  a more transparent proof, in
that  its main idea is summarized in a separate theorem (Theorem
\ref{coreofkapimplieshind} below). Another advantage of the
present proof is that one can see what is happening in the graph
itself, rather than in the bipartite conversion.

The basic notion in the proof of the theorem is that of
popularity of vertices in a hindrance. A vertex is ``popular" if it has a
large  in-fan of $\cy$-alternating paths, where $\cy$ is the
warp appearing in the hindrance, and ``large" means
reaching ``stationarily many" points $x_\alpha$.
Let us first illustrate this idea in a very simple case - the
simplest type of unlinkable webs:

\begin{theorem}\label{cardinalities}
A bipartite web $(D,A,B)$ in which $|A| > |B|$ contains a
hindrance.
\end{theorem}

\begin{proof}
 The argument is easy when $B$ is finite, so assume that $B$ is
infinite, and write $|B|=\kappa$. Call a vertex $b \in B$ {\em
popular} if $|N(b)| > \kappa$. Let $U$ be the set of unpopular
elements of $B$. Then $|N(U)| \le \kappa$, and hence in the web
$(D-U-N(U),A \setminus N(U),B \setminus U)$ every vertex in $B
\setminus U$ is of degree larger than $\kappa$, while of course
$|B \setminus U| \le \kappa$. Hence there exists a matching $F$ of
$B \setminus U$ properly into $A \setminus N(U)$. The warp $F \cup
\{(a) \mid ~a \in N(U)\}$ is then a hindrance in $\Delta$.
\end{proof}

 Next we introduce a more general type of unlinkable webs:

\begin{definition}\label{unbalanced}
A web $(G,X,Y)$ is called $\kappa$-{\em unbalanced} if there exist
a function
 $f:~X \to \kappa$  and an injection $g:~Y \to \kappa$,  such that:

\begin{enumerate}
\item\label{xstationary}
$f[X]$ is $\kappa$-stationary.
\item \label{regressive}
$f(in(P))>g(ter(P))$ for every  $X$--$Y$-path $P$.

\end{enumerate}
\end{definition}

This is an ordinal version of the notion of a web in which the
source side has larger cardinality than the destination side. And
indeed, from Fodor's lemma there follows:

\begin{lemma}\label{unbalancedunlinkable}
A $\kappa$-unbalanced web is unlinkable. In fact, for every
$X$--$Y$-warp $\cw$, $f[in[\cw]]$ is non-stationary.
\end{lemma}

In particular, $f[X \cap Y]$ is non-stationary.

The core of the proof of Theorem
\ref{kappahindranceimplieshindrance} is in showing that
$\kappa$-unbalanced webs are hindered, which is of course a
special case of our main theorem, Theorem \ref{halltypetheorem}.
But we shall need a bit more.

Given such a web, a set $S$ of vertices is called {\em popular} if
either $S \cap X \neq \emptyset$, or there exists an $S$-joined
family of  $X$-$S$-paths $\cp$, such that $f[in[\cp]]$ is
$\kappa$-stationary. It is called {\em strongly popular} if there
exists an  $X$-$S$-warp $\cp$, such that $f[in[\cp]]$ is
$\kappa$-stationary (in particular, if $f[X \cap S]$ is
stationary). A vertex $v$ is called ``popular" if $\{v\}$ is
popular.

\begin{theorem}\label{coreofkapimplieshind}
Let $\Lambda=(G,X,Y)$ be a $\kappa$-unbalanced web, with $f$ and
$g$ as above. Then there exists an $X$--$Y$-separating set $S$
such that:
\begin{enumerate}
\item\label{spopular}~
Every vertex $s$ of $S$ is popular in
$\Lambda[RF^\circ(S)\cup \{s\}]$,
i.e., either $s \in X$ or there exits an $X$-starting $s$-in-fan $\cp$
in $G[RF^\circ(S)\cup \{s\}]$,
where $f[in[\cp]]$
is stationary.
\item \label{snotstronglypopular}
$S$ is not strongly popular.

\item\label{ssmall}~
$|S\setminus X|\le \kappa$.

\end{enumerate}
\end{theorem}

For the proof we shall need two results from \cite{aharonikonig}:

\begin{lemma}\label{choicefromnonstationary}
If ~ $\Xi_u,~u \in U$ are non-stationary subsets of $\kappa$ whose
union is stationary, then there exists a choice $g(u)$ of one ordinal
from each $\Xi_u$ such that $g[U]$ is stationary.
\end{lemma}

\begin{lemma}\label{choicefrommorethankappa}
With the notation above,
let $C$ be a set of vertices satisfying
$|C| > \kappa$ and let
$\cf_v$ be an $X$-$v$ fan for every $v \in C$. Then
there exists an  $X$--$C$-warp $\cf$ such that $in[\cf] \supseteq
in[\cf_v]$ for some $v \in C$.
\end{lemma}

{\bf Remark:} As noted in \cite{aharonikonig}, Lemma
\ref{choicefrommorethankappa} follows easily from Theorem
\ref{main} (assuming it is proved). In fact, Theorem \ref{main}
has the following stronger corollary (written below in terms of the reverse web):

\begin{corollary}[of Theorem \ref{main}]\label{corofmain}
Assume that the web $\Gamma=(G,A,B)$ is unlinkable, and let
$\cf_a$ be an $a$-$B$-fan for every $a \in A$. Then there exists
an $A$--$B$-warp $\cf$ such that $ter[\cf] \supseteq ter[\cf_a]$
for some $a \in A$.
\end{corollary}

{\em Proof of Corollary \ref{corofmain}} Assuming the validity of
Theorem \ref{main}, there exist a family $\cp$ of disjoint paths
and an $A$--$B$-separating set $S$ such that $S$ consists of a
choice of one vertex from each $P \in \cp$. Since, by assumption,
$\Gamma$ is unlinkable, there exists $a \in A \setminus in[\cp]$.
Then $\cp[RF(S)] \hetz \cf_a$ is the desired warp $\cf$. $\enp$

{\em Proof of Theorem \ref{coreofkapimplieshind}} Let $POP$ be the
set of popular vertices of $\Lambda$, and let $UNP = V \setminus
POP$. Let $U_0=Y \cap UNP,~P_0=Y \cap POP$. Define inductively sets
$U_i, P_i~(i < \omega)$ as follows: $U_{i+1}=N^-(U_i)\cap
UNP,~P_{i+1}=N^-(U_i) \cap POP$. Finally, let $S = \bigcup_{i <
\omega}P_i$.

Since $X \subseteq POP$, we have $U_i \cap X =\emptyset$. Let $P$
be an $X$--$Y$-path having $k$ vertices. By the definition of the
sets $U_i$, if $P$ avoids $S$, then $V(P) \subseteq \bigcup_{i<k}
U_i$, thus $in(P) \not \in X$, a contradiction. This shows that
$S$ is separating.

\begin{assertion}\label{uiunpopular} $U_i$ is unpopular.\end{assertion}

\begin{proof}
 By induction on $i$. Suppose, first, that $U_0$ is popular.
Let $\cf$ be a $U_0$-joined family of  $X$-$U_0$-paths, such that
$f[in[\cf]]$ is stationary. For every $u \in U_0$ write $\cf_u=\{P
\in \cf,~ter(P)=u\}$. For every $\alpha \in f[in[\cf]]$ choose a
path $P \in \cf$ such that $f(in(P))=\alpha$, and define
$h(\alpha)=g(ter(P))$ (since $ter(P) \in U_0 \subseteq Y$, the
value $g(ter(P))$ is defined). By Definition
\ref{unbalanced}(\ref{regressive}), $h$ is regressive. Hence, by
Fodor's lemma (Theorem \ref{fodor}) there exist a stationary
subset $\Psi$ of $f[in[\cf]]$ and an ordinal $\beta$ such that
$h(\alpha)=\beta$ for every $\alpha \in \Psi$. This means that
there exists a vertex $u \in U_0$ such that $f[in[\cf_u]]$ is
stationary, contradicting the fact that $U_0 \subseteq UNP$.

Let now $k>0$, assume that the assertion is true for $i=k-1$, and
assume, for contradiction, that $U_k$ is popular. Let $\cf$ be a
$U_k$-joined family of $X$-$U_k$-paths, such that $f[in[\cf]]$ is
stationary. Again,  for every $u \in U_k$ write $\cf_u=\{P \in
\cf,~ter(P)=u\}$, and $\Xi_u= f[in[\cf_u]]$. Since $U_k \subseteq
UNP$, each set $\Xi_u$ is non-stationary. By Lemma
\ref{choicefromnonstationary}, there exists a choice of a path
$P(u)\in \cf_u$ for every $u\in U_k$, such that $f[in\{P(u)\mid ~u
\in U_k\}]$ is stationary. Since $U_k\subseteq N^-(U_{k-1})$, by
adding edges joining  $U_k$ to $U_{k-1}$, the family $\{P(u):~u
\in U_k\}$ can be extended to a $U_{k-1}$-joined family of paths.
But this contradicts the fact that  $U_{k-1}$ is unpopular.
\end{proof}

\begin{assertion}\label{pinotstronglypopular} $P_i$ is not strongly popular, for any $i <
\omega$.\end{assertion}

\begin{proof}

Assume that there exists an $X$-$P_i$-warp $\cp$ with $f[in[\cp]]$
stationary (this happens, in particular, if $f[P_i \cap X]$ is
stationary). The case $i=0$ follows from Lemma
\ref{unbalancedunlinkable}, since $P_0 \subseteq Y$. For $i>0$,
since $P_i \subseteq N^-(U_{i-1})$, the warp $\cp$ can be extended
to a
 $U_{i-1}$-joined family of paths $\cf$, with $in[\cf]=in[\cp]$.
This contradicts Assertion \ref{uiunpopular}.
\end{proof}

\begin{assertion} \label{pilessthankappa} $|P_i \setminus X| \le \kappa$
for every $i< \omega$.
\end{assertion}

\begin{proof}
Every point $p \in P_i \setminus X$ has a $p$-joined $X$-$p$ warp
$\cw_p$ such that $f(in[\cw_p])$ is stationary. If $|P_i
\setminus X| >\kappa$ then by Assertion \ref{choicefrommorethankappa} there exists
an $X$-$P_i$-warp $\cw$ such that $in[\cw] \supseteq in[\cw_p]$ for some $p \in P_i$,
implying that $in[\cw]$ is stationary, and hence that $P_i$ is strongly popular. This contradicts
Assertion \ref{pinotstronglypopular}.
\end{proof}

We are now ready to conclude the
proof of Theorem \ref{coreofkapimplieshind}.
Assertion \ref{pilessthankappa} yields condition (\ref{ssmall}) of
the theorem, and Assertion \ref{pinotstronglypopular} implies
condition (\ref{snotstronglypopular}). It remains to show
condition (\ref{spopular}), namely that a point $s \in S$ is not
only popular in $\Lambda$, but also in $\Lambda[RF^\circ(S)
\cup\{s\}]$. If $s \in X$ then there is nothing to prove.
Otherwise, there exists  an $s$-joined family $\cf$ of
$X$-$s$-paths such that $f[in[\cf]]$ is stationary. For each $i$
let $\cf_i$ be the set of those paths $P \in \cf$ on which there
exists a vertex $x\neq s$ in $P_i$ such that $xP$ meets $S$ only
at $x$. Since no $P_i$ is strongly popular, $f[in[\cf_i]]$ is
non-stationary for every $i<\omega$. Hence, by Lemma
\ref{stationaryisbig}, $f[in[\bigcup_{i<\omega} \cf_i]]$ is
non-stationary. Thus the set $\cf'$ of paths from $\cf$ meeting
$S$ only at $s$ satisfies the property that $f[in[\cf']]$ is
stationary. $\enp$

Clearly, the properties of the set $S$ in Theorem
\ref{coreofkapimplieshind} imply that $S$ is linkable in
$\overleftarrow{G}$ properly into $X$, which yields Theorem
\ref{halltypetheorem} for $\kappa$-unbalanced webs.\\

{\em Proof of Theorem \ref{kappahindranceimplieshindrance}}.~

By assumption, there exists in $\Gamma$ a $\kappa$-hindrance $\cl$
for some regular cardinal $\kappa$. We shall use for $\cl$ the
notation of Section 7. By Lemma \ref{phiastationary}, we may
assume that $\Phi_G=\Phi_G(\cl)$ is stationary.

Let $\cy = \cy(\cl)$. We wish to turn $\cy$ into a hindrance. In
fact, it almost {\em is} a hindrance: $ter[\cy]$ is
$A$--$B$-separating, and any $\alpha \in \Phi=\Phi(\cl)$ gives
rise to a path in $\ci\ce(\cy)$. The problem is that there are
paths in $\cy$ that ``hang in  air", namely they start at vertices
$y_\beta$. We wish to ``ground" such paths, using reverse
$\cy_G$-alternating paths from such vertices $y_\beta$ to some
$x_\alpha,~\alpha \in \Phi_G \setminus \Phi^\infty$ or to some
infinite path $H_\alpha$, $\alpha \in \Phi_G \cap \Phi^\infty$.
Applying such a path to $\cy$ ``connects $y_\beta$ to the ground".
We shall be able to do this only for ``popular" vertices
$y_\beta$, in a sense to be defined below. But using Theorem
\ref{coreofkapimplieshind}, we shall find that this suffices.

For every $\alpha \in \Phi_G \cap \Phi^\infty(\cl)$ let $x_\alpha$
be a new vertex added, which represents the infinite path
$H_\alpha$. Let $X^\infty$ be the set of vertices thus added.
 Let $X=X^{fin}(\cl) \cup X^\infty$ and $Y=Y(\cl)\cap V[\ce(\cy)]$ (see Notation
 \ref{yalphaandxalpha} for the definitions of $X^{fin}(\cl)$ and of
 $Y=Y(\cl)$.)
  To understand
 the choice of the definition of $Y$, note that only paths in $\ce(\cy)$
need to be ``connected to the ground", to obtain a wave.
  For each $\alpha \le \kappa$ write $T_\alpha
=T_\alpha(\cl)$. Write $T=T_\kappa$, namely $T=ter[\ce(\cy)]$.

Let $\tilde{D}=D[RF(T)]$. Let $F$ be the graph whose vertex set is
$RF(T) \cup X^\infty$, and whose edge set is $E(\tilde{D}) \cup
\{(x_\alpha,v) \mid u \in RF(T),~x_\alpha \in X^\infty,~(u,v) \in
E(D) ~\mbox{for some}~u \in V(H_\alpha)\}$. Let $\Theta$ be the web
$(F,X,Y)$, and let $\Lambda=\Lambda_\Theta(\cy)$, as defined in
Section \ref{conversion}. As recalled, $\Lambda$ is the web of
$\cy$-alternating paths in $\Theta$.

\begin{remark}
For the sake of clarity, we shall redefine the web $\Lambda$
explicitly. The definition of $\Lambda$ below is quite complex.
However, it is quite natural when viewed in the bipartite
conversion of $\Theta$, and it is advisable to keep in mind this
conversion. For example, it is helpful to remember that $X$
consists in the bipartite conversion of ``men", and hence can be
connected only to ``women". Since every edge $(u,v) \in E[\cy]$
corresponds to the edge $(m(u),w(v))$ in the bipartite conversion,
this means that $x \in X$ can be connected in $\Lambda$ only to
$v$.
\end{remark}

The vertex set of $\Lambda$ is $V_\Lambda=X \cup Y \cup (RF(T)
\setminus V[\cy]) \cup E[\cy]$.

The edge set of $\Lambda$ is constructed by the rule that an edge
$(u,v) \in E[\cy]$ sends an edge somewhere (namely, a vertex or an
edge) if $u$ sends there an edge in $D$, and it receives an edge
from somewhere if $v$ receives an edge from there (corresponding
to an edge ending at $w(v)$). We shall also have edges between two
consecutive edges $(u,v)$ and $(v,w)$ of $\cy$, the edge being
directed from the latter to the former (in the bipartite
conversion this means ``directed from the man to the woman". In
alternating paths terminology, this corresponds to the fact that
alternating paths go backwards on paths from $\cy$). Another rule
is that $X$-vertices only send edges, and $Y$ vertices only
receive edges. Finally, a vertex $x_\alpha \in X^\infty$ sends
edges in $\Lambda$ to all vertices (and, consequently, to edges)
to which some vertex on $H_\alpha$ sent an edge in $D$.

Formally, write:
$$E_{VV} = \{(u,v) \mid~ u \in (RF(T) \setminus
V[\cy]) \cup X,~ v \in (RF(T) \setminus V[\cy]) \cup Y , ~(u,v)
\in E(D)\}$$
$$E_{EV} =\{(e,w) \mid e=(u,v)\in E[\cy], ~w \in (RF(T) \setminus V[\cy]) \cup Y,~ (u,w) \in
E(D)\}$$
$$E_{VE} =  \{(w,e) \mid e=(u,v)\in E[\cy],~ w \in (RF(T) \setminus V[\cy])\cup X,~ (w,v) \in E(D)\}$$
$$E_{EE}= \{(e,f) \mid e=(u,v), ~f=(w,z) \in E[\cy], u=z {\mbox { or }} (v,w) \in E(D)\}$$
$$E_{\infty V} = \{(x_\alpha,v)\mid x_\alpha \in X^\infty, ~ v \in  (RF(T) \setminus
V[\cy]) \cup Y, (u,v) \in E(D) \mbox{ for some } u \in H_\alpha\}
$$
$$E_{\infty E} = \{(x_\alpha,e)\mid x_\alpha \in X^\infty, ~ e = (w,v) \in E[\cy], ~
(u,v) \in E(D) \mbox{ for some } u \in H_\alpha\} $$

Let $E_\Lambda=E_{VV} \cup E_{EV} \cup E_{VE} \cup E_{EE} \cup
E_{\infty V} \cup E_{\infty E}$. Let $D_\Lambda$ be the digraph
$(V_\Lambda,E_\Lambda)$, and define the web $\Lambda$ as
$(D_\Lambda, X,Y)$. For each $x=x_\alpha \in X$ define
$f(x)=\alpha$, and for each $y=y_\beta \in Y$ let $g(y)=\beta$.

\begin{assertion}\label{lambdaunbalanced}
$\Lambda$ is $\kappa$-unbalanced, as is witnessed by $f$ and $g$.
\end{assertion}

\begin{proof}
Condition (\ref{xstationary}) of Definition \ref{unbalanced} is
true since $f[X]=\Phi(\cl)$.  Condition (\ref{regressive}) is
tantamount to the fact that $g(ter(Q))<f(in(Q))$ for every
$X$--$Y$-alternating path $Q$ in $\Theta$. If $in(Q) \in X^{fin}$
then this follows from Lemmas \ref{stayingwithinroofedtalpha} and
\ref{appear_in_time}. If $in(Q) = x_\alpha \in X^\infty$, and the
first edge in $Q$ is $(x_\alpha,u)$, then in $D$ there exists an
edge $(v,u)$ for some $v \in H_\alpha$. Then $v \in RF(T_\gamma)$
for some $\gamma \le \alpha$, and thus, again by Lemma
\ref{appear_in_time}, $g(ter(Q))<\gamma$, yielding $g(ter(Q))<
\alpha$.
\end{proof}

Let  $S$ be an  $X$--$Y$-separating set as in Theorem
\ref{coreofkapimplieshind}. Write $S_V=S \cap V(D),~S_E=S \cap
E[\cy]$. Also write $\Theta-S$ for the web obtained from $\Theta$
by deleting $S_V$ from its vertex set, and deleting $S_E$ from its
edge set.

The fact that $S$ is $X$--$Y$-separating in $\Lambda$ implies that
there are no augmenting $\cy$-alternating paths in $\Theta-S$.
Namely:

\begin{assertion}\label{thetaminuss}
There are no $S$-avoiding $\cy$-alternating paths in $D$ from $X$
to $Y$.
\end{assertion}

%

Let $\cg=\cy -S_E$, namely the set of fragments of $\cy$ resulting
from the deletion of edges in $S_E$.


\begin{remark}\label{alternatingpathsmaystartatxalpha}
To understand the next assertion, it should be kept in mind that
there are $\cy$-alternating paths that start at some $x_\alpha$,
and as their first step  go backwards on an edge belonging to
$E[\cy]$. This type of alternating paths is again best understood
in terms of the bipartite conversion. In the bipartite conversion,
the first edge of the corresponding alternating path starts with
the edge $(m(x_\alpha),w(x_\alpha))$, which does not belong to
$E[\cy]$, as is the customary definition of alternating paths.
\end{remark}

\begin{assertion}\label{noblockingpointsonhalpha}
Let $H=H_\alpha$ be a path belonging to $\cg_G^f$ ($H$ is then a
finite path in $\ci\ce(\cy)$ not containing an edge from $S_E$),
such that $x=ter(H) \in X \setminus S$. Then there is no
$\cy$-alternating path avoiding $S$ from a vertex of $H$ to
$Y\setminus S$.
\end{assertion}

\begin{proof} Suppose that there exists such a path $Q$. Let $u$ be the  last
vertex on $Q$ lying on $H$. Then the path $\overleftarrow{H}uQ$ is
a $\cy$-alternating $X$--$Y$-path avoiding $S$ (see the remark
above), contradicting the fact that $S$ is separating in
$\Lambda$.
\end{proof}

\begin{notation} Denote by $\ch_\emptyset$ the set of paths $H=H_\alpha \in \cg_G$ such that
either:\\

(i)~$H$ is finite and $ter(H)\not \in S$, or:\\

(ii)~$H$ is infinite and no $\cy$-alternating,
$S$-avoiding path starts at a vertex of $H$ and ends at $Y
\setminus S$.
\end{notation}

Let $\cg'=\cg \setminus \ch_\emptyset$.

Let $RR$ be the set of vertices $v$ such that there exists an
$S$-avoiding $\cg$-alternating path starting at $v$ and
terminating at $Y \setminus S$.  Assertion
\ref{noblockingpointsonhalpha} implies:

\begin{assertion}\label{rrimpliesgprime}
If $P \in \cg$ and $V(P) \cap RR \neq
 \emptyset$ then $P \in \cg'$.
\end{assertion}

For each $P \in \cg'$ define $bl(P)$ to be:
\begin{itemize}
\item
 the first vertex on $P$ belonging to $RR$ if $V(P) \cap RR \neq
 \emptyset$, and:
\item
$ter(P)$, if $V(P) \cap RR \neq
 \emptyset$.
\end{itemize}

Let $BL= \{bl(P) \mid P \in \cg'\}$ and $BB=S_V \cup BL$.


\begin{assertion}\label{bbisabseparating}
$BB$ is $A$--$B$-separating.
\end{assertion}

(Remark: The idea of the proof is borrowed from the proof
of Theorem \ref{blockingsimply}.)

\begin{proof}
Since $T$ is $A$--$B$-separating, it suffices to show that $BB$ is
$A$-$T$-separating. Let $R$ be an $A$-$T$-path in $D$, and assume,
for contradiction, that $V(R) \cap BB = \emptyset$. Write
$t=ter(R)$. Since $t \in T=\ce(ter[\cy])$, and since by assumption
$t \not \in S_V$, it follows that $t=ter(P)$ for some path $P \in
\cg$. Since $P$ is finite, and since $ter(P)\in \ce(ter[\cy])$
(namely, $P$ cannot be some $H_\alpha$), $P \in \cg'$. Let $q
=bl(P)$. Since $t \not \in BB$, it follows that $t >_P q$. Let $Q$
be a $\cg$-alternating path from $q$ to $Y \setminus S$.

Assume, first, that $R$ does not meet any path of $\cg$ apart from
$P$. Then, in particular, $in(R) \not \in V[\cy]$, and hence
$in(R) \in X$. If $R$ does not meet $Q$, then the path
$Rt\overleftarrow{P}qQ$ is an $S$-avoiding $\cy$-alternating path
from $A$ to $Y$, contradicting Assertion \ref{thetaminuss}. If $R$
meets $Q$, and the last vertex on $R$ belonging to $Q$ is, say,
$v$ then $RvQ$ is an  $S$-avoiding $\cy$-alternating path from $A$
to $Y$, again providing a contradiction.

Thus we may assume that $R$ meets another path from $\cg$, besides
$P$. Let $P_1$ be the last path different from $P$ met by $R$, and
let $t_1$ be the last vertex on $R$ lying on $P_1$.  The path
$t_1Rt\overleftarrow{P}Z$ (or a "shortcut" of it, as in the
previous paragraph) witnesses the fact that $t_1 \in RR$, and
hence by Assertion \ref{rrimpliesgprime} $P_1 \in \cg'$. Let $q_1
=bl(P_1)$. Since by assumption $v_1 \not \in BB$, it follows that
$t_1 >_{P_1} q_1$. Let $Q_1$ be an $S$-avoiding $\cg$-alternating
path from $q_1$ to $Y \setminus S$. If $R$ does not meet any other
path, besides $P$ and $P_1$, belonging to $\cg$ then the path
$Rt_1\overleftarrow{P_1}q_1Q_1$ (or a shortcut of it) is an
$S$-avoiding $X$--$Y$ $\cg$-alternating path, contradicting
Assertion \ref{thetaminuss}. Thus we may assume that $R$ meets
still another path from $\cg$. Continuing this argument, we
eventually must reach a contradiction, since $R$ is finite.
\end{proof}

\begin{assertion} \label{mostlylong}
Let $p \in RF(T)$, and let $\cj$ be an $X$-$p$-in-fan of
$\cy$-alternating paths in $\Theta$, such that each path in $\cj$
meets some path in $\cy_H$ not containing $p$. Then $f[in[\cj]]$
is non-stationary.
\end{assertion}

\begin{proof}
Assume for contradiction that $f[in[\cj]]$ is stationary. For
each $P \in \cj$ choose $\beta=\beta(P)$ such that $P$ meets the
path $\cy(y_\beta)$. As before, by choosing a subfamily of $\cj$
if necessary, we may assume that $f$ is injective on $in[\cj]$.
Hence the function $h$ on $f[in[\cj]]$ defined by
$h(\alpha)=\beta(P)$ for that $P\in \cj$ for which
$f(in(P))=\alpha$, is well defined. By an argument as in the proof
of Assertion \ref{lambdaunbalanced}, $h(\alpha) < \alpha$, namely
$h$ is regressive. By Fodor's Lemma, this implies that
$f^{-1}(\beta)$ is of size $\kappa$ for some $\beta$. But this is
clearly impossible, since only finitely many paths from $\cj$ can
meet $\cy(y_\beta)$.
\end{proof}

\begin{assertion}\label{meeting_mainly_non_se_paths}
Let $p \in RF(T)$, and let $\cj$ be an $X$-$p$-fan of
$\cy$-alternating paths in $\Theta$, such that each path in $\cj$
meets a path in $\cg_H$ (namely, a fragment of $\cy -S_E$ hanging
in air) not containing $p$. Then $f[in[\cj]]$ is non-stationary.
\end{assertion}

\begin{proof}
Suppose that $f[in[\cj]]$ is stationary. Let $P \in \cj$. Choose a
path $W \in \cg_H$ that $P$ meets, and let $e$ be the last edge of
$P$ lying on $W$. Denote by $s$ the edge in $S_E$ such that
$head(s) = in(W)$. Going from $s$ along $W$ to $e$ and then
continuing along $P$ yields then a $\cy$ alternating path $Q(P)$
starting at $s$ and ending at $ter(P)$. Since the paths $Q(P)$ are
all disjoint, it follows that $S_E$ is strongly popular. But this
contradicts property (\ref{ssmall}) of $S_E$, as guaranteed by
Theorem \ref{coreofkapimplieshind}.
\end{proof}

\begin{assertion}\label{meetingaboveblockingpoint}
Let $Q$ be an $X$-starting $\cy$-alternating path avoiding
$S$. Suppose that $Q$ meets a path $P$ from $\cg$, and let $p$ be
the last point on $P$ belonging to $Q$ (thus $p=tail(e)$ for some
edge $e \in E(P) \cap E(\overleftarrow{Q})$). Then $p \le_P
bl(P)$.
\end{assertion}

\begin{proof} Assume that $bl(P) <_P p$. By the definition of $bl(P)$,
there exists  a $\cy$-alternating path $R$, starting at $bl(P)$,
ending in $Y$ and avoiding $S$. Then the $\cy$-alternating path
$Qp\overleftarrow{P}bl(P)R$ (or part of it, if $R$ meets $Q$,) is
an $S$-avoiding $X$--$Y$ $\cy$-alternating path, contradicting the
fact that $S$ is $X$--$Y$-separating in $\Lambda$.
\end{proof}

\begin{assertion}\label{separatingsubsetofbb}
There exists in $\Gamma$ a warp $\cv$ such that $in[\cv] \subseteq
A$ and $ter[\cv] =BB$.
\end{assertion}

\begin{proof}
  Let $\tilde{S}=S_V \setminus X \cup \{head(e) \mid e \in S_E\}$. Order
the points of $\tilde{S}$ as $(s_\theta:~\theta < \lambda)$, where
$\lambda \le \kappa$. By the properties of $S$, each $s_\theta$
has an $X$-$s_\theta$-fan $\cf_\theta$ in $\Theta-S$ of size
$\kappa$ of $\cy$-alternating paths, such that $f[in[\cf_\theta]]$
is stationary. By Assertion \ref{mostlylong} we may also assume
that no path in $\cf_\theta$ meets a path from $\cy_H$, namely:

(i)~All  paths in $\cf_\theta$ meet (apart from possibly at $s_\theta$) only paths from $\cy_G$.\\

By Assertion \ref{meeting_mainly_non_se_paths} we may further
assume that no path in $\cf_\theta$ meets a path in $\cg_H$,
namely:

(ii)~ All paths in $\cf_\theta$ meet (apart from possibly at
$s_\theta$) only paths from
$\cg_G$.\\

%

By induction on $\theta$, choose for each $s_\theta$ a
$\cy$-alternating path $Q_\theta \in \cf_\theta$, ending at
$s_\theta$ and satisfying:

(a)~ $Q_\theta$ does not meet any path from $\cy_G$ met by any
$Q_\delta,~\delta < \theta$.\\

(b)~$Q_\theta$ does not meet (apart from possibly at $s_\theta$)
any path from $\cy_H$.\\

(c)~$Q_\theta$ does not meet (apart from possibly at $s_\theta$)
any path from $\cg_H$.\\

Since the paths $Q_\theta$ avoid $S$, they are not only
$\cy$-alternating, but also $\cg$-alternating.  We now apply all
$Q_\theta$'s to $\cg$. Let $\cz$ be the resulting warp. We wish to
form a corresponding warp in $D$. The paths in $\cz$ which are not
contained in $D$ are paths $Z$ such that $in(Z)=x_\alpha \in
X^\infty$. Such a path was obtained by the application of an
alternating path $Q_\theta$ such that $in(Q_\theta)=x_\alpha$. Let
$(x,v)$ be the first edge of $Q_\theta$. By the definition of
$E(\Lambda)$, this means that $(p,v) \in E(D)$ for some $p \in
V(H_\alpha)$. Replace then $Z$ by $H_\alpha pZ$.

Denote by $\cu$ the resulting warp in $D$. Conditions (a), (b) and
(c) imply that there are no non-starting paths in $\cu$ and $in[\cu] \subseteq A$.
Assertion \ref{meetingaboveblockingpoint} together with condition (a)
imply that each path from $\cu$ intersects $BB$ at most once.
Assertion \ref{meetingaboveblockingpoint} also implies $BB \subseteq V[\cu]$.
Therefore, by pruning the warp $\cu$ we can obtain a warp $\cv$ with
$in[\cv] \subseteq A$ and $ter[\cv] = BB$ as required.

\end{proof}

Since $BB$ is separating, $\cv$ is a wave. By the equivalent
formulation of the main theorem, given in Conjecture
\ref{halltypeloose}, to complete the proof of the theorem it is
enough to show that $\cv$ is non-trivial,
which is clear. In fact, more than that is true: $\ce(\cv)$ is a
hindrance, in a strong sense. Since $S$ is not strongly popular in
$\Lambda$, the set $\{f(ter(Q_\theta) \mid \theta < \lambda\}$ is
non-stationary. Thus, the set $\Xi=\{\alpha \mid x_\alpha \not \in
ter[\cv]\}$ is stationary. Each $\alpha \in \Xi$  corresponds to
some (finite or infinite) path $H_\alpha$, unreached by any
$Q_\theta$, and thus belonging to $\ci\ce(\cv)$.

This completes the proof of Theorem
\ref{kappahindranceimplieshindrance}. To prove Theorem \ref{halltypetheorem}, and
thereby Theorem \ref{main}, it remains to prove the ``linkability theorem", Theorem \ref{C}.

\section{Proof of the Linkability Theorem}

Define the {\em height}  of a set $Y$ of vertices to be the
minimal cardinality of a subset $X$ of $V \setminus A$ for which
there exists a wave $\cw$ in $\Gamma\quo X$, such that $Y
\subseteq RF_\Gamma(ter[\cw])$.
The height of $\Gamma$ is defined as the
height of $V$.

\begin{definition}
A warp $\cw$ is a {\em half-way linkage} if it is an
$A$--$C$-linkage, with $ter[\cu]\subseteq C$, for some minimal
separating set $C$ for which $\Gamma \quo  C$ is unhindered.  Such
a set $C$ is called a {\em stop-over set} of $\cw$. Note that in
this definition $C$ is not uniquely determined by
$\cw$. The {\em altitude} of $\cw$ is the minimal height of such a
set $C$.
\end{definition}

We shall prove:

\begin{theorem}\label{halltype3}
Suppose that $\Gamma$ is unhindered.
Let $A' \subseteq A$ be a set of cardinality $\lambda$. Then
\begin{itemize}
\item ($\clubsuit$) If $(D,A \setminus A',B)$ is linkable then so is the web $(D,A,B)$.
\item ($\clubsuit\clubsuit$) There exists a half-way linkage
of altitude at most $\lambda$, linking $A'$ to $B$.
\end{itemize}
\end{theorem}

Theorem \ref{C} follows from ($\clubsuit$) upon taking $A'=A$.


To gradually impart the ideas of the proof of Theorem
\ref{halltype3}, let us first prove a few low cardinality cases.

\subsection*{Proof of ($\clubsuit$) for $\lambda = \aleph_0$}

This is the main result of \cite{countablelike}. The proof there
is very laborious, circumventing as it does Theorem
\ref{safelinking}. With the aid of the latter, ($\clubsuit$)
follows in the countable case by a classic ``Hilbert hotel''
argument. Let $\cf$ be a linkage in the web $(D,A \setminus
A',B)$. Let $A_0 =A'$. Choose a vertex $a \in A_0$, and using
Theorem \ref{safelinking} link it to $B$ by a path $P_1$, such
that $\Gamma-P_1$ is unhindered. Let $A_1=A_0 \cup in[\cf\langle
V(P_1)\rangle]$ (namely, $A_1$ is obtained by adding to $A_0$ all
initial points of paths from $\cf$ met by $P_1$). Choose a vertex
from $A_1$, different from $a$, and link it to $B$ by a path $P_2$
in $\Gamma-P_1$, such that $\Gamma-P_1-P_2$ is unhindered. Let
$A_2=A_1 \cup in[\cf\langle V(P_2)\rangle]$. Continuing this way,
and choosing wisely the order of the elements to be linked by
$P_i$, all elements of all $A_i$'s serve as $in(P_j)$ for some
$j$, and thus the set $A''=\bigcup A_i$ is linked to $B$ by the
warp $\cp=\{P_0,P_1,\ldots\}$, and all paths in $\cf\langle A
\setminus A'' \rangle$ are disjoint from all paths in $\cp$. Thus
$\cf\langle A \setminus A'' \rangle \cup \cp$ is a linkage of $A$.

\subsection*{Proof of ($\clubsuit\clubsuit$) for $\lambda = \aleph_0$
and $|V| = \aleph_1$}

Order the elements of $V$ as $(v_\theta:~\theta < \aleph_1)$.
Construct an $\aleph_1$-ladder $\cl$,  at each stage $\alpha$
choosing $y_\alpha$ to be the first $v_\theta$ not belonging to
$RF(T_\alpha)$ and choosing $\cw_\alpha$ to be a hindrance in
$\Gamma_\alpha$ if such exists. The construction of $\cl$
terminates after $\zeta\leq \aleph_1$ steps.

By the choice of the vertices $y_\alpha$, we have:

\begin{assertion}\label{everybodyroofedt}
$V = \bigcup_{\alpha \in \Sigma} RF(T_\alpha) = RF(\cl)$.
\end{assertion}

Write $\cy=\cy(\cl)$ and for $\alpha \le \zeta$  write
$\cy_\alpha=\cy_\alpha(\cl)$ (thus $\cy=\cy_\zeta$) and $T_\alpha
=T_\alpha(\cl)$.

Assume, first, that $\zeta$ is countable. By Assertion
\ref{everybodyroofedt}
$RF[T_\zeta] = V$ and hence
$T_\zeta = \ce(V) = B$.
Together with Lemma
\ref{most_y_reach_salpha} (applied with $\alpha = \zeta$)
this implies that $\cy \langle \sim  B
\rangle$ is countable. Thus, $A \setminus in[\cy \langle B
\rangle]$ is countable. Hence, by the case of ($\clubsuit$) proved
above, $\Gamma$ is linkable, which clearly implies
($\clubsuit\clubsuit$).

Thus we may assume that $\zeta=\aleph_1$. Since $\Gamma$ is
unhindered, by Theorem \ref{kappahindranceimplieshindrance} $\cl$
is not an $\aleph_1$-hindrance, and hence there exists a closed
unbounded set $\Sigma$ not intersecting $\Phi(\cl)$. By Lemma
\ref{hinderedimplieshindered}, $\Sigma \cap \Phi_h(\cl)
=\emptyset$, namely:


\begin{assertion}\label{obs2}
$\Gamma_\alpha$  is unhindered for every $\alpha \in \Sigma$.
\end{assertion}

 Assertion \ref{everybodyroofedt} implies:

\begin{assertion}\label{countableroofed}
For every countable set of vertices $X$  there exists $\gamma(X) \in \Sigma$
 such that $X \subseteq RF(T_{\gamma(X)})$.
\end{assertion}

\begin{assertion}\label{mostysmeetsbeta}
$\cy \langle T_\alpha \rangle \setminus \cy \langle T_\beta
\rangle$ is countable for every $\alpha, \beta \in \Sigma$.
\end{assertion}

\begin{proof}

If $\beta < \alpha$ then $\cy \langle T_\alpha \rangle \setminus
\cy \langle T_\beta \rangle$ consists of those paths in $\cy$ that
start at some $y_\gamma$ for some $\beta \leq \gamma < \alpha$,
and thus it is countable. For $\alpha < \beta$, we have $\cy
\langle T_\alpha \rangle \setminus \cy \langle T_\beta \rangle
\subseteq \ci\ce(\cy_\beta)$, and hence the assertion follows from
Lemma \ref{fewhindered}.
\end{proof}

In particular, $\cy_G \setminus  \cy\langle T_\alpha
\rangle=\cy\langle T_0 \rangle \setminus  \cy\langle T_\alpha
\rangle$ is countable for every $\alpha \in \Sigma$~$\quad$
(remember that ``$\cy_G$" stands for ``$\cy \langle A \rangle$").

Write $A_0=A'$. Choose $a_0 \in A_0$, and using Theorem
\ref{safelinking} link it to $B$ by a path $P_0$, such that
$\Gamma-P_0$ is unhindered. Let $\gamma_0=\gamma(V(P_0))$. (See
Assertion \ref{countableroofed} for the definition of $\gamma$.)
Let $A_1=A_0 \cup in[\cy_G \langle V(P_0) \rangle] \cup in[\cy_G
\setminus \cy \langle T_{\gamma_0} \rangle]$. By Assertion
\ref{mostysmeetsbeta} $A_1$ is countable.

Choose $a_1 \in A_1 \setminus\{a_0\}$, and find an $a_1$-$B$ path
$P_1$ such that $\Gamma -P_0 -P_1$ is unhindered. Let
$\gamma_1=\max(\gamma(V(P_0)),\gamma(V(P_1)))$, and $A_2=A_1 \cup
in[\cy_G\langle V(P_1) \rangle] \cup in[\cy_G \setminus \cy
\langle T_{\gamma_1} \rangle]$.

Continue this way $\omega$ steps. Let $X = \cup_{i < \omega}
V(P_i)$,  and  $\gamma=\sup_{i < \omega} \gamma_i$. Since $\Sigma$
is closed, $\gamma \in \Sigma$. By Lemma \ref{noinfinitepathsiny}
every path $P \in \cy_G
 \setminus \cy \langle T_\gamma \rangle$  must belong to $\cy_G
 \setminus \cy \langle T_{\gamma_i} \rangle$ for some $i < \omega$
 and then, by the definition of the sets $A_i$, we have $in(P) \in A_{i+1}$.
 Note that each path $P_i$ ends at some vertex in $B \cap
 RF(T_\gamma)$ and since a vertex in $B$ can only be roofed by itself,
 this vertex must be in $T_\gamma$.

Choosing  the vertices $a_i$ in an appropriate order, we can see
to it that $\{a_i : ~ i < \omega \}=A' \cup in[\cy_G
 \setminus \cy \langle T_\gamma \rangle] \cup in[\cy\langle
X \rangle]$. Write $\cp=\{P_i:~i < \omega\}$, and let $\cv= \cp
\cup \cy \langle \sim X \rangle [RF(T_\gamma)] \langle A \rangle$.
Then $\cv$ is an $A$-$T_\gamma$-linkage linking $A'$ to $B$. By
Assertion \ref{obs2}, $\Gamma\quo T_\gamma$ is unhindered and
therefore, taking $C=T_\gamma$ in the definition of ``half-way
linkage" shows, by Lemma \ref{talpha_minimal}, that  $\cv$ is a
half-way linkage. The warp $\cy_\gamma \quo Y_\gamma(\cl)$ is a
wave in $\Gamma \quo Y_\gamma(\cl)$, whose terminal points set
contains $T_\gamma$, showing that $\cv$ has countable altitude.

This concludes the proof of ($\clubsuit\clubsuit$) for $\lambda =
\aleph_0$ and $|V| = \aleph_1$.

\subsection*{Proof of ($\clubsuit$) for $\lambda = |V| = \aleph_1$}

This was proved in \cite{aharoni97},
assuming Theorem \ref{safelinking}.
The arguments given here are
more involved, but fit better our general proof scheme.

We may clearly assume that $A' = A$. Again, construct an
$\aleph_1$-ladder $\cl$, for which Assertion
\ref{everybodyroofedt} holds. Let $\Sigma$ be defined as above
(once again using Theorem \ref{kappahindranceimplieshindrance}).

In the construction of $\cl$, we take each $\cw_\alpha$ to be a
hindrance in $\Gamma_\alpha$, if such exists. By Corollary
\ref{maximalhindrance}, we may also assume that $\cw_\alpha$ is a
maximal wave in $\Gamma_\alpha$ ($\fextended$-maximal and thus
also $\le$-maximal). The maximality of $\cw_\alpha$ implies:

\begin{assertion} For all $\alpha < \aleph_1$, every
wave in $\Gamma_\alpha$ is roofed by $T_{\alpha+1}$.
\end{assertion}

which implies:

\begin{corollary}\label{maxwave2} Whenever $\alpha < \beta < \aleph_1$, every
wave in $\Gamma_\alpha$ is roofed by $T_\beta$.
\end{corollary}

\begin{assertion}
\label{sbetaplus1}
 If $\alpha< \zeta$ and $X \subseteq RF_{\Gamma_\alpha}(T_\zeta)$ then every wave
in $\Gamma_\alpha\quo X$ is roofed by $T_{\zeta+1}$.
\end{assertion}

\begin{proof}
Let $\cv$ be a wave in $\Gamma_\alpha \quo X$.
Then $\cv \quo T_\zeta$ is a wave in
$(\Gamma_\alpha \quo X) \quo T_\zeta = \Gamma_\zeta$.
By Corollary
\ref{maxwave2}, the wave $\cv \quo T_\zeta$ is roofed by
$T_{\zeta+1}$, which implies that $\cv$ is roofed by
$T_{\zeta+1}$.
\end{proof}


The core of the proof is in the following:

\begin{assertion}
\label{alphabetalinkage} Let  $\alpha$ be an ordinal in
$\Sigma$, and let $U$ be a countable subset of $T_\alpha$. Then
there exist $\beta > \alpha$ in $\Sigma$ and a $T_\alpha$-$T_\beta$ linkage
$\ct$ linking $U$ to $B$, such that all but at most countably many
paths of $\ct$ are contained in paths of $\cy$.
\end{assertion}

\begin{proof} By the special case of ($\clubsuit\clubsuit$) proved above, there
exists in $\Gamma_\alpha$ a half-way linkage $\cu$ of altitude
$\aleph_0$, linking $U$ to $B$. Let $C$ be a stop-over set of
$\cu$, of height $\aleph_0$. We claim that there exists $\beta >
\alpha$ in $\Sigma$ such that $C \subseteq RF(T_\beta)$. The fact
that $\cu$ has altitude $\aleph_0$ means that $C$ is roofed by a
wave in $(\Gamma \quo T_\alpha) \quo  X$ for some countable set
$X$. Take $\beta \in \Sigma$ such that  $\beta > \max(\alpha,
\gamma(X))$ (where $\gamma(X)$ is defined as in Assertion
\ref{countableroofed}, which is valid also in the present case).
By Assertion \ref{sbetaplus1} we know that every wave in $(\Gamma
\quo T_\alpha) \quo  X$ is roofed by $T_\beta$ and thus also $C$
is roofed by $T_\beta$.

\newcommand{\yab}{{\cy \langle T_\alpha \rangle \langle T_\beta \rangle}}

By Lemma \ref{sandwich}, the set $C$ is $T_\alpha$--$T_\beta$-separating,
and thus

\begin{equation}
\label{yabcontainedinyc} \yab \subseteq \cy\langle C \rangle.
\end{equation}

Note that Assertion \ref{mostysmeetsbeta} holds here
(with the same proof as in the previous case), and together with Equation
(\ref{yabcontainedinyc}), it yields:

\begin{equation}
\label{cytalphaminuscyc}
 |\cy \langle T_\alpha \rangle
\setminus \cy \langle C \rangle| \leq \aleph_0.
\end{equation}

Let $J$ be the graph on $V(D)$ whose edge set is $E[\cu] \cup
E[\cy]$. By (\ref{cytalphaminuscyc}), at most countably many
connected components of $J$ contain vertices of $U$ or paths from
$\cy \langle T_\alpha \rangle \setminus \cy \langle C \rangle$. In
all other connected component of $J$ we can replace the paths of
$\cu$ by the segments of the paths of $\cy$ between $T_\alpha$ and
$C$ while maintaining the properties of $\cu$ as being a
$T_\alpha$-$C$ linkage linking $U$ to $B$. Therefore we may assume
that all but countably many paths in $\cu$ are contained in paths
of $\cy$.

Similarly to (\ref{cytalphaminuscyc}) we have:

\begin{equation}
\label{cytbetaminuscyc}
 |\cy \langle T_\beta \rangle
\setminus \cy \langle C \rangle| \leq \aleph_0.
\end{equation}

This implies that there exists a warp $\cf$, whose paths are parts
of paths of $\cy$, linking all but countably many vertices of
$ter[\cu]$ to $T_\beta$.

We may clearly assume (and hence will assume) that each path  $P
\in \cu$ meets $C$ only at $ter(P)$ and therefore $V[\cu]
\setminus ter[\cu] \subseteq RF^\circ(C)$. However, a path $F \in
\cf$
may intersect $C$ many times and may pass through $RF^\circ(C)$.
We wish to use $F$ in the construction of the desired linkage
$\ct$, which explains the necessity of the term $V[\cf]$ in the
following definition: define  $\Delta$ as the web $(D[(RF(T_\beta)
\setminus RF^\circ(C)) \cup V[\cf]],ter[\cu],T_\beta)$. Clearly,
$\Delta \quo (C \setminus ter[\cu]) = (\Gamma \quo C)
[RF(T_\beta)]$.  By Observation \ref{waveinduced} since $\Gamma
\quo C$ is unhindered so is $\Delta \quo (C \setminus ter[\cu])$,
and hence by Corollary \ref{hindranceinquotient} $\Delta$ is
unhindered.

We now apply the case $\lambda = \aleph_0$ of ($\clubsuit$) to
$\Delta$ and $A' = ter[\cu] \setminus in[\cf]$. This gives a
linkage $\cq$ of $ter[\cu]$ to $T_\beta$. By arguments similar to
those given above, we may assume that all but countably many paths
of $\cq$ are contained in paths of $\cy$. The concatenation $\cu *
\cq$ is then the linkage $\ct$ desired in the assertion.
\end{proof}

We now use Assertion \ref{alphabetalinkage} to prove
($\clubsuit$). The general idea of the proof is to link ``slices"
of the web, lying between $T_\alpha$'s, for ordinals $\alpha \in
\Sigma$. Assertion \ref{alphabetalinkage} is used to avoid the
generation of infinite paths in this process. By Lemma
\ref{phiinftycontainedinphi}, paths belonging to $\cy$ do not
become infinite along this procedure. Thus we have to be careful
only about paths not contained in paths from $\cy$. Using the
assertion, at each stage  we can take care of such paths, by
linking their terminal points to $B$.

Formally, this is done as follows. Write $A$ as $\{a_\alpha: ~
\alpha < \omega_1\}$, and let $U_0 = \{a_0\}$. Use the assertion
to find $\sigma_1 < \omega_1$ in $\Sigma$ and an
$A$-$T_{\sigma_1}$ linkage $\ct_0$, linking $a_0$ to $B$,  such
that at most countably many paths of $\ct_0$ are not contained in
a path of $\cy$. Let $U_1$ be the set of end vertices of such
paths, together with the end vertex of the path in $\ct_0$
starting at $a_1$.

We  use the assertion in this way, to define inductively ordinals
$\sigma_\alpha \in \Sigma$ and
$T_{\sigma_\alpha}$-$T_{\sigma_{\alpha + 1}}$ linkages
$\ct_\alpha$ linking $U_\alpha$ to $B$. Having defined these up to
and including $\alpha$, we write $\ct_{\le\alpha}=*(\ct_\theta:
\theta \le \alpha)$ and $\ct_{<\alpha}=*(\ct_\theta: \theta <
\alpha)$. Let $U_{\alpha + 1}$ consist of the end vertices of all
paths in $\ct_{\le \alpha}$ not contained in a path of $\cy$, together
with the end vertex of the path in $\ct_{\le\alpha}$ starting at
$a_{\alpha+1}$.

\begin{assertion} $\ct_{<\alpha}$ is an $A-S_{\sigma_\alpha}$ linkage.
\end{assertion}

\begin{proof}
For successor $\alpha$, this follows by induction from the
definitions. For limit $\alpha$, this follows from Lemma
\ref{noinfinitepathsiny}, and the fact that, by our construction,
all paths in $\ct_{<\alpha}$ not contained in a path from $\cy$
terminate in $B$.
\end{proof}

For limit $\alpha$ we take $U_\alpha = ter[\ct_{<\alpha} \langle
\{a_\alpha\} \rangle]$ and $\sigma_\alpha = \sup_{\theta < \alpha}
\sigma_\theta$.

Since $a_\alpha$ is linked to $B$ by $\ct_\alpha$, the
concatenation $\ct$ of $(\ct_\alpha : ~ \alpha < \omega_1)$ is the
desired $A$--$B$ linkage.

This concludes the proof of ($\clubsuit$) for $\lambda = |V| =
\aleph_1$.
\\
We now go on to the proof of $(\clubsuit)$ and $(\clubsuit \clubsuit)$ in the general case.\\


{\bf Proof of ($\clubsuit$)  (assuming ($\clubsuit$) and ($\clubsuit \clubsuit$) for
cardinals smaller than $\lambda$)} \\{\bf Case I: ~~$\lambda$ is
regular.}

Let $\cf$ be a linkage in the web $(D,A \setminus A',B)$.
Similarly to the $\lambda=\aleph_1$ case, we construct a
$\lambda$-ladder $\cl$ and a choose a closed unbounded set $\Sigma
\subseteq \lambda$ disjoint from $\Phi(\cl)$. At each stage
$\alpha$ we take $\cw_\alpha$ to be a maximal hindrance in
$\Gamma_\alpha$, if $\Gamma_\alpha$ is hindered. Then Corollary
\ref{maxwave2} and Assertion \ref{sbetaplus1} are valid also here.

Let $\cy=\cy(\cl)$. We then have the analogue of Assertion
\ref{mostysmeetsbeta}:

\begin{assertion}\label{mostysmeetsbetagen}
$|\cy \langle T_\alpha \rangle \setminus \cy \langle T_\beta
\rangle| < \lambda$ for every $\alpha, \beta \in \Sigma$.
\end{assertion}

(For the notation used, see Convention \ref{notationonladders}.)




The difficulty we may face is that possibly $|V| > \lambda$. This
means that Assertion \ref{everybodyroofedt} may fail, namely we
cannot guarantee that every vertex is roofed by some $T_\alpha$.
We can only hope to achieve this for $\lambda$ many vertices.
Fortunately, this suffices. Along with the construction of the
rungs $R_\alpha$ of $\cl$, we shall define sets $Z_\alpha$ of
cardinality  at most $\lambda$, each of whose elements we shall
wish to roof by $T_\beta$ for some $\beta > \alpha$.

Having defined $Z_\theta$, we enumerate its elements as
$(z_\theta^\beta : ~ \beta <  |Z_\theta| \le \lambda)$.

To define $Z_\alpha$, we do the following. Assume that the rungs
$R_\beta$ of $\cl$ as well as the sets $Z_\beta$ have been defined
for $\beta < \alpha$. Write $Z_{<\alpha} = \bigcup_{\theta <
\alpha} Z_\theta$ and $Z_{<\alpha}^{<\alpha} = \{z_\beta^\gamma :
~ \beta < \alpha, ~ \gamma < \alpha \}$.

Let $(\gamma, \delta)$ be a pair of ordinals such that $\alpha=
\max(\gamma,\delta)$. Consider two cases:

\begin{itemize}
\item
$\Gamma_\delta$ is unhindered. Apply then
($\clubsuit\clubsuit$), which by the inductive hypothesis is
true when $|A'| < \lambda$, to the web $\Gamma_\delta$ with $A'=
T_\delta \cap Z_{<\gamma}^{<\gamma}$. This yields the existence of
a half-way linkage $\ca=\ca_{\delta,\gamma}$ in $\Gamma_\delta$,
linking $T_\delta \cap Z_{<\gamma}^{<\gamma}$ to $B$. Furthermore,
$\ca$ is  of height less than $\lambda$, namely it is roofed by
some wave in $\Gamma_\delta \quo  X_{\delta,\gamma}$ for some set
$X_{\delta,\gamma}$ of cardinality less than $\lambda$.
\item
$\Gamma_\delta$ is hindered. In this case let
$X_{\delta,\gamma}=\emptyset$.
\end{itemize}

Let $(\beta, \gamma, \delta)$ be a triple of ordinals such that
$\delta < \beta$ and $\alpha = \max(\beta, \gamma)$. Consider the
following two cases:

\begin{itemize}
\item
There exists a $T_\delta $-$T_\beta$-linkage linking $T_\delta
\cap Z_{<\gamma}^{<\gamma}$ to $B$, in which all paths are
contained in paths of $\cy_\beta$ except for a set of size smaller
than $\lambda$. In such a case choose such a linkage and denote it
by $\cu_{\beta, \gamma, \delta}$. Write
$\cu_{\beta,\gamma,\delta}^m$ for the set of paths in
$\cu_{\beta,\gamma,\delta}$ not contained in a path of $\cy$ (the
``m" standing for ``maverick").

\item
There does not exist such a linkage. Write then
$\cu_{\beta,\gamma,\delta}^m=\emptyset$.

\end{itemize}

Let $Z_0 = A'$ and for $\alpha > 0$ let

$$Z_\alpha= Z_{<\alpha} \cup V(H_\alpha) \cup \{y_\alpha\}
\cup V[\cf \langle Z_{< \alpha} \rangle] \cup V[\cy \langle Z_{<
\alpha} \rangle] \cup
 \bigcup_{\begin{array}{cl}
        \delta \leq \alpha  \\
        \gamma \leq \alpha
           \end{array}} X_{\delta,\gamma}
\cup \bigcup_{\begin{array}{cl}
        \delta < \beta \leq \alpha  \\
        \gamma \leq \alpha
           \end{array}} V[\cu_{\beta,\gamma,\delta}^m]$$

Let $Z=\bigcup_{\alpha< \lambda}Z_\alpha$. By the regularity of
$\lambda$ we have:

\begin{assertion}
\label{ztagsubsetzalphaalpha} Every subset $U$ of $Z$ of
cardinality less than $\lambda$ is contained in
$Z_{<\alpha}^{<\alpha}$ for some $\alpha < \lambda$.
\end{assertion}

Choosing carefully the vertices $y_\alpha$ in the ladder $\cl$, we
can see to it that the following weaker version of Assertion
\ref{everybodyroofedt} holds:

\begin{assertion}
\label{everyzroofedt} $Z \subseteq RF(\cl)$.
\end{assertion}

We now have the analogue of Assertion \ref{alphabetalinkage}, with
practically the same proof:

\begin{assertion}
\label{alphabetalinkage2} For every $\alpha \in \Sigma$ and every
 subset $U$ of $T_\alpha \cap Z$ having cardinality less than
$\lambda$, the following is true: there exist $\beta > \alpha$ and
a $T_\alpha$-$T_\beta$ linkage $\ct$ linking $U$ to $B$, such that
all but fewer than $\lambda$ paths of $\ct$ are contained in paths
of $\cy$, and $V(P) \subseteq Z$ for each  path $P \in \ct$ not
contained in a path of $\cy$.
\end{assertion}

>From here the proof continues in a way  similar to that of the
$\aleph_1$ case.  We define inductively ordinals  $(\sigma_\alpha
: ~ \alpha < \lambda)$, warps $\ct_\alpha$ and subsets $U_\alpha$
 of $T_{\sigma_\alpha}$, as follows.
 Enumerate $Z \cap A$ as $(z_\alpha: ~ \alpha < \lambda)$ and let
$U_0 = \{z_0\}$, $\sigma_0 = 0$. Assume now that $\sigma_\alpha$
and $U_\alpha$ have been defined. Use Assertion
\ref{alphabetalinkage2} to find an ordinal
$\beta=\sigma_{\alpha+1}> \sigma_\alpha$ in  $\Sigma$, and a
$T_{\sigma_\alpha} $-$T_{\sigma_{\alpha + 1}}$-linkage
$\ct_\alpha$,  linking $U_\alpha$ to $B$ and satisfying the
conditions stated in the assertion.

Let $U_{\alpha + 1}$ consist of the terminal vertex of the path in
$*(\ct_\theta : \theta \leq \alpha)$ starting at $z_{\alpha+1}$,
together with the  terminal points of all those paths in
$\ct_\alpha$ that are not contained in a path of $\cy$.

For limit $\alpha$ let $U_\alpha = ter[*(\ct_\theta :~ \theta <
\alpha) \langle \{z_\alpha\} \rangle]$ and $\sigma_\alpha =
\sup_{\theta < \alpha} \sigma_\theta$.

Having defined all these for all $\alpha < \lambda$, we define
$\ct = *(\ct_\alpha : ~ \alpha < \lambda)$. For each $\beta$, the
vertex $z_\beta \in Z \cap A$ is linked to $B$ by $*(\ct_\alpha :~
\alpha \le \beta)$, and thus it is linked to $B$ by $\ct$.
Every $a \in A \setminus Z$ is the initial point of some path $P \in \cf$.
By the way we chose $Z$, we have
$V[\ct \langle Z \rangle] \cap V[\cf \langle \sim Z \rangle] = \emptyset$
and therefore the set
$\ct \langle Z \rangle \cup \cf \langle \sim Z \rangle$ is a warp. This is
the desired $A$--$B$-linkage,
completing the proof of ($\clubsuit$).\\
\
\\

{\bf Proof of ($\clubsuit$), Case II: $\lambda$ is singular.}

\begin{definition}
Given a set $\cp$ of paths, two vertices $u,v$  are said to be
{\em competitors in} $\cp$ if there exist $P,Q \in \cp$ such that
$in(P) = u$, $in(Q)=v$ and $V(P) \cap V(Q) \neq \emptyset$.
\end{definition}

Note that if $\cp$ is the union of $\mu$ warps, then each vertex
has at most $\mu$ competitors.

Let $\cf$ be a linkage in $(D,A \setminus A',B)$. Let
$\mu=cf(\lambda)$ and let $(\kappa_\alpha : ~ \alpha < \mu)$ be a
sequence converging to $\lambda$. We may assume that $\kappa_0 >
\mu$.

Call a matrix of sets {\em increasing} if each row and each column
of the matrix is ascending with respect to the relation of
containment.

\begin{assertion}\label{matrices}

There exist two $\mu \times \omega$ matrices: an increasing matrix
of sets $(A_\alpha^k : ~ \alpha < \mu, ~ k < \omega)$ and a matrix
of half-way linkages $(\cw_\alpha^k : ~ \alpha < \mu, ~ k <
\omega)$, jointly satisfying  the following properties:

\begin{enumerate}
\renewcommand{\theenumi}{\roman{enumi}}
\item\label{rightcardinalities} $| A_\alpha^k | = \kappa_\alpha$.
\item\label{allthere0} $ \bigcup_{\alpha < \mu} A_\alpha^0 = A'$.
\item\label{takecarek}
$\cw_\alpha^k$ links $A_\alpha^k$ to $B$.
\item\label{competitorsk} If $a \in A_\alpha^k$
then all competitors of $a$ in $\cf \cup \bigcup_{\beta < \mu}
\cw_\beta^k$ are in $A_\alpha^{k+1}$.
\item\label{increasinghalfways}
For every $\alpha < \mu$ the sequence $(\cw_\alpha^k: ~ k <
\omega)$ is  $\fextended$-increasing (as a sequence of warps).
\end{enumerate}

\end{assertion}
\begin{proof}
We first choose $(A_\alpha^0: \alpha < \mu)$ that satisfy
conditions (\ref{rightcardinalities}) and (\ref{allthere0}). We
use ($\clubsuit\clubsuit$) of the induction hypothesis to obtain
half-way linkages $(\cw_\alpha^0: \alpha < \mu)$ that satisfy
(\ref{takecarek}). Denote the stop-over set of $\cw_\alpha^0$ by $C_\alpha^0$.
We now define $A_\alpha^1$ to be the set of all
competitors of members of $A_\alpha^0$ in $\cf \cup \bigcup_{\beta
< \mu} \cw_\beta^0$. We then use ($\clubsuit\clubsuit$) for the
webs $\Gamma \quo  C_\alpha^0$  to get $(\cw_\alpha^1: \alpha <
\mu)$ that satisfy conditions (\ref{takecarek}) and
(\ref{increasinghalfways}). We continue this way, where at each
step we define $A_\alpha^{k+1}$ to be the set of all competitors
of members of $A_\alpha^k$ in $\cf \cup \bigcup_{\beta < \mu}
\cw_\beta^k$ and we use ($\clubsuit\clubsuit$) to get
$(\cw_\alpha^{k+1}: \alpha < \mu)$ that satisfy conditions
(\ref{takecarek}) and (\ref{increasinghalfways}). Condition
(\ref{rightcardinalities}) is satisfied since no vertex has more
than $\mu$ competitors at any stage.

\end{proof}

\begin{assertion}

There exist an ascending sequence of subsets $(A_\alpha
: ~ \alpha < \mu)$ of $A$ and a sequence of warps $(\cw_\alpha : ~
\alpha < \mu)$, satisfying together the following properties:

\begin{enumerate}
\item\label{takecare} $\cw_\alpha$ links $A_\alpha$ to $B$.
\item\label{allthere} $ \bigcup_{\alpha < \mu} A_\alpha \supseteq A'$.
\item\label{competitors} If $a \in A_\alpha$ then all competitors of $a$
in $\cf \cup \bigcup_{\beta < \mu} \cw_\beta$ are also in
$A_\alpha$.
\end{enumerate}

\end{assertion}

\begin{proof}
Let $(A_\alpha^k)$ and $(\cw_\alpha^k)$ be as in Assertion
\ref{matrices}. Take $A_\alpha = \bigcup_{k < \omega} A_\alpha^k$
and $\cw_\alpha = \lim_{k < \omega} \cw_\alpha^k$. Conditions
(\ref{takecarek}) and (\ref{increasinghalfways}) imply
(\ref{takecare}), condition (\ref{allthere0}) implies
(\ref{allthere}) and condition (\ref{competitorsk}) implies
(\ref{competitors}) because every two competitors in $\cf \cup
\bigcup_{\beta < \mu} \cw_\beta$ are competitors in $\cf \cup
\bigcup_{\beta < \mu} \cw_\beta^k$ for some $k$.

\end{proof}

 We can now conclude the proof of ($\clubsuit$).
For every $a \in \bigcup_{\alpha < \mu} A_\alpha$ use the path to
$B$ in $\cw_\alpha$ to link $a$ to $B$, where $\alpha$ is minimal
with respect to the property that $a \in A_\alpha$. Such a path
exists by condition (\ref{takecare}). For every $a \in A \setminus
\bigcup_{\alpha < \mu} A_\alpha$, we know by condition
(\ref{allthere}) that $a \in A \setminus A' = in[\cf]$, and hence
we can link $a$ to $B$ by the path in $\cf$ starting at $a$.
Condition (\ref{competitors}) guarantees that these paths are
disjoint.\\
\
\\

{\bf Proof of ($\clubsuit\clubsuit$) for general $\lambda$
(assuming ($\clubsuit$) for cardinals $\le \lambda$)}

Recall that in the case  $\lambda = \aleph_0$ and $|V| = \aleph_1$
we used an $\aleph_1$-ladder. Analogously, for general $\lambda$
we construct a $\lambda^+$-ladder, $\cl$.

As before, since by Theorem \ref{kappahindranceimplieshindrance} $\cl$ is not a $\lambda^+$-hindrance,
there exists a closed unbounded set $\Sigma$,
disjoint from $\Phi(\cl)$. Replacing $\lambda$ by $\lambda^+$, we
then have the analogues of Corollary \ref{maxwave2} and Assertions
\ref{sbetaplus1},
 \ref{mostysmeetsbetagen} and \ref{everyzroofedt}.

The basic idea of the proof is relatively simple. We wish to use
($\clubsuit$)  for $\lambda$, which is true by the inductive
assumption, in order to ``climb" $\cl$. This is done as follows:
Order $A'$ as $(a_i \mid i < \lambda)$. Use Theorem
\ref{safelinking} to link $a_0$ to $B$ by a path $P$ so that
$\Gamma - P$ is  unhindered. Choose $\alpha_1 \in \Sigma$ such
that $V(P) \subseteq RF(T_{\alpha_1})$. Then use Lemma
\ref{most_y_reach_salpha} and the fact that ($\clubsuit$) holds
for $\lambda$, to complete $P$ to a linkage $\ck_1$ of $A$ into
$T_{\alpha_1}$. Then repeat the procedure with the web
$\Gamma_{\alpha_1}$ replacing $\Gamma$, and the element in
$T_{\alpha_1}$ to which $a_1$ is linked by $\ck_1$ replacing
$a_0$. After $\lambda$ such steps, $A'$ is linked to $B$, and $A$
is linked to some $T_\gamma$.

As usual, the  problem is the possible generation of infinite
paths. To avoid this, we have to anticipate which vertices may
participate in infinite paths, and link them to $B$ by the
procedure described above.  The trouble is that we can take care
in this way only of $\lambda^+$ such vertices. It is possible for
a vertex from $A'$ to have degree larger than $\lambda^+$, and
then it may be necessary to add more than $\lambda^+$ vertices to
the set $Z$ of vertices ``in jeopardy". The concept used to solve
this problem is that of {\em popularity} of vertices, having in
this case a slightly different meaning from the ``popularity" of
the previous section.
 ``Popularity"
of a vertex $z$ means that there exist many $z$-joined
$\cy$-s.a.p's emanating from $z$, and going to infinity or to $B$. (In
this sense the concept was used in
\cite{countablelike} and \cite{aharonidiestel}. A similar notion, solving
a similar problem, was used in \cite{inftutte}).
A popular vertex does not need to be taken care of immediately, since
it can be linked at a later stage, using its popularity. Thus we
have to perform the closure operation only with respect to
non-popular vertices, and this indeed will necessitate adding only
$\lambda^+$ vertices to $Z$.

A first type of vertices which should be considered ``popular" are
those that do not belong  to $RF^\circ(T_\alpha)$ for any $\alpha<
\lambda^+$. Note that for each vertex $v$, the set $\{\theta : ~ v
\in T_\theta\}$ is an interval, namely it is either empty or of
the form $\{\theta:~ \alpha \leq \theta < \beta\}$ for some
$\alpha < \beta \leq \lambda^+$. Let $T_{\lambda^+}$ be the set of
vertices for which this  set is unbounded in $\lambda^+$. By Lemma
\ref{slambdaplus} we have:

\begin{assertion}\label{unbounded}
$T_{\lambda^+}= RF(\cl) \setminus RF^\circ(\cl)$.
\end{assertion}

As in the proof of ($\clubsuit$) for regular $\lambda$, the
construction of $\cl$ is accompanied by choosing sets $Z_\alpha$
of size at most $\lambda^+$, of elements that have to be linked to
$B$ in a special way. Let $Z_0=\emptyset$.

Let $\alpha \le \lambda^+$ (for some definitions below we shall
need to refer also to the case $\alpha=\lambda^+$), and assume
that we have defined $R_\beta$ (the rungs of the ladder $\cl$) as
well as $Z_\beta$ for all $\beta <\alpha$. Write $Z_{<\alpha}=
\bigcup_{\beta<\alpha}Z_\beta$.

\begin{definition}\label{hammocks}
Let $u \in Z_{< \alpha} \cap RF^\circ(T_\alpha),~ v \in  Z_{<
\alpha} \cap RF(T_\alpha) \cup \{\infty\}$. A $(u,v,\alpha)$-{\em
hammock} is a set of pairwise internally disjoint
$\cy_\alpha$-s.a.p's from $u$ to $v$. A $(u,v,\lambda^+)$-hammock
is plainly called a $(u,v)$-hammock.
\end{definition}

\begin{definition} \label{maxuptolambda}

Let $\kappa$ be a cardinality. We say that a
$(u,v,\alpha)$-hammock $\ch$ is {\em maximal up to $\kappa$} if
one of the following two (mutually exclusive) possibilities
occurs:
\begin{itemize}
\item $\ch$ is a $(u,v,\alpha)$-hammock which is maximal with respect to inclusion and $|\ch|\le\kappa$,
or:
\item $|\ch| = \kappa$ and there exists a $(u,v,\alpha)$-hammock of
size $\kappa^+$.
\end{itemize}
\end{definition}

For the construction of $Z_\alpha$ we now choose a
$(u,v,\alpha)$-hammock maximal up to $\lambda^+$, for every $u \in
Z_{< \alpha} \cap RF^\circ(T_\alpha)$ and every $v \in Z_{<\alpha}
\cup \{\infty\}$, and put its entire vertex set into $Z_\alpha$.



Clearly, a $(u,v,\alpha)$-hammock that is maximal up to
$\lambda^+$ contains a $(u,v,\alpha)$-hammock that is maximal up
to $\mu$  for every cardinal $\mu < \lambda^+$. Hence, choosing
the elements of $Z_\alpha$ carefully, we can see to it that the
set $Z = Z_{\lambda^+}$ satisfies:

\begin{assertion} \label{maximalhammockzmu} For every
$\alpha < \lambda^+$,
$u \in Z \cap RF(T_\alpha)$, every $v \in (Z \cap RF^\circ(T_\alpha)) \cup \{\infty\}$,
and every $\mu \leq \lambda^+$ there exist a
$(u,v,\alpha)$-hammock maximal
up to $\mu$, whose vertex set is contained in $Z$.
\end{assertion}

By Theorem \ref{safelinking} it is also possible to choose the
elements of $Z_\alpha$ so as to guarantee:

\begin{assertion}
\label{safelinkinginz} For every $\alpha < \lambda^+$ such that
$\Gamma_\alpha$ is unhindered,  and every $v \in T_\alpha \cap Z$,
there exists in $\Gamma_\alpha$ a $v$-$B$-path $P$ such that
$\Gamma_\alpha - P$ is unhindered and $V(P) \subseteq Z$.
\end{assertion}

Yet another condition that can  be taken care of  is:

\begin{assertion}
$$V[\cy \langle Z \rangle] \subseteq Z ~ .$$
\end{assertion}

Choosing the vertices $y_\alpha$ of the ladder $\cl$ as members of
$Z$, we can ensure:
\begin{assertion}\label{zroofed}
$Z \subseteq RF(\cl)$.
\end{assertion}

Assertion \ref{zroofed} will be used to pick objects (like paths
or hammocks) contained in $Z$ within $RF(\cl)$. This will be done
without further explicit reference to the assertion.

The description of the construction of $\cl$ is now complete. We
now show how this construction and the fact that $\Phi=\Phi(\cl)$
is not stationary can be used to prove the linkability of
$\Gamma$. As already mentioned, we choose a closed unbounded set
$\Sigma$ disjoint from $\Phi$.


\begin{definition}\label{popular}
A vertex $u$ is said to be {\em popular} if either $u \in
T_{\lambda^+}$, or there exists a $(u,\infty)$-hammock of
cardinality  $\lambda^+$. The set of popular vertices is denoted
by $POP$.
\end{definition}

\begin{remark}\label{goodhamocks}
By Lemma  \ref{stayingwithinroofedtalpha}, if $u \in
RF(T_\alpha)$, then all $\cy$-alternating paths starting at $u$
are contained in $V^\alpha$, and are thus
$\cy_\alpha$-alternating. Since for each $\alpha < \lambda^+$ we
have $|\cy_\alpha\langle \sim A \rangle| \le \lambda$ and
$|\cy^\infty_\alpha| \le \lambda$, we can assume that all s.a.p's
in the hammock witnessing the popularity of $u$ are, in fact,
$(\cy_\alpha \langle A \rangle)^f$-alternating.
\end{remark}


Let $IE$ be the set of pairs $(u,v)$ of vertices
in $Z$ having a
$(u,v)$-hammock of cardinality at least $\lambda^+$ (``IE" stands
for ``imaginary edges"). Let $SIE$ be the set of all pairs $(u,v)$
for which such a hammock exists in which all s.a.p's are
non-degenerate (see Definition \ref{degenerate}), and let $WIE =
IE \setminus SIE$ (``SIE'' / ``WIE'' stand for ``strong / weak
imaginary edges''). Let $D'$ be the graph $(V,E(D) \cup IE)$. Note
that possibly $E \cap IE \neq \emptyset$, i.e., there may exist
edges that are both ``real" and ``imaginary".


For a warp $\cw$ in $D'$, we define the {\em real part} $Re(\cw)$
of $\cw$ to be the warp in $D$ whose vertex set is $V[\cw]$ and
whose edge set is $E[\cw] \cap E(D)$. If $u=tail(e)$ for an edge
$e \in E[\cw] \cap IE$, we write $\cw_u$ for the warp obtained
from $\cw$ by removing $e$. Also, if $u \in ter[\cw]$ we write
$\cw_u = \cw$.

Let us pause to explain the intuition behind these definitions.
Consider a warp $\cw$ in $D'$ and an imaginary edge $e = (u,v)$ in
it. We should think of  $e$ as a reminder that we should apply
some s.a.p in order to continue the real path ending at $u$ at
some later stage of our construction. Since there are $\lambda^+$
possible such s.a.p's, not all of them will have been destroyed by
the time that it is the turn of $u$ to be linked. Similarly, a
popular vertex $v \in ter[\cw]$ can wait patiently for its turn to
be linked. A vertex $v \in T_{\lambda^+}$ can be linked to $B$ by
applying Assertion \ref{safelinkinginz} for some $\alpha$ which
can be as large as we wish. If there exists a $(v,\infty)$-hammock
of cardinality $\lambda^+$ then, when it is $v$'s turn to be
linked, we can use one of the $(v,\infty)$-s.a.p's to link $v$ to
$T_\alpha$ for some large $\alpha < \lambda^+$.

Let us now return to the rigorous proof.

\begin{definition}\label{goodpair}
Given $\alpha \in \Sigma$, a warp $\cw$ in $D'$ is called an
$\alpha$-{\em linkage blueprint} (or $\alpha$-LB for short) if:

\begin{enumerate}

\item \label{wisroofed} $V[\cw] \subseteq RF_\Gamma(T_\alpha)$.

\item  \label{wlinksa} $  in[\cw  \cup  (\cy \langle
T_\alpha \rangle \setminus \cy \langle V[\cw] \rangle)]  \supseteq
A$.

\item $ V[\cw]  \subseteq Z$.

\item $ |\cw| \leq \lambda$.

\item \label{infinitestrong}
 Every infinite path in $\cw$ contains infinitely many
strong imaginary edges.

 \item \label{salphaorpopular} $ter[\cw] \subseteq POP
\cup T_\alpha$.

\end{enumerate}
\end{definition}

\begin{definition} An $\alpha$-LB $\cw$ satisfying
$ter[\cw] \cap T_\alpha \subseteq  T_{\lambda^+}$  is called a
{\em stable} $\alpha$-LB.
\end{definition}



$\alpha$-linkage blueprints are used to outline a way in which
$\cy$ can be altered, via the application of s.a.p's, so as to
yield an $A$-$T_\alpha$-linkage. An edge  $(u,v)\in E[\cw] \cap
IE$ is going to be replaced by a future application to $\cy$ of a
$(u,v)$-s.a.p. Furthermore, by Definition
\ref{goodpair}(\ref{salphaorpopular}), terminal vertices of $\cw$
 not belonging to $T_{\alpha}$ are popular, again meaning
that they can be linked to $T_\alpha$ by the future use of
s.a.p's.

\begin{assertion}
\label{verygoodtogood} Let $\cv$ be an $\alpha$-LB and let $u \in
ter[Re(\cv)]$. Then there exists an $\alpha$-LB $\cg$ extending
$\cv_u$, such that $Re(\cg)$ links $u$ to $T_\alpha$, and
$ter[Re(\cv)] \subseteq ter[Re(\cg)] \cup \{u\}$.
\end{assertion}

(See Definition \ref{extendingwarps} of a warp being an extension
of another warp. Note that in this case, the extension will not necessarily
be a forward extension.)

\begin{proof}

Let $U=\cv(u)$, namely the path in $\cv$ containing $u$.
Consider first the case that $u \in ter[\cv]$. We may clearly
assume that $u \not \in T_\alpha$, as otherwise we could take $\cg
= \cv$. By
  Definition \ref{goodpair}(\ref{salphaorpopular}), it follows
that  $u \in POP$. Since $u \not\in T_{\lambda^+}$, by Assertion
\ref{maximalhammockzmu}  there exists a $(u,\infty)$-hammock $\ch$
of size $\lambda^+$ contained in $Z$.
Since  $|\cy_\alpha \langle \sim A \rangle | \le \lambda$ and
since by Lemma \ref{most_y_reach_salpha} also $|\cy_\alpha \langle
\sim T_\alpha \rangle| \le \lambda$, it follows that $\ch$
contains a $\cy\langle A,T_\alpha \rangle$-s.a.p $Q$, that does
not meet $V[\cv]$ apart from at $u$. Let $\cj=\cy \triangle Q$.
Then $\cg = \cv \diamond \cj$ is the desired $\alpha$-LB (the
``$\diamond$" operation is defined in Definition \ref{star}).

Assume next that $u \not\in ter[\cv]$. Let $(u,v)$ be the edge in
$E[U]$ having $u$ as its tail. Then $(u,v) \in IE$, meaning that
there exists a $(u,v)$-hammock $\ch$ of size $\lambda^+$,
contained in $Z$. Again, there exists a s.a.p $Q \in \ch$ such
that $V(Q) \setminus \{u\}$ avoids $\cy_\alpha \langle V[\cv]
\rangle \cup  \cy \langle \sim T_\alpha \rangle$ and $in[\cj]
\subseteq A$. Let $\cj=\cy \triangle Q$. If $(u,v) \in SIE$ we can
also assume that $\cj$ links $u$ to $T_\alpha$ and hence $\cv \diamond
\cj$ is the desired warp $\cg$. If $(u,v) \in WIE$, let $\cg_1 =
\cv \diamond \cj$, let $P_1$ be the path in $Re(\cg_1)$ containing
$u$ (thus
 $P_1$ goes through $v$, and then continues along $U$,
 until it reaches either $ter(U)$ or the next imaginary edge on $U$),
and let $u_1 = ter(P_1)$. Apply the same construction, replacing
$u$ by $u_1$, to obtain an $\alpha$-LB $\cg_2$. By part
\ref{infinitestrong} of definition \ref{goodpair} we know that
this process will terminate after a finite number of steps. The
warp $\cg_i$ obtained at that stage is the desired warp $\cg$.
\end{proof}

We shall need to strengthen Assertion \ref{verygoodtogood} in two
ways. One is that we wish to link $u$ to $B$,  not merely to
$T_\alpha$. The other is that we wish $\cg$ to be a {\em stable}
linkage-blueprint. The next assertion takes care of both these
points:

\begin{assertion}
\label{goodtoverygood} If $\cv$ is an $\alpha$-LB and $z \in
T_\alpha \cap ter[\cv]$ then there exist an ordinal $\beta >
\alpha$ and a stable $\beta$-LB $\cu$ extending $\cv$, such that:

\begin{enumerate}
\item
 $Re(\cu)$ links $z$ to $B$.

\item
$ter[Re(\cv)] \subseteq ter[Re(\cu)] \cup T_\alpha$.

\item \label{slpremains}

$ter[\cv] \cap \tlp \subseteq ter[\cu] \cup \{z\}$.

\end{enumerate}

\end{assertion}

\begin{proof}
By Assertion \ref{safelinkinginz} there exists in $\Gamma_\alpha$
a $z$-$B$-path $P$ contained in $Z$, such that $\Gamma_\alpha - P$
is unhindered.

\begin{claim}\label{existsx}
There exist a set $X$ of vertices of size at most $\lambda$, and
an ordinal $\beta > \alpha$, satisfying:

\begin{enumerate}
\item \label{containvandtercv}
$V(P) \cup (ter[\cv] \cap T_\alpha) \subseteq X \subseteq Z \cap
RF(T_\beta)$.
\item \label{xinslp} $X \cap T_\beta \subseteq T_{\lambda^+}$.
\item \label{consistencyofx} $V[\cy \langle X \rangle] \subseteq X$.
\item \label{xcontainsfreepoints}
$V[\cy \langle T_\alpha \rangle \setminus \cy \langle T_\beta
\rangle] \cup V[\cy \langle T_\beta \rangle \setminus \cy \langle
T_\alpha \rangle] \subseteq X$.
\item \label{maximalhammock} For every $u \in X \setminus T_{\lambda^+}$ and $v \in X \cup
\{\infty\}$ there exists a $(u,v)$-hammock maximal up to $\lambda$
contained in $X$.
\end{enumerate}
\end{claim}

The construction of $X$ and $\beta$ is done by a closing-up
process. By Assertion \ref{maximalhammockzmu}, for every $u \in Z
\setminus T_{\lambda^+}$ and $v \in Z \cup \infty$ there exists
 a $(u,v)$-hammock $H_{u,v}$ contained in $Z$ that is maximal up to
$\lambda$. Let $M_{u,v}=V[H_{u,v}]$. For $u \in Z \cap
T_{\lambda^+}$ let $\gamma_u = \min \{\theta :~ u \in T_\theta\}$.
For $u \in Z \setminus T_{\lambda^+}$ define $\gamma_u = \min
\{\theta :~ u \in RF^\circ(T_\theta)\}$. For every $\gamma <
\lambda^+$ let $H_\gamma = V[\cy \langle T_\alpha \rangle
\setminus \cy \langle T_\gamma \rangle] \cup V[\cy \langle
T_\gamma \rangle \setminus \cy \langle T_\alpha \rangle]$

Let $\beta_0 = \alpha$ and let $X_0 = V(P) \cup (ter[\cv] \cap
T_\alpha)$. For every $i < \omega$, let $\beta_{i+1} = \sup
\{\gamma_x : ~ x \in X_i\}$ and let
$$X_{i+1} =
\bigcup_{\begin{array}{cl}
        u \in X_i \setminus T_{\lambda^+}\\
        v \in X_i \cup \{\infty\}
           \end{array}} M_{u,v}
\cup H_{\beta_i} \cup V[\cy \langle X_i \rangle] \; .$$

Taking $X = \bigcup_{i < \omega} X_i$ and $\beta = \sup_i \beta_i$
proves the claim.

\begin{claim}\label{fillingbyimaginary}
Let $Q$ be a $(u,v)$-s.a.p, where
 $u \in Z\setminus \tlp$ and $v \in Z \cup \{\infty\}$. If
$V(Q) \cap X  \subseteq \{u,v\}$  then:
\begin{enumerate}
\item
If $v \in Z$ then  $(u,v) \in IE$.
\item
If $v =\infty$ then $u \in POP$. \end{enumerate}
\end{claim}

To prove (1), assume that $(u,v) \not \in IE$. By the properties
of $X$ there exists a maximal $(u,v)$-hammock $H$ lying within
$X$. By the maximality of $H$, the s.a.p $Q$ must meet some path
belonging to  $H$, contradicting the assumption that $V(Q) \cap X
= \{u,v\}$. The proof of (2) is similar.



Returning to the proof of the assertion, apply now ($\clubsuit$)
to the web $\Gamma^\beta_\alpha -P$, to obtain a
$T_\alpha$-$T_\beta$-linkage $\cw$ containing $P$. Let $\ca=\cv
\cup (\cy \langle T_\alpha \cap X, \sim
V[\cv]\rangle)[RF(T_\alpha)]$
and $\cc=\ca \diamond \cw[X]$.
If $V[\cw \langle X \rangle] \subseteq X$ then
we can take $\cu = \cc$ to be our desired $\beta$-LB.
Unfortunately, there is no way to guarantee
$V[\cw \langle X \rangle] \subseteq X$.
Therefore, there might be paths in $\cw$ with some
vertices in $X$ and some vertices not in $X$.
In this case there may be vertices in $ter[\cw[X]]$
which are not in $ter[\cw]$, and thus $\cc$ might
not be a linkage-blueprint, failing to satisfy
part \ref{salphaorpopular} of Definition
\ref{goodpair}. This is the reason we need to use imaginary edges.
We use imaginary edges to ``mend'' the holes in $\cw[X]$.
This is done according to the behavior of $\cw$ outside of $X$.


Define  $\cz=\cw \downharpoonright X$, namely the fractured warp consisting
of the ``holes" formed in $\cw$ by the removal of $X$ (thus
$E[\cz]=E[\cw] \setminus E[\cw[X]]$). By Theorem \ref{czcysaps} and Remark \ref{fracturedczcy}
there exists an assignment of an element $v=v(u) \in ter[\cz] \cup
\{\infty\}$ and a $(u,v(u))$-$[\cz,\cy]$-s.a.p $Q(u)$ to every $u
\in in[\cz]$, such that $v(u_1) \neq v(u_2)$ whenever $u_1 \neq
u_2$ and $v(u_1),v(u_2) \in ter[\cz]$.

The desired warp $\cu$ is now defined by $ISO(\cu) = ISO(\cv)$ and
$E[\cu]=E[\cc] \cup
\{(u,v(u)) \mid u \in in[\cz],~Q(u) ~~~ {\mbox {is finite}}\}$. By
part (1) of Claim \ref{fillingbyimaginary} for every $u$ such that
$v(u) \in ter[\cz]$ the edge $(u,v(u))$ belongs to $IE$, and thus
$E[\cu] \subseteq E \cup IE$. By part (2) of the claim, every $u
\in in[\cz]$ for which $v(u) =\infty$ is popular, and thus
$ter[\cu] \subseteq POP$. By Lemma \ref{safeisnondegenerate},
whenever $Q(u)$ is finite and degenerate $u$ and $v(u)$ lie on the
same path from $\cw$. Since $\cw$ is f.c., this implies that every
infinite path in $\cu$ contains infinitely many non-degenerate
edges, as required in the definition of linkage-blueprints. Put
together, this shows that $\cu$ is a $\beta$-LB. By Claim
\ref{existsx}(\ref{xinslp}) it is stable.
\end{proof}


\begin{definition}\label{realextension}
For $\alpha \leq \beta < \lambda^+$, we say that a $\beta$-LB $\cu$
is a {\em real extension} of an $\alpha$-LB $\cv$ if $Re(\cu)$ is
an extension of $Re(\cv)$ and
$V[\cv] \subseteq (ter[\cu] \cap ter[\cv]) \cup
tail[E[\cu] \cap E[\cv]] \cup V[Re(\cu) \langle B \rangle]$
We write then $\cv \sqsubseteq \cu$.
\end{definition}

 We shall later ``grow" blueprints $\cv_\alpha$, ordered by the
 ``$\sqsubseteq$" order.
The requirement
$V[\cv] \subseteq (ter[\cu] \cap ter[\cv]) \cup
tail[E[\cu] \cap E[\cv]] \cup V[Re(\cu) \langle B \rangle]$
should be thought of as follows.
Let $R \in Re(\cv)$ (so $ter(R)$ is either
a vertex in $ter[\cv]$ or is the tail of some imaginary edge)
and let $R' \in Re(\cu)$ be the path containing it.
One of the following two happens.
\begin{itemize}
\item $ter(R) \in ter[Re(\cu)]$, so $ter(R) = ter(R')$,
meaning that $R$ was not ``continued forward'',
\item $ter(R) \in V[Re(\cu)  \langle B \rangle]$,
so $ter(R') \in B$, meaning
that $R$ was ``continued all the way to $B$''.
\end{itemize}
The third possibility, that $R$ is continued, but
not all the way to $B$, should be disallowed
in order to avoid infinite paths.

One can easily check that $\sqsubseteq$ is a partial order. The next assertion
 states that it behaves well with respect to taking limits:

\begin{assertion}
\label{verygoodlimit}

Let $\alpha < \lambda^+$ be a limit ordinal
and let $(\beta_\theta \mid \theta \leq \alpha)$ be an ascending sequence
of ordinals satisfying
$\beta_\alpha = \sup_{\theta < \alpha} \beta_\theta < \lambda^+$.
Let  $\cv_\theta$ be
a stable $\beta_\theta$-LB for every $~\theta < \alpha$, where
$\cv_\mu \sqsubseteq \cv_\nu$ whenever
$\mu<\nu<\alpha$.
Let the warp $\cv_\alpha = \lim_{\theta < \alpha} \cv_\alpha$. Namely,
$V[\cv_\alpha] = \bigcup_{\theta < \alpha} V[\cv_\theta]$
and
$E[\cv_\alpha] =
\bigcup_{\beta < \alpha} \bigcap_{\theta \geq \beta} E[\cv_\theta]$
Then $\cv_\alpha$ is a stable $\beta_\alpha$-LB, that is a real extension of all
$\cv_{\beta_\theta},~\theta< \alpha$.
\end{assertion}

Checking most of the properties of an $\alpha$-LB for
$\cv_\alpha$ is easy. The only
non-trivial part is part (\ref{salphaorpopular}) of the
definition, which follows from the stability of the warps
$\cv_\theta$.

We can now combine Assertions \ref{verygoodtogood} and
\ref{goodtoverygood}, to obtain the following:

\begin{assertion}
\label{verygoodtoverygood} Let $\cv$ be a stable $\alpha$-LB and
let $u \in ter[Re(\cv)]$. Then there exist $\beta > \alpha$ and a
stable $\beta$-linkage-blueprint $\cu$, such that:
\begin{enumerate}
\item
$\cv \sqsubseteq \cu$.

\item
$Re(\cu)$ links $u$ to $B$,  and:

\item $ter[Re(\cv)] \subseteq ter[Re(\cu)] \cup \{u\}$.
\end{enumerate}
\end{assertion}

\begin{proof}

By Assertion   \ref{verygoodtogood}, there exists an $\alpha$-LB
$\cg$ extending  $\cv_u$, and satisfying $ter[Re(\cv)] \subseteq
ter[Re(\cg)] \cup \{u\}$. Let $z$ be the terminal vertex of the
path in $Re(\cg)$ containing $u$. Use Assertion
\ref{goodtoverygood} to obtain an ordinal $\beta > \alpha$ and a
stable $\beta$-LB $\cu$ extending $\cg$, such that $Re(\cu)$ links
$z$ to $B$, and $ter[Re(\cg)] \subseteq ter[Re(\cu)] \cup
T_\alpha$. Thus $ter[Re(\cv)] \subseteq ter[Re(\cu)] \cup T_\alpha
\cup \{u\}$.

To show that $ter[Re(\cv)] \subseteq ter[Re(\cu)] \cup \{u\}$ it
suffices to prove that $ter[Re(\cv)] \cap T_\alpha \subseteq
ter[Re(\cu)] \cup \{u\}$. Note that $ter[Re(\cv)] \cap T_\alpha
\subseteq ter[\cv] \cap T_\alpha$. Since $\cv$ is a stable
$\alpha$-LB, we have $ter[\cv] \cap T_\alpha \subseteq \tlp$.
By part
(\ref{slpremains}) of Assertion \ref{goodtoverygood}, we have
$ter[Re(\cv)] \cap T_\alpha \subseteq ter[Re(\cu)] \cup \{u\}$.
One can easily check the $\cu$ is a real extension of $\cv$,
proving the assertion.

\end{proof}



We can now conclude the proof of ($\clubsuit\clubsuit$). We shall
do this by applying  Assertion \ref{verygoodtoverygood} $\lambda$
times. Observe first that $\langle A' \rangle$ is a $0$-LB. By
Assertion \ref{goodtoverygood},  it can be extended to a stable
$\sigma_0$-LB $\cv_0$, for some $0 < \sigma_0 < \lambda^+$. Choose
now some $u_0 \in ter[Re(\cv_0)]$. By Assertion
\ref{verygoodtoverygood}, there exists a stable $\sigma_1$-LB
$\cv_1$ for some $\sigma_1 > \sigma_0$, such that $\cv_0
\sqsubseteq \cv_1$ and $Re(\cv_1)$ links $u_0$ to $B$. We continue
this way. For each $\alpha < \lambda$ we choose $u_\alpha \in
ter[Re(\cv_\alpha)]$ and use Assertion \ref{verygoodtoverygood} to
find a stable $\sigma_{\alpha+1}$-LB such that $\cv_{\alpha}
\sqsubseteq \cv_{\alpha+1}$ and $Re(\cv_{\alpha+1})$ links
$u_\alpha$ to $B$. For limit ordinals $\alpha \leq \lambda$ define
$\sigma_\alpha = \sup_{\theta <\alpha} \sigma_\theta$ and define
$\cv_\alpha$
as in Assertion \ref{verygoodlimit}, so $\cv_\alpha$  is a stable
$\sigma_\alpha$-LB.

Choosing  the vertices $u_\alpha$ appropriately, we can procure
the following condition:

$$\{u_\alpha : ~ \alpha < \lambda\} = \bigcup_{\alpha < \lambda}
ter[Re(\cv_\alpha)] \setminus B \; .$$

This implies that $\cv_\lambda = Re(\cv_\lambda)$ and
$ter[\cv_\lambda] \subseteq B$. Let $\ch$ be the warp obtained by
adding to $\cv_\lambda$ all paths of $\cy$ not intersecting
$V[\cv_\lambda]$ and let $\sigma=\sigma_\lambda$. Then $\ch$ is an
$A$-$T_{\sigma}$-linkage linking $A'$ to $B$. Since $\Gamma\quo
T_\sigma$ is unhindered, $\ch$ is a half-way linkage, as required
in the theorem. $\enp$

\section{Open problems in infinite matching theory}
The Erd\H{o}s-Menger conjecture pointed at the way duality should
be formulated in the infinite case: rather than state equality of
cardinalities, the conjecture stated the existence of dual objects
satisfying the so-called ``complementary slackness conditions".
There are still many problems of this type that are open. One of
the most attractive of those is the ``fish-scale conjecture",
named so because of the way its objects can be drawn
\cite{aharonikorman}:

\begin{conjecture}\label{fishscale}
In every poset  not containing an infinite antichain there exist a
chain $C$ and a decomposition of the vertex set into antichains
$A_i$, such that $C$ meets every antichain $A_i$.
\end{conjecture}

The dual statement, obtained by replacing the terms
``chain" and ``antichain", follows from the
infinite version of K\"{o}nig's theorem
\cite{steffensinfinitedilworth, aharonisurvey}. It is likely that,
if true, Conjecture \ref{fishscale} does not have much to do with
posets, but with a very general property of infinite hypergraphs.

\begin{definition} Let $H=(V,E)$ be a hypergraph. A {\em matching} in
$H$ is a subset of $E$ consisting of disjoint edges. An {\em edge
cover} is a subset of $E$ whose union is $V$. A matching $I$ is
called {\em strongly maximal} if $|J \setminus I| \le |I \setminus
J|$ for every matching $J$ in $H$. An edge cover $F$ is called
{\em strongly minimal} if $|K \setminus F| \ge |F \setminus K|$
for every edge cover $K$ in $H$.
\end{definition}

As noted above, our main theorem is tantamount to the fact that
the hypergraph of vertex sets of $A$--$B$-paths in a web possesses
a strongly maximal matching. Call a hypergraph {\em finitely
bounded} if its edges are of size bounded by some fixed finite
number. Call a hypergraph $H$ a {\em flag complex} if it is closed
down, namely every subset of an edge is also an edge, and it is
$2$-determined, namely if all $2$-subsets of a set belong to $H$
then the set belongs to $H$.

\begin{conjecture}\label{stronglymax}
$\;$

\begin {enumerate}

\item Every finitely bounded hypergraph contains a strongly maximal
matching and a strongly minimal cover.
\item
Any flag complex contains a strongly minimal cover.

\end{enumerate}
\end{conjecture}

 Conjecture
 \ref{fishscale} would follow by a compactness argument from part (2) of this
 conjecture. For graphs part (1)  of the conjecture follows
from the main theorem of \cite{inftutte}.

The mere condition of having only finite edges does not suffice
for the existence of a strongly maximal matching, as was shown in
\cite{ahlswedekhach}. In the example given there,  for every
matching $M$ there exists a matching $M'$ with $|M \setminus
M'|=2,~ |M' \setminus M|=3$.

\begin{problem}[Tardos]
Is it true that in every hypergraph with finite edges there exists
a matching $M$ such that  no matching $M'$ exists for which $|M
\setminus M'|=1,~ |M' \setminus M|=2$?
\end{problem}

{\bf Acknowledgement} We are grateful to the members of the Hamburg
University Combinatorics seminar led by   Reinhard Diestel, for a
careful reading of a preliminary draft of this paper, and for
pointing out many inaccuracies. In particular, Henning Bruhn and
Maya Stein contributed to the presentation and correctness of
Sections 8 and 9, and Philipp Spr\"ussel made helpful remarks
concerning Section 2. We are also grateful to Eliahu Levy and Mat\'e
Vizer for many useful comments.


\end{document}